\documentstyle[11pt]{article}

\begin{document}
\noindent
\begin{center}
  {\LARGE Symplectic surgery and Gromov-Witten invariants of Calabi-Yau
3-folds I}
  \end{center}

  \noindent
  \begin{center}
   {\large An-Min Li}\footnote{partially supported by a NSFC grant and a
   Qiu Shi grant}\\[5pt]
      Department of Mathematics, Sichuan University\\
        Chengdu, PRC\\[5pt]
    {\large Yongbin Ruan}\footnote{partially supported by a NSF grant and a
    Sloan fellowship}\\[5pt]
      Department of Mathematics, University of Wisconsin-Madison\\
      Madison, WI 53706\\[5pt]
\end{center}

              \def \J{{\cal J}}
              \def \Map{Map(S^2, V)}
              \def \M{{\cal M}}
              \def \A{{\cal A}}
              \def \B{{\cal B}}
              \def \C{{\bf C}}
              \def \Z{{\bf Z}}
              \def \R{{\bf R}}
              \def \P{{\bf P}}
              \def \I{{\bf I}}
              \def \N{{\bf N}}
              \def \T{{\bf T}}
              \def \O{{\cal O}}
              \def \Q{{\bf Q}}
              \def \D{{\bf D}}
              \def \H{{\bf H}}
              \def \S{{\cal S}}
              \def \e{{\bf E}}
              \def \k{{\bf k}}
              \def \U{{\cal U}}
              \def \E{{\cal E}}
              \def \F{{\cal F}}
              \def \L{{\cal L}}
              \def \K{{\cal K}}
              \def \G{{\bf G}}

\tableofcontents
\section{Introduction}
    This is the first of a series of papers devoted to the study of how
Gromov-Witten invariants of 3-folds transform under surgery.
Special attention will be paid to Calabi-Yau 3-folds.
 In \cite{R1}, \cite{RT1},\cite{RT2}, Ruan-Tian established
the mathematical foundations of the
theory of quantum cohomology or Gromov-Witten invariants for
the semipositive symplectic manifolds.
Recently, semipositivity condition has been removed
by the work of many authors: \cite{LT2},  \cite{B}, \cite{FO},\cite{LT3}, \cite{R5}, \cite{S}.
The focus now is on  calculation and applications.
As quantum cohomology was developed, many examples were
calculated by direct computation. In particular, the localization technique has had many
successes
\cite{Gi}, \cite{LLY}.
In this article, we
take a new direction by computing the change of Gromov-Witten
invariants under some well-known surgeries in algebraic geometry.
Then, our formula will calculate the Gromov-Witten invariants of
any 3-fold obtained by performing these surgeries from known examples.
Another outcome of this paper is the appearance of the relation between
the
naturality problem of quantum cohomology and birational geometry.
This newly discovered relation is both surprising and intriguing. We
hope to investigate this more in the future.

To motivate our choice of surgeries, let's review the
classification of Calabi-Yau (CY) 3-folds. A first step  is to classify CY 3-folds in the
same birational class.
Any two CY 3-folds are
birationally equivalent to each other iff they can be connected by a sequence of
{\em flops} \cite{Ka}, \cite{K}.
A flop is a kind of surgery: two  CY 3-folds
$M$ and $M_f$ are said to be related
by a flop if there is a a singular CY 3-fold $M_s$
obtained from  $M$ by a ``small contraction''
(contracting finitely many rational curves)
with $M_f$ obtained from $M_s$ by blowing up these singularities differently.
Flopping can be defined for any 3-fold once the corresponding small contraction exists. A conjecture by Morrison \cite{Mo1} is that
flopping among CY 3-folds induces an isomorphism on quantum cohomology. We will give an affirmative
answer to Morrison's conjecture for any 3-fold. It is well-known that
flopping
does not induce an isomorphism on ordinary cohomology. Birational
geometry is a central  topic in algebraic geometry. It is often difficult
to construct birational invariants. A proof of Morrison's conjecture would provide
the first truly quantum birational invariant.
It seems to us that it will be an important
problem to study quantum cohomology under other type of
birational transformations, for example flips. This requires an extension of quantum
cohomology to orbifolds.
We shall leave this to  future research.

For nonbirationally equivalent CY  3-folds, an influential conjecture by
Miles Reid (Reid's fantasy) states that any two Calabi-Yau 3-folds are
connected to each other by a sequence of so-called contract-deform or
deform-resolve surgery.
Contract-deform  means that we contract certain curves or divisors to obtain
a singular Calabi-Yau 3-fold and deform the result to a smooth CY 3-fold.
Deform-resolve is the  reverse.
Reid's original conjecture
was stated in the category of non-K\"ahler manifolds.
To stay in the K\"ahler category,
the contractions must be {\em extremal} in the sense of Mori.
The smoothing theory
for the case of ordinary double points was treated in \cite{F}, \cite{Ti1}.
In the general case, we refer to \cite{Gro} for  references.
An  extremal contraction-smoothing or its opposite
surgery is called {\em an extremal transition or transition}. We call
a transition {\em small} if the corresponding contraction is
small. Furthermore, a small transition can be performed over any
3-fold.
A modified version of Reid's conjecture
is that any two smooth CY 3-folds are connected to each other by a sequence
of flops or extremal transitions. Large classes of
CY 3-folds are indeed connected to each other in this way.

The classification of CY 3-folds involves the study of surgeries.
It would be desirable to study the effect of
these surgeries on mirror symmetry. This paper takes a  step in this direction.
Recall that the mirror
symmetry conjecture asserts that every CY 3-fold $X$ has a mirror
partner $Y$ such that $h_{1,1}(X)=h_{1,2}(Y),
h_{1,2}(X)=h_{1,1}(Y)$. Furthermore, the quantum cohomology of $X$
is the same as the Yukawa coupling of $Y$. It was known that the most
general form of the mirror symmetry conjecture is false due to the
existence of rigid CY 3-folds with $h_{1,2}=0$. The conjectured
mirror partner of $X$ will have $h_{1,1}=0$. Hence, it can not be
K\"ahler. We believe that the most difficult problem in mirror
symmetry is to find the precise category of CY 3-folds where
mirror symmetry holds.
An interesting speculation \cite{Mo2} here was that,
except for the obvious counterexamples,
each extremal transition has a mirror surgery which preserves mirror symmetry.
It can be summarized as follows: {\em Mirror surgery conjecture: Every extremal
transition $L$ has a mirror surgery $L_m$ with the following property: If we have
a mirror pair $(X, Y)$ and
perform an extremal transition $L$ on $X$ and obtain $\tilde{X}$, then one of the following
is true: (1) $\tilde{X}$ has no large complex structure limit. In this case, $X$ has no
mirror. (2) The mirror of $\tilde{X}$ can be obtained by performing $L_m$ on $Y$.}
Indeed, all the known examples of rigid CY 3-folds can be obtained by extremal
transitions.
 This gives a nice explanation of the failure of
mirror symmetry for rigid CY 3-folds. This conjecture grew from a discussion of the second
author with P. Aspinwall.
A closely related conjecture was also proposed by D. Morrison \cite{Mo2}.
Once we understand how relevant
invariants change under extremal transitions, we can extend mirror symmetry
to a large classes of CY 3-folds, and hopefully find the precise category
where mirror symmetry holds. One can view the mirror surgery conjecture as a combination of
the classification problem and mirror symmetry. Any results on the mirror surgery conjecture would yield a deeper
understanding to both the classification problem and mirror symmetry.
If we want to prove Morrison's conjecture or  the mirror surgery
conjecture,
it is clear that we have to calculate the change of GW-invariants under
flops and transitions. In this paper, we study flops and small transitions over any 3-fold.
In  subsequent papers, we will focus on CY 3-folds to deal with other type of transitions where
there is a good classification theory.

 The invariants we consider here are primitive GW-invariants $\Psi^M_{(A,g,m)}(
\overline{\M}_{g,m}; \{\alpha_i\})$ for the stable range $2g+m\geq 3$.  It was conjectured \cite{RT2} that any GW-invariants
can be reduced to primitive ones. We shall drop $\overline{\M}_{g,m}$
to simplify the notation. If $m=0$, We will  drop $m$ as well.
For primitive invariants, it is also convenient to drop the stable
range condition. The construction is standard and we leave it to
the readers. Then, one can eliminate the divisor class $\alpha\in H^2(M, \R)$ by the relation

$$
\Psi^M_{(A,g,m+1)}(\alpha,\alpha_1, \cdots, \alpha_m)=
   \alpha(A)\Psi^M_{(A,g,k)}(\alpha_1, \cdots, \alpha_m),
\leqno(1.1)
$$
for $A\neq 0$.

Choose a basis $A_1, \cdots, A_m$ of $H_2(M, \Q)$.
For $A=\sum_i a_i A_i$,
we define the formal product
$q^A=(q_{A_1})^{a_1}\cdots (q_{A_m})^{a_m}$.
We can define a quantum 3-point function
$$
\Psi^M_w(\alpha_1, \alpha_2, \alpha_3)=
      \sum_j\frac{1}{j!}\sum_{A} \Psi^M_{(A,0,j+3)}(\alpha_1, \alpha_2, \alpha_3,w, \cdots, w)q^A,
\leqno(1.2)
$$
where $w$ appears $m$ times. Here, we view $\Psi^M$ as a power series in the formal variables $p_i=q^{A_i}$.
Clearly,
an isomorphism on $H_2$ will induce a change of variables $p_i$.
One can define the quantum product from the quantum 3-point function. They
contain the same information. For our purpose, it is
convenient to work directly with the quantum 3-point function.
\vskip 0.1in
\noindent
{\bf Definition 1.1: }{\it Suppose that
$$\varphi: H_2(X, \Z)\rightarrow H_2(Y, \Z),\ H^{even}(Y, \R)\rightarrow H^{even}(X, \R)$$
are vector space  homomorphisms such that the maps on $H_2, H^2$
are dual to each other. We say $\varphi$ is  natural with respect to (big) quantum cohomology if $\varphi^* \Psi^X_0=\Psi^Y_0$
( $\varphi^* \Psi^X_{\varphi^*w}=\Psi^Y_w$) up to a change of formal variable $q^A\rightarrow q^{\varphi(A)}$.
If $\varphi$ is also an isomorphism, we say $\varphi$
induces an isomorphism on (big) quantum cohomology or they have the same (big) quantum cohomology.}
\vskip 0.1in

Here, two power series $F,G$ are treated as the same if $F=H+F', G=H+G'$
such that $G'$ is an analytic continuation of $F'$.
For example, we can expand $\frac{1}{1-t}=\sum_{i=0}t^i$ at $t=0$
or $\frac{1}{1-t}=\frac{1}{-t(1-t^{-1})}=-\sum_{i=0} t^{-i-1}$
at $t=\infty$. Hence, we will treat $\sum_{i=0}t^i, \sum_{i=0} t^{-i-1}$
as the same power series.

When $X, Y$ are 3-folds, such a $\varphi$ is completely determined  by  maps on $H_2$.
For example, the dual map of $\varphi: H_2(X, \Z) \rightarrow H_2(Y, \Z)$ gives
a map $H^2(Y, \R)\rightarrow H^2(X, \R)$. A map $H^4(Y, \R)\rightarrow H^4(X, \R)$
is Poincar\'e dual to a map $H_2(Y, \R)\rightarrow H_2(X, \R)$. In the case of
flops, the natural map $H_2(X, \Z) \rightarrow H_2(Y, \Z)$ is an isomorphism.
Therefore, we can take the map $H_2(Y, \Z)\rightarrow H_2(X, \Z)$ as its inverse.
The maps on $H^0, H^6$ are obvious.

 Suppose that $M_f$ is obtained after a flop on $M$.
There is a natural isomorphism (see section 2)
$$
\varphi: H_{2}(M, \Z)\rightarrow H_{2}(M_f, \Z).\leqno(1.3)
$$
The manifolds $M$ and $M_f$ have the same set of exceptional curves. Suppose
that $\Gamma$ is an exceptional curve and  $\Gamma_f$ is the
corresponding exceptional curve on $M_f$. Then,
$$
\varphi([\Gamma])=-[\Gamma_f].\leqno(1.4)
$$
Our first theorem is that

\vskip 0.1in

\noindent {\bf Theorem A }{\it
If $A\neq n [\Gamma]$ for any exceptional curve  $\Gamma$,
then
$$
\Psi^M_{(A,g,m)}(\{\varphi^*\alpha_i\})=\Psi^{M_f}_{(\varphi(A),g,m)}(\{\alpha_i\}).\leqno(1.5)
$$
Moreover,
$$
\Psi^M_{(n[\Gamma],g)}=\Psi^{M_f}_{(n[\Gamma_f],g)}.\leqno(1.6)
$$
}

\vskip 0.1in

  When $M$ is a Calabi-Yau 3-fold, Theorem A takes a particularly
  simple form.
\vskip 0.1in
\noindent
{\bf Corollary A.1: }{\it Suppose that $M$ is a Calabi-Yau 3-fold.
If $A\neq n [\Gamma]$ for any exceptional curve  $\Gamma$,
then
$$
\Psi^M_{(A,g)}=\Psi^{M_f}_{(\varphi(A),g)}.
$$
Moreover,
$$
\Psi^M_{(n[\Gamma],g)}=\Psi^{M_f}_{(n[\Gamma_f],g)}.\leqno(1.7)
$$
}

\vskip 0.1in

  Using our formula in genus zero case and Morrison's argument \cite{Mo1},
  we have the following corollary.

\vskip 0.1in

\noindent{\bf Corollary A.2: }{\it
$
\varphi
$
induces an isomorphism on quantum cohomology.}

\vskip 0.1in

Recall that a minimal model is a projective variety with terminal
singularities and nef canonical bundle. It is known that in
higher dimension there are many different minimal models in the
same birational class. However, in dimension three, they are
related by flops. Then the above corollary yields

\vskip 0.1in
\noindent
{\bf Corollary A.3: }{\it Any  two smooth
minimal models in dimension three have the same quantum
cohomology.}
\vskip 0.1in

The second author conjectures that {\em two smooth minimal models over any dimension
have isomorphic quantum cohomology \cite{R6}.}
When $M$ is a CY-3-fold, the above corollary implies  Morrison's
conjecture. In the case of a CY 3-fold, Corollary A.2 admits another
interpretation. Instead of considering $M, M_f$ as different
manifolds, we can use $\phi$ to map the cohomology of different
birational models into a single space called the movable cone. Then,
the Corollary A.2 can be restated as follows:  the quantum 3-point function
extends analytically over the movable cone. In fact, this was how
Morrison stated his conjecture.

      Let $M_e$ be the 3-fold obtained after a small transition.
The small transition induces a surjective map
$$
\varphi_e: H_2(M, \Z)\rightarrow H_2(M_e, \Z).\leqno(1.8)
$$
We choose a right inverse, which is a map $H_2(M_e, \R)\rightarrow H_2( M, \R)$.
Using the general gluing theory we establish, we obtain the following formula:

\vskip 0.1in

\noindent
{\bf Theorem B }{\it Let $M_e$ be obtained by a small transition.
For $B\neq 0$,
}
$$
\Psi^{M_e}_{(B,g,m)}(\{\alpha_i\})=\sum_{\varphi_e(A)=B} \Psi^M_{(A, g,m)}(\{\varphi^*(\alpha_i)\}).
\leqno(1.9)
$$

\vskip 0.1in
Again, Theorem B takes a simple form for Calabi-Yau 3-fold.
\vskip 0.1in

\noindent {\bf Corollary B.1 }{\it Let $M_e$ be obtained by a small transition performed on a Calabi-Yau 3-fold $M$.
Then,for $B\neq 0$}
$$
\Psi^{M_e}_{(B,g)}=\sum_{\varphi_e(A)=B} \Psi^M_{(A, g)}.
$$

\vskip 0.1in
    Then, we can use the genus zero case of Theorem B and
generalize Tian's argument \cite{Ti} to compare quantum cohomology.

 \vskip 0.1in

\noindent{\bf Corollary B.2: }{\it
$
\varphi$ is  natural with respect to big quantum cohomology.}

\vskip 0.1in

    We remark that, in the case of Calabi-Yau 3-folds, the genus zero cases of both
Theorems~A and~B have been studied in physics.
However, it is rather surprising to us that it exhibits such a
naturality relation. From known computation, this is the first
case with this property.
Higher genus GW-invariants are not enumerative invariants in general.
The nature of these invariants is quite mysterious even today.
We were surprised that
the same formula is also true for higher genus GW-invariants.
Several years ago, there was
a mysterious Kodaira-Spencer quantum field theory in physics which dealt with the higher
genus case \cite{BCOV}. We hope
that our calculation will shed some light on the structure of higher genus
invariants.

One of the main ideas of this paper is a symplectic interpretation of
the above algebraic geometric surgeries. Then, we can use symplectic
techniques, which are more flexible. We know of two instances where surgeries were used to calculate GW-invariants.
In  an approach quite different from ours, G. Tian studied the degeneration
of rational curves under  symplectic  degeneration \cite{Ti}
which is an analogue of  the degeneration in algebraic  geometry.
The method McDuff used in \cite{M2}, \cite{Lo}, is similar to that of this paper.
However, she studied the holomorphic curves  completely inside the symplectic
manifold with contact boundary, which is much easier to deal with. The main
difficulty in this paper is to handle  holomorphic curves which intersect
the  boundary.

When we search for the formula for general surgery, Theorems~A
and~B are rather misleading. The natural invariants appearing in
general surgery are {\em log GW-invariants} instead of {\em
absolute GW-invariants}. To describe log GW-invariants, let's
first review the definition of symplectic cutting \cite{L}.

Suppose that $H: M^0\rightarrow \R$ is a proper Hamiltonian function such that the
Hamiltonian vector field $X_H$ generates a circle action, where $M^0\subset M$ is an open domain. By adding a constant,
we can assume that $0$ is a regular value. Then $H^{-1}(0)$ is a smooth
submanifold preserved by circle action. The quotient $Z=H^{-1}(0)/S^1$ is the famous
symplectic reduction. Namely, it has an induced symplectic structure.
We can cut $M$ along $H^{-1}(0)$. Suppose that we obtain two
disjoint components $M^{\pm}$ which have the boundary $H^{-1}(0)$.
We can collapse the $S^1$-action on $H^{-1}(0)$ to obtain $\overline{M}^{\pm}$
containing a real codimension two submanifold $Z$. Moreover, there
is a map
$$\pi: M\rightarrow \overline{M}^+\cup_Z \overline{M}^-,\leqno(1.10)$$
where $\overline{M}^+\cup_Z \overline{M}^-$ is the union of $\overline{M}^{\pm}$
along $Z$. It induces a homomorphism
$$\pi^*: H^*(\overline{M}^+\cup_Z \overline{M}^-, \Q)\rightarrow H^*(M, \Q).\leqno(1.11)$$
It was shown by Lerman \cite{L} that the restriction of the
symplectic structure $\omega$ on $M^{\pm}$ extends to
symplectic structures $\omega^{\pm}$ over $\overline{M}^{\pm}$
such that $\omega^+|_Z=\omega^-|_Z$ is the induced symplectic structure from symplectic reduction.
By the Mayer-Vietoris sequence, a pair of cohomology classes $(\alpha^+, \alpha^-)\in (H^*(\overline{M}^+, \Z),
H^*(\overline{M}^-, \Z))$ with $\alpha^+|_Z=\alpha^-|_Z$
defines a cohomology class  of $\overline{M}^+\cup_Z
\overline{M}^-$, denoted by $\alpha^+\cup_Z \alpha^-$. It is
clear from Lerman's construction that $\pi^*(\omega^+\cup_Z
\omega^-)=\omega$.

Suppose that $B$ is a real codimension two symplectic submanifold of $M$.
By a result of Guillemin and Sternberg \cite{GS},
the symplectic structure of a tubular
neighborhood of $B$ is modeled on  a neighborhood of $Z$ in either
$\overline{M}^+$ or $\overline{M}^-$.

  We can define {\em a log GW-invariant $\Psi^{(M,Z)}$ }
by counting the number of log stable holomorphic maps intersecting
$Z$ at  finitely many points with  prescribed tangency. Let
$T_m=(t_1, \cdots, t_m)$ be a set of positive integers such that
$\sum_i t_i=Z^*(A)$, where $Z^*$ is the Poincare dual of $Z$. We
order them such that $t_1=\cdots=t_l=0$ and $t_i>0$ for $i>l$.
Consider the moduli space $\M_A(g,T_m)$ of genus $g$
pseudo-holomorphic maps $f$ such that $f$ has marked points $(x_1,
\cdots, x_m)$ with the property that $f$ is tangent to $Z$ at
$x_i$ with order $t_i$. Here, $t_i=0$ means that there is no
intersection. Then, we compactify $\M_A(g,T_m)$ by
$\overline{\M}_A(g,T_m)$, {\em the space of log stable maps}. We
have evaluation map
    $$ e_i: \overline{\M}_A(g,T_m)\rightarrow M$$
    for $i\leq l$ and
    $$ e_j: \overline{\M}_A(g,T_m)\rightarrow Z \leqno(1.12) $$
    for $j>l$. Roughly, the log GW-invariants are
defined as
    $$ \Psi^{(M,Z)}_{(A,g, T_m)}(\alpha_1, \cdots,
\alpha_l; \beta_{l+1}, \cdots, \beta_m)=
\int^{vir}_{\overline{\M}_A(g,T_m)}\prod_i e^*_i\alpha_i\wedge
\prod_j e^*_j\beta_j. \leqno(1.13) $$
To carry out the virtual integration, we need to
construct a virtual neighborhood (see Section 5) of
$\overline{\M}_A(g, T_m)$. Then we take the ordinary integrand (1.13) over
the virtual neighborhood.

To justify our notation, we remark that when $T_m=(0,...,0,1,\cdots, 1)$,
$\Psi^{(M,Z)}$ is different from an ordinary or absolute
GW-invariant in general. So $\Psi^{(M, Z)}$ is not a generalized
GW-invariant. \vskip 0.1in \noindent {\bf Theorem C(Theorem 5.3):
} \vskip 0.1in \noindent {\it (i).
$\Psi^{(M,Z)}_{(A,g,T_m)}(\alpha_1,..., \alpha_l; \beta_{l+1},...,
\beta_{m})$ is well-defined, multilinear and skew symmetric with respect to
$\alpha_i$ and $\beta_j$ respectively.
\vskip 0.1in \noindent (ii).
$\Psi^{(M,Z)}_{(A,g,T_m)}(\alpha_1,...,\alpha_l; \beta_{l+1},...,
\beta_{m})$ is independent of the choice of forms $\alpha_i,
\beta_j$ representing the cohomology classes $[\beta_j],
[\alpha_i]$,  and the choice of virtual neighborhoods. \vskip
0.1in \noindent (iii). $\Psi^{(M,Z)}_{(A,g,T_m)}(\alpha_1,\ldots,
\alpha_l; \beta_{l+1},\ldots, \beta_{m})$ is independent of the
choice of almost complex structures on $Z$ and on the normal
bundle $\nu_Z$, and hence is an invariant of $(M,Z)$.} \vskip
0.1in \noindent

{\bf Theorem D(Theorem 5.6, 5.7): }{\it Suppose that
$\alpha^+_i|_Z=\alpha^-_i|_Z$ and hence $\alpha^+_i\cup_Z
\alpha^-_i\in H^*( \overline{M}^+ \cup_Z \overline{M}^-, \R)$. Let
$\alpha_i=\pi^*(\alpha^+_i\cup_Z \alpha^-_i)$ (1.11). There is a
gluing formula to relate $\Psi^M(\{\alpha_i\})$ to the log
invariants $\Psi^{(\overline{M}^+, Z)}(\{\alpha^+_i\})$ and
$\Psi^{(\overline{M}^-, Z)}(\{\alpha^-_i\})$.}

Although it is not needed in this paper, it is also important to
consider descendant GW-invariant in general quantum cohomology
theory. Our gluing theory extends to this more general setting with
the same proof. Let's give a brief description on the formulation.
Recall that the cotangent space of
    each marked point $x_i$ generates an orbifold line bundle $\L_i$ over
the moduli space of stable map $\overline{\M}_A(g,k)$. The
descendant GW-invariant is defined as
$$\Psi^M_{(A,g,k)}(\O_{h_1}(\alpha_1), \cdots, \O_{h_k}(\alpha_k))
=\int^{vir}_{\overline{\M}_A(g,k )}\prod_i
C^{h_i}_1(\L_i)e^*_i\alpha_i. \leqno(1.14)$$ One can define
descendant log-GW-invariants
$\Psi^{M,Z}_{(A,g,T_k)}(\O_{h_1}(\alpha_1),\cdots,
\O_{h_l}(\alpha_l); \O_{h_{l+1}}(\beta_{l+1}), \cdots,
\O_{h_m}(\beta_m))$ in the same fashion and Theorem C extends to
this more general descendant log-GW-invariants. In gluing theory,
we only need $ \Psi^{X,Z}_{(A,g,T_k)}(\O_{h_1}(\alpha_1),\cdots,
\O_{h_l}(\alpha_l); \beta_{l+1}, \cdots, \beta_m)$. Theorem D can
be easily generalized (by the identical proof) to express
descendant invariant $\Psi^M_{(A,g,k)}=(\{\O_{h_i}(\alpha_i)\})$
in terms of descendant log-invariant $\Psi^{(\overline{M}^{\pm},
Z)}(\{\O^{\pm}_{h_i}( \alpha_i)\})$ by the same formula as the
case of ordinary GW-invariants.

     We also remark that different $\alpha^{\pm}_i$ may give
the same $\alpha_i$. Then we obtain different gluing formulas for
the same invariant. For example, if $\alpha_i$ is Poincar\'e dual
to a point, its support could be in $\overline{M}^+$ or
$\overline{M}^-$. This is very important in the applications.
\vskip 0.1in

The exact statement of the general formula is  rather complicated.
We refer the reader to section 5 for the details.

   In a subsequent paper \cite{LQR},
we will apply our general gluing formula to study
the change of GW-invariants for other types of transition. The other types of transitions are closely
related to ordinary blow-up. It is known that the blow-up formula of quantum cohomology is complicated. Hence,
one can not expect to have such a simple formula as in the case of small transition. However, our general
gluing formula gives a complete answer to the case of the other type of transition for Calabi-Yau 3-folds.
We refer to our next paper \cite{LQR} for the detail.

    Let's briefly describe the idea of the proof. The first step is to
reinterpret the flop and extremal transiton as a combination of
contact surgery and symplectic cutting (section 2). Then, we  stretch
symplectic manifolds along  a hypersurface admitting a
local $S^1$-Hamiltonian action. Then we study the behavior of pseudoholomorphic curves under
the stretching. In the case of a contact manifold, every contact manifold with a fixed contact form possesses
a unique vector field called the Reeb vector field.
Hofer observed that the boundary of a finite energy pseudoholomorphic curve
will converge to a periodic orbit of the Reeb vector field.
Furthermore, Hofer and his collabrators established the
analysis of the moduli space of pseudoholomorphic curves whose ends converge to
periodic orbits. We will generalize it to the case of symplectic
cutting. We will do so by casting it into the language of stable maps (section 3).
Furthermore,
we establish  the gluing theorem which is the reverse of the stretching
construction (section 4).
The last piece of information we need is a vanishing theorem for certain
relative
GW-invariants. This is done by a simple index calculation.

   We should point out that the log  GW-invariants
also appeared in the work of Caporaso-Harris \cite{CH}.
They are  closely related to the blow-up formula for Seiberg-Witten
invariants.

    This paper has been revised many times due to referee's suggestions. In many
ways, current version is much better written than original
version. We would like to thank referee for the insistence which
made this possible. But our underline approach has always been the
same. Namely, we explored the local hamiltonian circle action near
the symplectic submanifold and established a gluing theorem in
this setting. In the original version of our paper, in addition to
the results of this paper, we also proved a compactness theorem
for contact manifolds whose periodic orbits are of Bott type. The
contact case is irrelevant to our paper. It was included in our
original version only because it can be proved by the same
technique. Latter, Hofer informed us that it has been dealt with
already in \cite{HWZ1}. Then, the contact case was taken out in
the revised version. There are many papers
  in physics that discuss both flop and extremal transitions . We refer to
\cite{Mo2} for the relevant references. But they  discuss only
genus zero invariants. The paper by P. Wilson \cite{Wi3}
calculates GW-invariants of extremal rays. His  paper is
complementary to ours. During the preparation of this paper, we
received an article \cite{BCKS} which is related to our results.
We are working on the setting of the local Hamiltonian
$S^1$-action, which is not contact in general. In the setting of
contact manifold, a general gluing formula over a contact boundary
requires a definition of contact Floer homology, which is
developed by Eliashberg and Hofer. We are also informed that Ionel
and Parker are developing a gluing theory independently using a
different approach. Apparently, they only used one component of
the moduli space of log-stable map to define their relative
invariants instead of using full moduli space as we did. Then,
their relative invariants is different from ours, which results in
different gluing formula.

   This paper is organized as follows. We will review the constructions of
various surgery operations in symplectic geometry and interpret both the flop and
extremal transitions as symplectic  surgeries. The gluing theory will be
established in
sections~3-5. Theorems~A and~B will be proved
in  section~6.

  The main results of this paper were announced in a conference in Kyoto in
December, 1996. The second author would like to thank S. Mori for
the invitation. The second author also wish to thank H. Clemens,
Y. Eliashberg, M. Gross, J. Koll\'ar, K. Kawamata, S. Katz, E.
Lerman, S. Mori, K. Oguiso, M. Reid, Z. Qin and P. Wilson for
valuable discussions. In the original version of our paper, Log
GW-invariants was called {\em relative GW-invariants}. Both author
would like to thank Qi Zhang to point out to us that our invariant
belongs to log category in algebraic geometry.  Both authors would
like to thank Yihong Gao for valuable discussions and the referee
for many suggestions to improve the presentation of this paper.
Thanks also to J. Robbin and  A. Greenspoon for editorial
assistance.

\section{Symplectic  surgery, flops and extremal transition}
    Symplectic surgeries have been extensively studied by a number of authors.
Many such surgeries  already appeared in \cite{Gr1}.
One of the oldest ones is the gluing along a contact boundary.
We do not know who was the first to propose  contact surgeries;
however, the second author benefited from a number
of stimulating discussions with Y. Eliashberg on contact surgeries over the
years.
During the last ten years, symplectic surgeries have been successfully used to
study symplectic topology, for example,
symplectic blow-up and blow-down by McDuff \cite{MS1} and
symplectic norm sum by Gompf \cite{Go2} and McCarthy and Wolfson \cite{MW}.
Very recently, Lerman \cite{L} introduced an operation called ``symplectic cutting''
which plays an important role in this article. For the
reader's convienence, we will give a quick
review of each of these operations.
The second half of this section is devoted to studying
algebro-geometric surgeries of flop and small transition in terms of
symplectic surgery.

\vskip 0.1in

\noindent
{\bf Definition 2.1  }{\it A contact structure $\xi$ on an odd dimensional
manifold
$N^{2n+1}$ is a codimension one distribution defined globably by a one-form
$\lambda$ ($Ker \lambda=\xi$) such that $\lambda\wedge (d\lambda)^{n}$ is a
volume
form. We call $\lambda$ a contact form.}

\vskip 0.1in

   Any contact manifold has a canonical orientation defined by the volume form
$\lambda\wedge (d\lambda)^n$. If $\lambda$ is a contact form, so is $f\lambda$
for a positive function $f$.
Any contact form defines the Reeb vector field $X_{\lambda}$ by the equation
$$
i_{X_{\lambda}}\lambda=1, i_{X_{\lambda}}d\lambda=0.\leqno(2.1)
$$
In general, the dynamics of $X_{\lambda}$ depends on $\lambda$.
There is a version of Moser's theorem for contact manifolds as follows.
Let $\lambda_t$ be a family of contact forms.
There is a family of diffeomorphisms
$$
\varphi_t: N\rightarrow N
$$
such that $\varphi_0=Id$ and $\varphi_t^*\lambda_t=f_t \lambda$, where $f_t$ is a
family of functions such that $f_0=1$.

\vskip 0.1in

\noindent
{\bf Example 2.2:  } The sphere $S^{2n-1}\subset {\bf C}^n$ with the standard symplectic
structure $\omega=\sum_i d x_i\wedge dy_i$ is a contact hypersurface.
The contact form is the restriction of $\lambda=\sum (x_i dy_i - y_i dx_i)$.
The closed orbit of the Reeb vector field is generated by complex multiplication
by $e^{i\theta}$. Furthermore, if $Q$ is a symplectic submanifold intersecting
$S^{2n-1}$ transversely, $Q\cap S^{2n-1}$ is a contact hypersurface of $Q$.

\vskip 0.1in

\noindent
{\bf Definition 2.3:  }{\it $(M, \omega)$ is a compact symplectic manifold with
boundary such that $\partial M=N$. $M$ has a contact boundary $N$ if in a
tubular neighborhood $N\times [0,\epsilon)$ of $N$, $\omega=d(f\lambda)$, where $f$ is
a function. Suppose that $X$ is an outward transverse vector field on $N$.
We say that $N$ is a convex contact boundary if $\omega(X_{\lambda}, X)>0$.
Otherwise, we call $N$ a concave contact boundary. In other words, the
induced orientation of a convex contact boundary coincides with the canonical
orientation. But the induced orientation of a concave contact boundary is
the opposite of the canonical orientation.}

\vskip 0.1in

Let $(M^+, \omega^+), (M^-, \omega^-)$ be symplectic manifolds with
contact boundary $N^+, N^-$. Furthermore, suppose that $N^+$ is a convex
boundary of $M^+$ and $N^-$ is a concave boundary of $M^-$. For any contact
diffeomorphisms $\varphi: N^+\rightarrow N^-$, we can glue them together to
form a  closed symplectic manifold $M_{gl}=M^+\cup_{\varphi} M^-$.
One can find  the details of the gluing construction from \cite{E}
(Proposition 3.1).

Recall that  GW-invariants only depend on the symplectic deformation class.

\vskip 0.1in

\noindent
{\bf Definition 2.4:  }{Two symplectic manifolds $(M, \omega), (M' \omega')$
with contact boundaries $(N, \lambda), (N', \lambda')$ are deformation
equivalent
if there is a diffeomorphism  $\varphi: M\rightarrow M'$, a family of symplectic
structures $\omega_t$ and a family of contact structures $\lambda_t$ such that
$(M, \omega_t)$ is a symplectic manifold with contact boundary $(N, \lambda_t)$
and
$$
\omega_0=\omega, \;\;\omega_1=\varphi^*\omega', \;\;\lambda_0=\lambda, \;\;
\lambda_1=\varphi^*\lambda'.\leqno(2.2)
$$
}

\vskip 0.1in

 The proof of the following Lemma is trivial. We omit it.

\vskip 0.1in

\noindent
{\bf Lemma 2.5:  }{\it If we deform the symplectic structures of $M^{\pm}$
according to  Definition 2.3, the glued-up manifolds $M_{gl}$ are
deformation equivalent.}

\vskip 0.1in

Another closely related surgery is symplectic cutting described in the introduction, where the Hamiltonian
vector field plays the role of the Reeb vector field.
The converse gluing process has been studied previously
by McCarthy and Wolfson \cite{MW}  in connection with symplectic norm sum.
By (1.10), we have a map
$$\pi: M\rightarrow \overline{M}^+\cup_Z \overline{M}^-.$$
It induces a map
$$\pi_*: H_i(M, \Z)\rightarrow H_i(\overline{M}^+\cup_Z \overline{M}^-, \Z).\leqno(2.3)$$
A class in $\ker \pi_*$ is called {\em a vanishing cycle}.  The virtual
neighborhood construction requires that the symplectic forms $\omega, \omega^+, \omega^-$
on $M, \overline{M}^+,
\overline{M}^-$ have integral periods. We can deform $\omega^+, \omega^-$ slightly (still denoted by
$\omega^+, \omega^-$ ) such that $\omega^+|_Z=\omega^-|_Z$ and
$\omega^{\pm}$ has rational periods. Then, we use the
Gompf-McCarthy-Wolfson gluing construction (an inverse operation of
the symplectic cutting) to obtain $\omega$. Then, $\omega$ has rational
periods since $\omega=\pi^*(\omega^+\cup_Z
\omega^-)$. To get integral periods, we just multiple $\omega$ by a
large integer.

\vskip 0.1in

\noindent
{\bf Example 2.6: }
 We consider the Hamiltonian action of $S^1$  on $({\bf C }^{n},\omega_0)$
corresponding to multiplication by $e^{it}$ with the Hamiltonian function
$H(z)= |z|^{2}-\epsilon$. Then we have a Hamiltonian action of $S^1$  on
$\left({\bf C }^{n}\times {\bf C},
-i(dw\wedge d\bar{w} + \Sigma^n_{i=1} dz_j\wedge d\bar{z}_j)\right)$ given by
$$
e^{it}(z,w)=\left( e^{it}z, e^{it}w\right)\leqno(2.4)
$$
with moment map
$$
\mu (z,w)=H(z) + |w|^2.\leqno(2.5)
$$
We have
$$
\mu^{-1}(0) = \left\{(z,w)| 0 < |w|^2 \leq \epsilon, H(z)=- |w|^2
\right\}\bigcup \left\{(z,0)|H(z)=0 \right\}\leqno(2.6)
$$
$$
\overline{M}^{+}= \mu^{-1}(0)/S^1.\leqno(2.7)
$$
where $M^{+}=\{z | |z|\leq \epsilon \}$ and $\overline{M}^{+}$ is
the corresponding symplectic cut. Choose the canonical complex structure $i$ in
${\bf C}^{n+1}={\bf C}^n\times {\bf C}$. The following fact is well-known
(see [MS1]) :
\vskip 0.1in
\noindent
{\it For ${\bf C}^{n}$ the symplectic cut $\overline{M}^{+}$ is ${\bf P}^{n}$
with symplectic form
$$-i \epsilon d\left(\frac{zd\bar{z}-\bar{z}dz}{1+\|z\|^2}\right)
. \leqno(2.8)$$ }

\vskip 0.1in

\noindent
{\bf Example 2.7: } We consider ${\cal O}(-1) + {\cal O}(-1)$ over ${\bf P}^1$. One
can construct a symplectic form as follows. Choose a hermitian metric on the vector bundle
${\cal O}(-1) + {\cal O}(-1)$. The $||z||^2$ for $z\in {\cal O}(-1) + {\cal O}(-1)$ is a
smooth function. $i\partial \bar{\partial} ||z||^2$ is a 2-form nondegenerate on the fiber.
Suppose that $\omega_0$ is a symplectic form on  ${\bf P}^1$.
$$
\omega=\pi^*\omega_0 + \epsilon i\partial \bar{\partial} ||z||^2 \leqno(2.9)
$$
is a symplectic form on the total space in a neighborhood of the zero section, where $\pi: {\cal O}(-1) + {\cal O}(-1)\rightarrow \P^1$
is the projection and $\epsilon$ is a small constant. The Hamiltonian function is
$$H(x,z_1,z_2) =|z_1|^2 + |z_2|^2 - \epsilon. \leqno(2.10)$$
It is a routine computation that the $S^1$-action is given by
$$
e^{it}(x,z_1,z_2)=\left(x, e^{it}z_1, e^{it}z_2\right).\leqno(2.11)
$$
We perform the symplectic cutting along the hypersurface $\widetilde{M}
=H^{-1}(0)$. By example 1,  we conclude that

\vskip 0.1in

\noindent
{\it For ${\cal O}(-1) + {\cal O}(-1)$ over ${\bf P}^1$ the symplectic cut
$\overline{M}^{+}$ is $P({\cal O}(-1) + {\cal O}(-1) + {\cal O})$. }

\vskip 0.1in

\noindent
{\bf Example 2.8: } Consider the algebraic variety $M$ defined by a homogeneous polynomial
$F$ in ${\bf C}^{n}$. Let $S^{2n-1}(\epsilon)$ be the sphere of radius
$\epsilon$. Let $\widetilde{M}=S^{2n-1}(\epsilon)\cap \{F=0\}.$ Consider the
Hamiltonian action of $S^1$  on $({\bf C }^{n},\omega_0)$
corresponding to multiplication by $e^{it}$ with the Hamiltonian
function
$$
H(z)= |z|^{2}-\epsilon.\leqno(2.12)
$$
Since $F$ is a homogeneous polynomial, $S^1$ acts on
$S^{2n-1}(\epsilon)\cap \{F=0\}$.
By example 1 we have

\vskip 0.1in

\noindent
{\it For $M=\{ F=0 \}$ the symplectic cut $\overline{M}^{+}$ is symplectic
deformation to the
variety defined by $F=0$ in ${\bf P}^{n}$. }

\vskip 0.1in

We have finished our digression about symplectic surgery. Next we study
flop and extremal transition from a symplectic point of view. Recall that
two projective manifolds $M$, $M'$ are birational equivalent iff some Zariski open
sets are isomorphic. If $M$, $M'$ are smooth CY 3-folds,
 $M$, $M'$ are related by a sequence of flops \cite{Ka}, \cite{K}. Namely,
there is a  sequence of smooth CY 3-folds $M_1, \ldots, M_k$
such that $M_1$ is obtained by a flop from $M$, $M_{i+1}$ is obtained by a
flop from $M_i$, and $M'$ is obtained by a flop from $M_k$. Furthermore,
 we can choose $M_i$ and the corresponding singular CY 3-folds
$(M_i)_s$ to be projective. Moreover, we can assume that the contractions
are primitive.
In this paper, we only deal with the small transition, which starts
from a small contraction. Suppose that $M, M_f$ is a pair of 3-folds connected by a
flop and $M_s$ is the corresponding singular projective 3-fold. The singularities of $M_s$
are rational double points. Wilson
observed that we can choose a local complex deformation around singularities
to split a rational double point into a collection of ordinary
double points (ODP). With its simultaneous resolution, it gives  local complex deformations of $M, M_f$
to deform an exceptional curve to a collection of rational curves with $\O(-1)+\O(-1)$
normal bundle. We can patch them with complex structures outside to define  almost
complex structures of $M, M_f$ which are tamed with symplectic structures. To relate $M, M_f$
as almost complex  manifolds, we can blow up the exceptional rational
curves (they are all $\O(-1)+\O(-1)$ curves) to obtain the same
almost complex  manifold $M_b$.  To summarize the construction,
\vskip 0.1in
\noindent
{\bf Proposition 2.9: }{\it As almost complex manifolds, $M, M_f$
described above have  the same  blow-up $M_b$.}
\vskip 0.1in
This is the model that we will use to compare invariants. Namely,
we will compare GW-invariants of $M, M_f$ with the log
GW-invariants of $M_b$. Clearly, $M, M_f$ have the same set of
exceptional curves. It is easy to see that the exceptional divisors
of $M_b$ is a disjoint union of $\P^1\times \P^1$'s.

   Next, we study the topology of flops. As we mentioned in the introduction, there is a homomorphism
$$
\varphi: H_2(M, \Z)\rightarrow H_2(M_f, \Z) \leqno(2.13)
$$
such that $\phi$ flips the sign of fundamental class of exceptional rational
curves.
The original proof used algebraic geometry. Here, we give a topological proof, which also fits with
our gluing construction.

 For any $A\in H_2(M, \Z)$, it
is represented
by  a 2-dimensional pseudo-submanifold $\Sigma$ of $M$. Since the exceptional set is a union of finitely many
curves, we can perturb  $\Sigma$ so that $\Sigma$ does not intersect
any  exceptional curve. Any
two different perturbations are pseudo-submanifold cobordant. Moreover, the
cobordism can also be pushed off the exceptional curves.
Hence, two different perturbations represent  the
same  homology class.
Let $\{\Gamma_k\}$ be the set of exceptional curves.
We have shown that the inclusion
$$
i: M-\cup \Gamma_k \rightarrow M
$$
induces an isomorphism on $H_2$. The same thing is true for $M_f$.
Moreover, $M-\cup \Gamma_k$ is the same as
$M_f-\cup (\Gamma_f)_k$ where $\{(\Gamma_f)_k\}$ is the corresponding set
of exceptional rational curves of $M_f$. Hence, $H_2(M, \Z)$ is isomorphic to
$H_2(M_f, \Z)$ and $\varphi$ is the isomorphism. Moreover, the same argument defines
an injective homomorphism
$$\phi_b: H_2(M, \Z)\rightarrow H_2(M_b, \Z) \leqno(2.14)$$ whose image consists of
elements with  zero intersection with exceptional divisors.
Let $\Gamma \subset M$,
$\Gamma_f\subset M_f$ be a pair of exceptional curves obtained by blowing
down the same exceptional divisor of $M_b$ over different rulings. We claim
$$
\varphi([\Gamma])=-[\Gamma_f].\leqno(2.15)
$$
Note that the normal bundle of an exceptional
divisor of $M_b$ is $\N=\O_1(-1)\otimes \O_2(-1)$ over $\P^1\times \P^1$,
where $\O_i(-1)$
means the $\O(-1)$ over the $i$-th factor of $\P^1\times \P^1$. Let
$\tau:  \P^1\rightarrow \P^1$ be the antipodal map which reverses the orientation.
Then the
restriction of $\N$ over $E=\{(x, \tau(x))\}$ is trivial. One can push $E$
off  $\P^1\times \P^1$.
We obtain the push-off of the exceptional curve of $M$ by projecting
the perturbation of $E$ to $M$. By the construction, $\Gamma, -\Gamma_f$
have the same push-off, where $-\Gamma_f$ means the opposite orientation. Hence
$$
\phi([\Gamma])=-[\Gamma_f].
$$

It was already observed in \cite{L} that the blow-up along a symplectic
submanifold
can be viewed as the symplectic cutting. For the reader's convenience, let's
construct the Hamiltonian $S^1$-action explicitly.
By the symplectic neighborhood theorem,
the symplectic structure of a neighborhood of $\Gamma$ is uniquely determined
by the  symplectic structure of $\Gamma$ and the almost complex structure of
the symplectic normal  bundle.
Moreover, a symplectic structure of $\Gamma$ is uniquely determined by
its volume.
It is enough to construct a specific symplectic structure on the total space of
$\O(-1)\oplus \O(-1)$ tamed to the complex structure and having the same volume as
that of $\Gamma$. Such a symplectic structure is constructed in
Example 2.7.
We perform the symplectic cutting  along the hypersurface
$N_{\epsilon}=\{||z||^2=\epsilon\}$ and obtain two symplectic manifolds
$\overline{M}^{+}$ and $\overline{M}^{-}$.
One is $M_b$. By  example 2,  the other one is the compactication of
$\O(-1)+\O(-1)$ at infinity
by $P(\O(-1)+\O(-1))$, which can be identified as the projectivization
$P(\O(-1)+\O(-1)+\O)$. We can summarize as
\vskip 0.1in
\noindent
{\bf Proposition 2.10: }{\it When we perform a symplectic cutting along
 the boundary of a tubular neighborhood of exceptional curves, we
 obtain $M_b$ and a collection of $P(\O(-1)+\O(-1)+\O)$'s. Here, $Z$ is a
 collection of $\P^1\times \P^1$'s.}
 \vskip 0.1in
The following lemma concerns  the vanishing two-cycles.
\vskip 0.1in
\noindent
{\bf Lemma 2.11: }{\it In the case of blow-up over a complex
codimension two submanifold,  there is no vanishing 2-cycle.}
\vskip 0.1in
{\bf Proof: } For any 2-cycle $\Sigma\subset
M$, using PL-transversality, we can assume that $\Sigma \subset
\overline{M}^-$. Hence, $\Sigma$ defines a homology class $\Sigma_b\in H_2(\overline{M}^-,
\Z)$.
However, $\overline{M}^-=M_b$ and there is a map $R: M_b\rightarrow
M$. Moreover, $R_*(\Sigma_b)=\Sigma.$ Hence, $\Sigma\neq 0$
implies $\Sigma_b\neq 0$.

   Next, we discuss the small transition. Again, Wilson's argument reduces to the
consideration of ODP.
Let $U_s\subset M_s$ be a neighborhood of the
singularity and $U_b\subset M_b$ be  its resolution.
We first standardize the symplectic
form over $0\in U_s\subset M_s$. Since $M_s$ is projective, we assume that $M_s$ is
embedded into a projective space. In the standard affine coordinates,
the Fubini-Study form can be written as
$$
\omega_0=\frac{i}{4\pi}d\left(\frac {\sum_i (x_id\bar{x}_i-\bar{x}_id x_i)}
{1+|x|^2}\right).\leqno(2.16)
$$
However, we must choose an
analytic change of coordinates to obtain the
standard form given by a quadratic equation.
Suppose that the change of coordinates is
$$
x_i=f_i(z_1, \cdots, z_4), \leqno(2.17)
$$
where $f_i(0)=0$ and $(\frac{\partial f_i}{\partial z_j}(0))$ is
nondegenerate. We use $L(f_i)$ to denote the linear term of $f_i$ and
$L(\omega_0)$ to denote the two-form obtained from $\omega_0$ by the linear
change of coordinates $(\frac{\partial f_i}{\partial z_j}(0))$.
Under such a change of coordinates,
$$\omega=L(\omega_0) + d(\tilde{\alpha}), \leqno(2.18)$$
where $\tilde{\alpha}=O(|z|^2)(dz_i+d\bar{z}_i)$.
Let $\beta_r$ be a cut-off function which equals  1 when $|z|>2r$ and
equals zero when $|z|<r$.
Let
$$
\omega_r=\omega-d(\beta_r\tilde{\alpha}).\leqno(2.19)
$$
$\omega_r$ is closed by construction and equal to $L(\omega_0)$ when
$|z|>2r$ and equal to $\omega$ when $|z|<r$.
Moreover,
$$
d(\beta_r\tilde{\alpha})=(d\beta_r)\tilde{\alpha}+\beta_r d\tilde{\alpha}
=O(r)(dz_i\wedge d\bar{z}_j).\leqno(2.20)
$$
Therefore, $\omega_r$ is nondegenerate for small $r$. Furthermore, by
Moser's theorem, it
is the same symplectic structure as $\omega$ on ${\bf P}^N$. Hence, we can assume
that $\omega=\omega_r$. We use the inverse of
$(\frac{\partial f_i}{\partial z_j}(0))$ to change $\omega$ to the standard
form $\omega_0$ and the equation to some homogeneous polynomial $F_2$ of
degree 2.

Let $S_r$ be the sphere of radius $r$. Here we choose $r$ small enough such
that $\omega=\omega_0$. Then, $S_r$ is a contact hypersurface of $\P^N$.
Let $N=S_r\cap\{F_2=0\}$. $N$ is a contact hypersurface of $M_s$.

   To obtain $M_t$, we have to deform the defining equation of $M_s$ to smooth
the singularity. This step can be described by a contact surgery.
Let $M^-_s=D_r\cap \{F_2=0\}$ where $D_r$ is the ball of
radius $r$ and $M^+_s=\overline{M_s-M^-_s}$.
The local smoothing is simply
$$\{F_2=t\}.\leqno(2.21)$$
Let
$$
M^-_{t_0}=D^r\cap \{F_2=t_0\}
$$
for small $t_0\neq 0$.
The global smoothing is well-understood by \cite{F}, \cite{Ti1} which requires some linear
conditions on the homology classes of the rational curves contracted. Here, we assume
that the global smoothing exists and that smoothing is $M_t$. We use the previous
construction to decompose
$$
M_t=M^+_t \cup_{id} M^-_t,\leqno(2.22)
$$
where $M^-_t=D_r\cap M_t$ and $M^+_t$ is the complement. However,
$M^+_t$ is symplectic deformation equivalent to $M^+_s$. $M^-_t$ is deformation
equivalent to $M^-_{t_0}$. Hence, $M_t$ is
symplectic deformation equivalent to
$$
M_e=M^+_s\cup M^-_{t_0}.\leqno(2.23)
$$
This is the symplectic model which we will use to calculate GW-invariants.
Clearly, if we choose $r'$ slightly larger than $r$, the manifold $M_e$ also admits
an $S^1$-Hamiltonian function near $M_e\cap S^{r'}$. Then we can perform the
symplectic cutting to obtain two symplectic manifolds. We observed
that
one of them is precisely $M_b$. The other one is symplectic deformation equivalent to
$$\{F_2-t_0 x^{2}_5=0\}\subset \P^4.$$
We summarize the previous construction as
\vskip 0.1in
\noindent
{\bf Proposition 2.12:}{\it We can perform a symplectic cutting on $M_e$
to obtain $M_b$ and a collection of quadric 3-folds.}
\vskip 0.1in

   By the previous argument, we can push  any 2-dimensional homology class of
$M$ off the exceptional loci, and hence off the singular points of $M_s$.
Therefore, we have a map
$$
\varphi: H_2(M, \Z)\rightarrow H_2(M_e, \Z).\leqno(2.24)
$$
The topological description of small transition along ODP was known to H. Clemens
a long time ago. Note that the boundary of a tubular neighborhood of $\O(-1)+\O(-1)$
is $S^2\times S^3$. The transition is an interchange of the handle $S^2\times D^4$
with $D^3\times S^3$. The singular 3-fold $M_s$ corresponds to the cone over $S^2\times S^3$.
In particular, Clemens showed that $M^-_{t_0}$ is homeomorphic to $S^3\times D^3$. It
doesn't  carry any nontrivial 2-dimensional homology class in the interior.
Therefore,  $\varphi$ (2.24) is surjective. Previous argument also
show
\vskip 0.1in
\noindent
{\bf Lemma 2.13: }{\it There
is no vanishing 2-cycle for small transition. }
\vskip 0.1in
Unfortunately, the authors do not know any direct symplectic
surgery from $M$ to $M_e$ without going through $M_b$.

Recall that a holomorphic small transition exists iff $\sum_i a_i [\Gamma_i]=0$ for
$a_i\neq 0$, where $\{\Gamma_i\}$ is the set of exceptional curves \cite{F}, \cite{Ti1}.

However, we have obtained a stronger result in the symplectic
category.
\vskip 0.1in
\noindent
{\bf Theorem 2.14: }{\it Suppose that $V_s$ is projective with rational double points. We always
 have a symplectic manifold $V_e$ obtained by gluing with local smoothing of singularities.
In the holomorphic case, $V_e$ is symplectic deformation equivalent to a holomophic small transition.
}
\vskip 0.1in
\noindent
{\bf Remark 2.15: }  Finally, we make a remark about general extremal transitions.
Many interesting examples are constructed by  nonprimitive extremal transitions.
A natural question is if we can decompose an extremal transition as a sequence of
primitive transitions. Our theorems actually give a criterion to this question.
After we perform a primitive extremal transition, we need to study the change of
the K\"ahler cone. Suppose that $E$ gives a different extremal ray. Then it is easy to
observe that $E$ will remain  an extremal ray if its GW-invariants are not
zero. Our theorem shows that the GW invariant can be calculated using the GW-invariant
of the Calabi-Yau 3-fold before the transition. Moreover,
GW-invariants of a crepant resolution depend only on a neighborhood of the exceptional
loci and can be calculated indepedently from Calabi-Yau 3-folds.

\vskip 0.1in

\noindent
{\bf Remark 2.16: }  We would like to make another remark about the
smoothing of singular Calabi-Yau 3-folds.
One can study the smoothing theory by studying the local smoothings of singularities
and
the global extension of a local smoothing. The singular CY 3-folds obtained by
contracting a smooth CY 3-fold can only have so called canonical singularities.

There are many interesting phenomena
in the smoothing theory of CY 3-folds with canonical singularities.
For example, M. Gross showed that the same singular CY 3-fold could have
two different smoothings \cite{Gro1}. They give a pair of examples
where two diffeomorphic CY 3-folds carry different quantum cohomology
structures \cite{R1}.
By our arguments, it is clear that if a  singular Calabi-Yau 3-fold
has a local smoothing for its singularity, it has a global symplectic smoothing.
This is not true in the algebro-geometric category by Friedman's results.
This indicates an exciting possibility that there is perhaps a theory of symplectic
Calabi-Yau 3-folds, which
is broader than  the algebro-geometric theory.
Such a theory would undoubtedly be  important for
the classification of CY 3-folds itself.

\section{\bf Compactness Theorems}

Let $(M,\omega)$ be a compact symplectic manifold of dimension
$2n+2$, $H:M\rightarrow {\bf R}$ a local Hamiltonian function such
that there is a small interval $I=(-\ell, \ell)$ of regular
values. Denote $\widetilde{M} =H^{-1}(0)$. Suppose that the
Hamiltonian vector field $X_H$ generates a circle action on
$H^{-1}(I)$. There is a circle bundle $\pi : \widetilde{M}
\rightarrow Z=\widetilde{M}/S^1 $ and a natural symplectic form
$\tau_0$ on $Z$. We may choose a connection 1-form $\lambda $ on
$\widetilde{M}$ such that $\lambda(X_H)=1$ and $d\lambda$
represents the first Chern class for the circle bundle (see
[MS1]). Denote $\xi=ker(\lambda)$. Then $\xi$ is an
$S^1$-invariant distribution and $(\xi,
\pi^{\ast}\tau_0)\rightarrow \widetilde{M}$ is a $2n$-dimensional
symplectic vector bundle. We identify $H^{-1}(I)$ with $I \times
\widetilde{M}$. By a uniqueness theorem on symplectic forms (see
[MS1], [MW]) we may assume that the symplectic form on
$\widetilde{M}\times I$ is expressed by $$ \omega =
\pi^{\ast}(\tau_0 + yd\lambda) - \lambda \wedge dy. \leqno(3.1)$$
We assume that the hypersurface $\widetilde{M}=H^{-1}(0)$ devides
$M$ into two parts $M^+$ and $M^-$, which can be written as
$$M^{+}_{0}\bigcup\left\{[-\ell,0)\times \widetilde{M}\right\},$$
$$M^{-}_{0}\bigcup\left\{(0,\ell]\times \widetilde{M}\right\},$$
where $M^{+}_{0}$ and $M^{-}_{0}$ are compact manifolds with
boundary. We mainly discuss $M^+$; the discussion for $M^-$ is
identical. Let $\phi^{+} :[0, \infty)\rightarrow [- \ell, 0)$ be a
function satisfying $$(\phi^{+})^{\prime}>0,
\;\;\phi^{+}(0)=-\ell, \;\;\phi^{+}(a)\rightarrow 0 \;\;as\;\; a
\rightarrow \infty.\leqno(3.2)$$ Through $\phi^{+}$ we consider
$M^+$ to be $M^{+} =  M^{+}_{0}\bigcup\{[0, \infty)\times
\widetilde{M}\}$ with symplectic form
$\omega_{\phi^+}|_{M^{+}_0}=\omega$, and over the cylinder ${\bf
R}\times \widetilde{M}$
$$\omega_{\phi^{+}} = \pi^{\ast}(\tau_0 +
\phi^{+}d\lambda) - (\phi^+)^{\prime}\lambda \wedge d a.
\leqno(3.3)$$
 Moreover, if we choose the origin of ${\bf R}$
tending to $\infty $, we obtain ${\bf R}\times \widetilde{M}$ in
the limit. Later we will fix the value of $\phi^{+}$ on $[0, 1]$.
Choose $\ell_0 < \ell $ and denote $$\Phi^+ =\left \{ \phi :[1,
\infty )\rightarrow [-\ell_0, 0) | \phi^{\prime}
> 0 \right \}.$$
Let $\ell_1 < \ell_2$ be two real numbers satisfying $ - \ell_0
\leq \ell_1 < \ell_2 \leq 0 .$ Let $\Phi_{\ell_1,\ell_2}$ be the
set of all smooth functions $\phi:{\bf R}\rightarrow
[\ell_1,\ell_2]$ satisfying $$\phi^{\prime}>0,\;\;
\phi(a)\rightarrow \ell_2\;\:\;{\rm as} \;a\rightarrow \infty ,\;\;
\phi(a)\rightarrow \ell_1\;\:\;{\rm as} \;a\rightarrow -\infty.$$ To
simplify notations we use $\Phi $ to denote both $\Phi^+$ and
$\Phi_{\ell_1,\ell_2}$, in case this does not cause confusion.
\vskip 0.1in \noindent We choose a compatible almost complex
structure $\widetilde{J}$ for $Z$ such that
$$g_{\widetilde{J}(x)}(h,k)=\tau_0(x)(h,\widetilde{J}(x)k),$$ for
all $h, k \in TZ $, defines a Riemannian metric. $\widetilde{J}$
and $g_{\widetilde{J}}$ are lifted in a natural way to $\xi $. We
define an almost complex structure $J$  on ${\bf R}\times
\widetilde{M}$ as follows: $$J\frac{\partial}{\partial
a}=X_{H},\;\;\; JX_{H}=-\frac{\partial}{\partial a},\leqno(3.4) $$
$$J\xi=\xi\;\;J|_{\xi}=\widetilde{J}.\leqno(3.5)$$ Since
$g_{\widetilde{J}}$ is positive, and $d\lambda$ is a 2-form on $Z$
(the curvature form), by choosing $\ell$  small enough we may
assume that $\widetilde{J}$ is tamed by $\tau_0 + yd\lambda $ for
$|y| < \ell $, and there is a constant $C>0$ such that
$$\tau_0(v,\widetilde{J}v)\leq C\left(\tau_0(v,\widetilde{J}v)+
yd\lambda(v, \widetilde{J}v)\right)\leqno(3.6)$$ for all $v\in TZ,
|y|\leq \ell $. Then $J$ is $\omega_{\phi}$-tamed for any $\phi
\in \Phi $ over the tube. \vskip 0.1in \noindent We denote by $N$
one of $M^{+}$, $M^-$ and ${\bf R}\times\widetilde{M}$. We may
choose an almost complex structure $J$ on $N$ such that\\ i) $J$
is tamed by $\omega_{\phi} $ in the usual sense, \par \noindent
ii) Over the tube ${\bf R }\times \widetilde {M},$  (3.4) and
(3.5) hold. \vskip 0.1in \noindent Then for any $\phi \in \Phi $
$$\langle v,w\rangle_{\omega_{\phi}} = \frac{1}{2}\left(
\omega_{\phi} (v,Jw) + \omega_{\phi} (w,Jv)\right) \;\;\;\;\;
\forall \;\; v, w \in TN\leqno(3.7)$$
defines a Riemannian metric
on $N$. Note that $\langle \;,\;\rangle_{\omega_{\phi}}$ is not
complete. We choose another metric $\langle \;,\;\rangle$ on $N$
such that $$\langle \;,\;\rangle = \langle
\;,\;\rangle_{\omega_{\phi}} \;\;\;\;on \;\; M^{+}_{0}$$ and over
the tubes $$\langle(a,v),(b,w) \rangle= ab + \lambda (v)\lambda
(w) + g_{\widetilde{J}}(\Pi v, \Pi w),\leqno(3.8)$$ where we
denote by $\Pi:T\widetilde{M}\rightarrow\xi$ the projection along
$X_H$. It is easy to see that $\langle \;,\;\rangle$ is a complete
metric on $N$. \vskip 0.1in \noindent Let $(\Sigma,i)$ be a
compact Riemann surface and $P\subset\Sigma$ be a finite
collection of points. Denote $\stackrel{\circ}{\Sigma}
=\Sigma\backslash P.$ Let ${u}:\stackrel{\circ}{\Sigma}
\rightarrow N$ be a ${J}$-holomorphic map, i.e., ${u}$ satisfies
$$d{u}\circ i={J}\circ d{u}.\leqno(3.9)$$ Following [HWZ1] we
impose an energy condition on $u$. For any $J$-holomorphic map
$u:\stackrel{\circ}{\Sigma}\rightarrow N$ and any $\phi \in \Phi $
the energy $E_{\phi}(u)$ is defined by
$$E_{\phi}(u)=\int_{\stackrel{\circ}{\Sigma}}u^{\ast}\omega_{\phi}.\leqno(3.10)$$
A $J$-holomorphic map $u:\stackrel{\circ}{\Sigma} \rightarrow N $
is called a finite energy $J$-holomorphic map if $$Sup_{\phi \in
\Phi }\left \{\int_{\stackrel{\circ}{\Sigma}}u^{*} \omega_{\phi}
\right \}<\infty. \leqno(3.11)$$ We shall see later that the
condition is natural in view of our surgery. For a $J$-holomorphic
map $u:\stackrel{\circ}{\Sigma} \rightarrow {\bf
R}\times\widetilde{M}$  we write $u=(a, \widetilde{u})$ and define
$$\widetilde{E}(u)=\int_{\stackrel{\circ}{\Sigma}}\widetilde{u}^{\ast}
(\pi^{\ast}\tau_0).\leqno(3.12)$$ Let $z=e^{s+2 \pi it}.$ One
computes over the cylindrical part $$u^{\ast}\omega_{\phi}=
(\left(\tau_0 + \phi d\lambda \right)\left((\pi\widetilde{u})_s,
(\pi\widetilde{u})_t)\right) + {\phi}^{\prime}(a^{2}_s + a^{2}_t
))ds\wedge dt,\leqno(3.13)$$ which is a nonnegative integrand .

In this section we shall prove a compactness theorem for $J$-
holomorphic curves in $M^{\pm}$, and a convergence theorem for
$J$-holomorphic curves when $M$ is stretched to infinity along $\widetilde{M}$.
There are various existing compactness theorems over compact manifolds. (see [RT1],[FO],[MS],[Ye]).
In our situation the manifold is not compact, so we have to analyse the
behaviour of sequences of holomorphic curves at infinity. In the
case of contact manifolds, a similar analysis has been done by
Hofer and his collabrators \cite{HWZ1}, \cite{HWZ2} \cite{HWZ3}. However, the
case under consideration is not contact in general and their
analysis does not directly apply. Here, we carry out the analysis in the case
of a Hamiltonian $S^1$-action where we follows closely Hofer's approach.
We remark that the existence of an $S^1$-action simplifies  many arguments
so that we have a  simpler proof than in the case of contact manifolds.

\subsection{\bf Convergence to periodic orbits}

We fix an integer $r>2$, let $W_r^2(S^1,\widetilde M)$ be the
space of $W_r^2$- maps from $S^1$ into $\widetilde M$, which is a
Hilbert manifold. For any $\gamma \in W_r^2(S^1,\widetilde M)$ the
tangent space $T_{\gamma}W_r^2(S^1,\widetilde M)$ is the space of
vector fields $\eta\in W_r^2(\gamma^*T\widetilde M)$ along
$\gamma$. The Riemannian metric $(\;,\;)$ on $W_r^2(S^1,\widetilde
M)$ is defined by
$$(\eta_1,\eta_2):=\int_{S^1}\langle\eta_1,\eta_2\rangle dt.
\leqno(3.14)$$ Let $x(t) \in Z $ be an orbit of the circle action.
For any integer $k$ let $x_{k}(t):=x(k t)$, which is called a
$k$-periodic orbit. We can consider $x_k(t)$ as a element in
$W_r^2(S^1,\widetilde M)$. Now we define an action functional in a
neighbourhood of $x_k(t)$ in $W_r^2(S^1,\widetilde M)$. Choose an
$\epsilon$-ball $O_{x_k,\epsilon}$ of  $0$ in the Hilbert space
$T_{x_k}W_r^2(S^1,\widetilde M)$ such that the exponential map
$$\exp_{x_k}:O_{x_k,\epsilon}\rightarrow W_r^2(S^1,\widetilde M)$$
identifies $O_{x_k,\epsilon}$ with a neighbourhood of $x_{k}(t)$
in $W_r^2(S^1, \widetilde M)$. Since $x_k:S^1\rightarrow
\widetilde M$ is an immersion and the points in $O_{x_k,\epsilon}$
are near $x_k$ with respect to the $W_r^2$-norm, we can assume
that for any $\gamma\in O_{x_k,\epsilon}$, $\gamma:S^1\rightarrow
\widetilde{M}$ is an immersion. For any $W_r^2$-loop $\gamma$ in
this neighbourhood let $\eta \in O_{x_k,\epsilon}$ be the
corresponding vector field along $x_k(t)$. Denote by $W: S^1\times
[0,1]\rightarrow \widetilde{M}$ the annulus defined by $\{
\exp_{x_k(t)}s\eta |0\leq s \leq 1, t\in S^1 \}$. We define an
action functional by $${\cal A}(\gamma)=-\int_{S^1\times
[0,1]}W^*\pi^* \tau_0.\leqno(3.15)$$ Suppose $\gamma_p$ is a
smooth curve in $W_r^2(S^1,\widetilde  M)$, with
$\gamma_0=\gamma.$ Then
$\eta=\left(\frac{d\gamma_p}{dp}\right)_{p=0}$ is a vector field
along $\gamma$. One can easily derive the first variational
formula: $$d{\cal A}(\gamma)\eta = \left.\frac{d{\cal
A}(\gamma_p)}{dp}\right|_{p=0}
=-\int_{S^1}\pi^{\ast}\tau_0(\dot{\gamma},\eta)dt=\int_{S^1}\langle\Pi\dot{\gamma},
\widetilde{J}\Pi\eta\rangle dt.\leqno(3.16)$$ It follows that
$\gamma $ is a critical point if and only if $\dot{\gamma}$ is
parallel to $X_H$ everywhere. So every $k$-orbit is a critical
point of ${\cal A}$. \vskip 0.1in \noindent Denote by $S_k$ the
set of $k$-periodic orbits. It is standard to prove that $S_k$ is
compact. Obviously, ${\cal A}(\gamma)$ is constant on the set
$S_k$. We now calculate the second variation for ${\cal A}$.
Consider a 2-parameter variation $(t,w,v)$ of $\gamma$, i.e., a
map $F: Q\rightarrow \widetilde{M}$ such that
$F(t,0,0)=\gamma(t)$, where $Q=S^1\times
(-\epsilon,\epsilon)\times (-\delta,\delta)$, and $F$ is smooth
with respect to $(w,v)$. Let $\tau,\eta, \zeta $ be vector fields
corresponding to the first, second and third variable of $F$,
respectively. Suppose that $\gamma$ is a critical point of ${\cal
A}$; then  $\tau(t,0,0)$ is parallel to $X_H $. We have
$$\frac{\partial{\cal A}}{\partial w}=
\int_{S^1}\langle\Pi\tau,\widetilde{J}\Pi\eta\rangle dt$$
$$\left.\frac{\partial^2 {\cal A}}{\partial v\partial
w}\right|_{(0,0)}= \int_{S^1}\langle\nabla_{\zeta}(\Pi\tau),
\widetilde{J}\Pi\eta\rangle dt$$ where $\nabla$ denotes the
covariant derivative with respect to the Levi-Civita connection of
$\langle\;,\;\rangle $ on $\widetilde{M}$. Denote $\zeta^{\perp}:=
\Pi\zeta$, the normal component of $\zeta$, and
$T^{\perp}_{\gamma}:=\{\zeta|\langle \zeta,\dot{\gamma}\rangle =0
\}$ . Note that at $(t,0,0)$
$$\nabla_{\zeta}(\Pi\tau)=\nabla_{\zeta}\tau - \nabla_{\zeta}
(\langle \tau, X_H \rangle X_H)\leqno(3.17)$$
$$=\nabla_{\tau}\zeta - \langle \tau, X_H \rangle
\nabla_{\zeta}X_H \;\;\;\mbox{mod} X_H$$
$$=\nabla_{\tau}\zeta^{\perp} - \langle \tau, X_H \rangle
\nabla_{\zeta^{\perp}}X_H \;\;\mbox{mod} X_H.$$ By our choice of
the Riemannian metric and the almost complex structure
$$\frac{D\zeta}{dt}\in T_{\gamma}^{\bot}\;\;\;\forall \zeta \in
T_{\gamma}^{\bot}.$$ Define the linear transformation
$S:T_{\gamma}^{\perp}\rightarrow T_{\gamma}^{\perp}$  by $\zeta
\rightarrow \widetilde{J}\nabla_{\zeta}X_H$. We get the following
variational formula: $$\left.\frac{\partial^2 {\cal A}}{\partial
v\partial w}\right|_{(0,0)} =-\int_{S^1}\langle
(\widetilde{J}\Pi\frac{D}{dt} + \langle \dot{\gamma}, X_H \rangle
S) \zeta^{\perp},\eta^{\perp}\rangle dt\leqno(3.18)$$ where $D/dt$
denotes the covariant derivative along $\gamma$. It is easy to
show that the second variation formula is a symmetric bilinear
form on the space of smooth vector fields along $\gamma$, denoted
by $I(\zeta,\eta)$. Let $P=-\widetilde{J}\Pi\frac{D}{dt} - \langle
\dot{\gamma}, X_H \rangle S$. We have
$$I(\zeta,\eta)=(P\zeta,\eta).\leqno(3.19)$$ Since the inclusion
$W_1^2(S^1) \rightarrow L^2(S^1)$ is compact, and $P$ is a compact
perturbation of the selfadjoint operator $J_0\frac{d}{dt}$, the
spectrum $\sigma(P)$ consists of isolated eigenvalues, which
accumulate at $\pm \infty$. \vskip 0.1in \noindent {\bf Remark 3.1
} {\it Let $\varphi : S^1\rightarrow S^1 $ be a diffeomorphism .
Then $x_k\circ \varphi $ is also a critical point of ${\cal A}$,
and the second variational formula (3.18) remains valid.} \vskip
0.1in Now we consider a smooth map $F: S^1\times
(-\epsilon,\epsilon) \rightarrow \widetilde{M}$ such that
$F(t,0)=\gamma(t)$ and such that at every point $(t,w)$,
$\Pi\tau=0$, i.e, for every $w$, $F(\cdot,w)$ is a critical point
of ${\cal A}$. It follows from (3.17) that $\zeta $ satisfies the
following equation: $$\nabla_{\tau}\zeta^{\perp} - \langle \tau,
X_H \rangle \nabla_{\zeta^{\perp}}X_H =0.  \;\;\leqno(3.20)$$ Put
$T_{x_k}^{\bot}=\{\eta\in T_{x_k}W_r^2(S^1,\widetilde
M)|\eta\bot\dot{x}_k \}.$ $T_{x_k}^{\bot}$ is a Hilbert space. Let
$C_r({\cal A})$ be the set of critical points of ${\cal A}$, and
denote $K=:C_r({\cal A})\bigcap T_{x_k}^{\bot}$. \vskip 0.1in
\noindent {\bf Proposition 3.2 }\\ 1) $K$ is a smooth submanifold
in $T_{x_k}^{\bot}$,\par \vskip 0.1in \noindent 2) At any $\gamma
\in K$, the restriction of the index $I(\;,\;)$ to the normal
direction of $K$ in $T_{x_k}^{\bot}$ is nondegenerate. \vskip
0.1in \noindent {\bf Proof:} It is easy to show that $K$ is a
smooth submanifold of dimension $2n$ in $T_{x_T}^{\bot}$ (for the
definition of submanifold in a Hilbert space see [Kli]). We prove
2). Let $\gamma \in K $. It follows from (3.20) that any tangent
vector $\zeta$ of $K$ at $\gamma$ satisfies
$$\widetilde{J}\frac{D\zeta^{\perp}}{dt} + \langle \dot{\gamma},
X_H \rangle S\zeta^{\perp} = 0.\leqno(3.21)$$ On the other hand,
any solution of the equation (3.21 ) is determined by its value
$\zeta^{\perp}(0)$ at $t=0$. Since $\zeta^{\perp}(0) \perp
\dot{\gamma}$ we have $dim ker(P)\leq 2n$. Hence the tangent space
of $K$ at $\gamma$ coincides with $ker P.$ The conclusion follows.
$\;\;\;\;\Box$ \vskip 0.1in \noindent {\bf Remark 3.3 }{\it The
set of critical points satisfying 1) and 2) in Proposition 3.2 is
usually called of Bott-type.} \vskip 0.1in \noindent {\bf
Proposition 3.4 } {\it Let $x_k(t) \in S_k$. There exists a
neighborhood $O$ of $\{x_k(t+d),0\leq d \leq 1\}$, in
$W_r^2(S^1,\widetilde{M})$ and a constant $C>0$ such that the
inequality $$\|\nabla {\cal A}\|_{L^2(S^1)}\ge C|{\cal A}|^{1\over
2}\leqno(3.22)$$ holds in $O$.} \vskip 0.1in \noindent {\bf
Proof:}  We first prove that (3.22) holds in a neighborhood of $0$
in $T_{x_k}^{\perp}$. Using the Morse lemma with parameters we can
find a neighborhood $O_{\epsilon^{'}}\subset O_{\epsilon}$ of $0$
such that the following holds: there are diffeomorphisms
$\varphi_x: N_x\rightarrow N_x$, which depend continuously on
$x\in K\cap O_{\epsilon^{'}}$, such that under the diffeomorphisms
the function ${\cal A}|_{N_x}$ has the form $${\cal
A}(y)=(P'(x)y,y),\;\;\;\forall y \in N_x$$ where $P'$ denotes the
restriction of $P(x)$ to $N_x$. It suffices to prove the
inequality for the quadratic function ${\cal A}(y)=(P'y,y)$. By
the nondegeneracy of $I(\;,\;)$ we can find a constant $C_2>0$
such that $|\lambda_i|\ge C_2$ for all $x\in K\cap O_{\epsilon'}$
and all eigenvalues $\lambda_i$. Then $$\|\nabla {\cal
A}\|_{L^2(S^1)}\ge \|\nabla_y{\cal A}\|_{L^2(S^1)}=
2(P'(x)y,P'(x)y)^{\frac{1}{2}}$$
$$\ge\sqrt{C_2}(P'(x)y,y)^{\frac{1}{2}}=\sqrt{C_2}|{\cal
A}|^{\frac{1}{2}},$$ where $\nabla_y{\cal A}$ denotes the gradient
of ${\cal A}|_{N_x}$ as a functional on $N_x$. Now let $\eta\in
W_r^2(x_k^{*}T\widetilde M)$ be a vector field along $x_k$ such
that $\|\eta\|_r$ is very small. For any $t\in S^1$ we can find a
unique $t^{\prime}\in S^{1}$ and a vector $\eta^{\prime}\in
T_{x_k(t^{\prime})} ^{\bot}\widetilde{M}$ such that
$$\exp_{x_k(t)}\eta=\exp_{x_k(t^{\prime})}\eta^{\prime}.$$ This
induces a map $\varphi:S^1\rightarrow S^1,\;t\longmapsto
t^{\prime}.$ When $\|\eta\|_r$ is very small $\varphi$ is a
diffeomorphism with $\varphi^{\prime}\approx 1$.  Since ${\cal A}$
is invariant under the diffeomorphism group $S^1\rightarrow S^1$
we can find a neighbourhood $O$ of $x_k(t+d),\;\:0\le d \le 1$, in
$W_r^2(S^1,\widetilde{M})$ such that (3.22) holds. The proposition
follows.$\Box$ \vskip 0.1in \noindent We are interested in the
behaviors of the finite energy $J$-holomorphic maps near a
puncture. There are two different types of puncture : the
removable singularities and the non-removable singularities. If
${u}$ is bounded near a puncture, then this puncture is a
removable singularity. In the following, we assume that all
punctures in $P$ are non-removable. Then ${u}$ is unbounded near
the punctures. By using the same method as in [H], one can prove
the following two lemmas: \vskip 0.1in \noindent {\bf Lemma 3.5 }
\\ {\bf (1)} Let $u=(a,\widetilde{u}):{\bf C} \rightarrow {\bf R}
\times \widetilde{M}$ be a $J$-holomorphic map with finite energy.
If $\int_{{\bf C}}\widetilde{u}^{\ast}(\pi^{\ast}\tau_0)=0$, then
$u$ is a constant.\\ {\bf (2)} Let $u=(a,\widetilde{u}):{\bf
R}\times S^1 \rightarrow {\bf R} \times \widetilde{M}$ be a
$J$-holomorphic map with finite energy. If $\int_{{\bf R}\times
S^1}\widetilde{u}^{\ast}(\pi^{\ast}\tau_0)=0$, then
$(a,\widetilde{u})=(ks+c, kt+d)$, where $k\in {\bf Z}$, $c$ and
$d$ are constants.
\vskip 0.1in
 \noindent
 {\bf Lemma 3.6 } Let
$u=(a,\widetilde{u}):{\bf C} \rightarrow {\bf R} \times
\widetilde{M}$ be a nonconstant $J$-holomorphic map with finite
energy. Put $z=e^{\delta+2\pi it}$. Then for any sequence
$s_i\rightarrow \infty $ , there is a subsequence, still denoted
by $s_i$, such that $$lim_{i\rightarrow \infty}\widetilde
u(s_i,t)=x(kt)$$ in $C^{\infty}(S^1)$ for some $k$-periodic orbit
$x_k$. \vskip 0.1in \noindent {\bf Theorem 3.7 }{\it Let
$u=(a,\widetilde{u}):{\bf C} \rightarrow {\bf R} \times
\widetilde{M}$ be a $J$-holomorphic map with finite energy. Put
$z=e^{s+2 \pi it}$. Then $$\lim_{s\rightarrow \infty}\widetilde
u(s,t)=x(kt)$$ in $C^{\infty}(S^1)$ for some $k$-periodic orbit
$x_k$. Moreover, the convergence is to be understood as
exponential decay uniformly in t.} \vskip 0.1in \noindent {\bf
Proof: } Denote
$$\widetilde{E}(s)=\int_s^{\infty}\int_{S^1}\widetilde{u}^{\ast}(\pi^{\ast}
\tau_0).$$ Then $$\widetilde{E}(s)=\int_s^{\infty}\int_{S^1}|\Pi
\widetilde{u}_t|^2dsdt,$$
$$\frac{d\widetilde{E}(s)}{ds}=-\int_{S^1}|\Pi
\widetilde{u}_t|^2dt. \leqno(3.23)$$ We choose a sequence $s_i$
such that $$lim_{i\rightarrow \infty}\widetilde
u(s_i,t)=x(kt)=:x_k(t)$$ for some $k$-periodic orbit $x$. For
every loop $\widetilde u(s,\cdot) \in O_{x_k,\epsilon}$ we draw an
annulus $W(s)$ as before and define the functional ${\cal
A}(s):={\cal A}(\widetilde u(s,\cdot))$. Since the annulus $W(s)$
varies continuous as $s$ varies we have
$\widetilde{E}(s)-\widetilde{E}(s_i)=|{\cal A}(s)-{\cal A}(s_i)|$,
where $s < s_i$. Let $i\rightarrow\infty;$ then
$$\widetilde{E}(s)=|{\cal A}(s)|.\leqno(3.24)$$ By Proposition
3.4, there is a $W_r^2$-neighbourhood $O$ of $\{x_k(t+d),0\leq d
\leq 1\}$ in which the inequality (3.22) holds. For $\widetilde
u(s,\cdot) \in O$ we have $$\frac{d\widetilde{E}(s)}{ds}=
-\int_{S^1}|\Pi \widetilde{u}_{t}|^{2}dt=-\|\bigtriangledown {\cal
A}(s)\|^2_{L^2(S^1)}\leq -C^2|{\cal
A}(s)|\leq-C^2\widetilde{E}(s)\leqno(3.25)$$
$$\frac{d\widetilde{E}(s)}{ds}\leq -C\|\Pi
\widetilde{u}_t\|_{L^2(S^1)} \widetilde{E}(s)^{1\over
2}.\leqno(3.26)$$ Suppose that $\widetilde u(s,\cdot) \in O$ for
$s_0\leq s \leq s_1$, and $\widetilde{E}(s)\neq 0\;\;\forall s \in
(s_0,s_1).$ It follows from (3.25) and (3.26) that
$$\widetilde{E}(s_1)\leq\widetilde{E}(s_{0})e^{-C^2(s_1-s_{0})}
\leqno(3.27)$$ $$\int_{s_0}^{s_1}\|\Pi
\widetilde{u}_t\|_{L^2(S^1)}ds\leq {1\over C}\left(\widetilde{E}
(s_0)^{1/2}-\widetilde{E}(s_1)^{1/2}\right)\leqno(3.28)$$ Denote
by $\widetilde{d}$ the distance function defined by the metric
$g_{\widetilde{J}}$ on $Z$. We have
$$\int_{S^1}\widetilde{d}(\widetilde u(s,t),\widetilde
u(s_i,t))dt\leq \int_{s}^{s_i}\|\Pi \widetilde{u}_t\|_{L^2}ds
\leqno(3.29)$$ $$\leq {1\over C}\left(\widetilde{E}(s)^{1/2}
-\widetilde{E}(s_i)^{1/2}\right)\leq {1\over
C}\left(\widetilde{E}(s) - \widetilde{E}(s_i)\right)^{1/2}.$$ We
show that for any $C^{\infty}$-neighbourhood $U$ of
$\{x_k(t+d),0\leq d \leq 1\}$ there is a $N>0$ such that if $s>N$
then $\widetilde u(s,\cdot)\in U.$ If not, we could find a
neighbourhood $U\subset O$ and a subsequence of $s_i$ (still
denoted by $s_i$) and a sequence $s_i^{\prime}$ such that
$$\widetilde u(s,\cdot)\in O, \;\; for \;\; s_i \leq s \leq
s_i^{\prime}, \leqno(3.30)$$ $$\widetilde
u(s_i^{\prime},\cdot)\notin U.\leqno(3.31)$$ By Lemma (3.6) and by
choosing a subsequence we may assume that $$\widetilde
u(s_i^{\prime},t)\rightarrow x^{\prime}(k^{\prime}t) \hskip 1cm
{\rm in}\,\, C^{\infty}(S^1,\widetilde{M})$$ for some
$k^{\prime}$-periodic solution $x^{\prime}(t)\in O.$ We may assume
that $O$ is so small that there is no $k^{\prime}$-periodic
solution in $O$ with $k^{\prime}\neq k$. Then, we have
$k^{\prime}=k, x^{\prime}\in S_k$. We assume that
$\widetilde{E}(s)\neq 0.$ From (3.29) we have
$$\int_{S^1}\widetilde{d}(\widetilde u(s_i^{\prime},t),\widetilde
u(s_i,t))dt \leq {1\over C}\left(\widetilde{E}(s_i)
\right)^{1/2},$$ where $\widetilde{d}$ denotes the distance
function defined by the metric $g_{\widetilde{J}}$ on $Z$. Taking
the limit $i\rightarrow \infty$, we get
$$\int_{S^1}\widetilde{d}(x^{\prime}(kt),x(kt))dt = 0.$$ It
follows that $$x^{\prime}(kt)= x(kt + \theta_0) $$ for some
constant $\theta_0$. This contradicts  (3.31). If there is some
$s_0$ such that $\widetilde{E}(s_0)=0$, then $|\Pi
\widetilde{u}_t|^2= |\Pi \widetilde{u}_s|^2 = 0 \;\;\;\forall
s\geq s_0.$ We still have a contradiction. \vskip 0.1in \noindent
We may choose a local coordinate system $y=(y_1,...,y_{2n} )$ on
$Z$ and a local trivialization of $\widetilde{M}\rightarrow Z$
around $x$ such that $x_k=\{ 0\leq \theta \leq 1, y= 0 \}$, and
$$\lambda = d\theta + \sum b_i(y) dy_i,\leqno(3.32)$$ where
$b_i(0)= 0$. Obviously, $\xi(\theta,0)$ is spanned by
$\frac{\partial} {\partial y_1},...,\frac{\partial}{\partial
y_{2n}}$. For $y$ small enough we may choose a frame
$e_1,...,e_{2n}$ for $\xi(\theta,y)$ as follows: in terms of the
coordinates $(\theta,y_1,...,y_{2n})$
$$e_i=(a_i(\theta,y),0,..,1,...,0),\;\; i=1,...,2n.$$ We write $$
u(s,t)=(a(s,t),\theta(s,t),y(s,t)).$$ Denote by $L$ the matrix of
the almost complex structure $\widetilde{J}$ on $\xi$ with respect
to the frame $e_1,...,e_{2n}$, and set $\widetilde{J}(s,t)
=L(u(s,t))$. We can write the equation (3.9) as follows:
$$a_s=\lambda(u_t)= \theta_t + \sum b_i(y)y_t \leqno(3.33)$$
$$a_t=-\lambda(u_s)=- \theta_s + \sum b_i(y)y_s \leqno(3.34)$$
$$y_s + \widetilde{J}(s,t)y_t=0. \leqno(3.35)$$ By the discussion
above, we have $y(s,t)\rightarrow 0$ in $C^{\infty}(S^1)$ as
$s\rightarrow \infty$. Hence we have $|y(s,t)|<1, |y_s(s,t)|<1,
|y_t(s,t)|<1$ for $s$ large enough. Then
$$\|y(s,t)\|^{2}_{L^2(S^1)}\leq \int_{S^1}|y(s,t)|dt \leq
\int_{S^1}\int_s^{\infty}|y_s(s,t)|dsdt $$ $$\leq
C_1\int_s^{\infty}\|\Pi \widetilde{u}_t\|_{L^2(S^1)}ds$$ for some
constant $C_1>0$. Together with (3.29) we have
$$\|y(s,t)\|_{L^2(S^1)}\leq {1\over
C}\widetilde{E}(s)^{1/4},\leqno(3.36)$$ $$\|\sum
b_i(y)y_s\|_{L^2(S^1)}\leq {1\over
C}\widetilde{E}(s)^{1/4},\leqno(3.37)$$ $$\|\sum
b_i(y)y_t\|_{L^2(S^1)}\leq {1\over
C}\widetilde{E}(s)^{1/4}.\leqno(3.38)$$ for some constant $C>0$
independent of $u(s,t)$. By using the same argument as in [HWZ1]
one can prove that there are constants $\delta>0$, $\ell_0$ and
$\theta_0$ such that for all $r=(r_1,r_2)\in {\bf Z^2 }$
$$|\partial^r[a(s,t)-ks-\ell_0]|\leq C_r
e^{-\delta|s|}\leqno(3.39)$$
$$|\partial^r[\theta(s,t)-kt-\theta_0]|\leq C_r e^{-\delta|s|}
\leqno(3.40)$$ $$|\partial^r w(s,t)|\leq C_r
e^{-\delta|s|}\leqno(3.41)$$ where $C_r$ are constants.
$\;\;\;\;\;\;\;\;\Box$ \vskip 0.1in \noindent We remark that since
$S_k$ is compact, we may choose $\delta, \;\;C_r$ to be
independent of $x_k \in S_k$ \vskip 0.1in \noindent Let $u:\Sigma
-\{p\}\rightarrow N$ be a $J$-holomorphic map with finite energy,
and $p$ be a nonremovable singularity. If , in terms of local
coordinates $(s,t)$ around $p$, $\lim_{s\rightarrow
\infty}\widetilde u(s,t)=x(kt)$, we say simply that $u(s,t)$
converges to the $k$-periodic orbit $x_k$. We call $p$ a {\em
positive} (resp. {\em negative} ) end, if $a(z) \rightarrow
\infty$ (resp. $-\infty$) as $z\rightarrow p$.

\subsection{\bf  Some technical lemmas }

Recall that if we collapse the $S^1$-action on $\widetilde{M}=H^{-1}(0)$
we obtain symplectic cuts $\overline{M}^{+}$ and $\overline{M}^{-}$. The
reduced space $Z$ is a codimension 2 symplectic submanifold of both
$\overline{M}^{+}$ and $\overline{M}^{-}$ (see [L]).
There is another way to look at this. The length of every orbit of the $S^1$ action on $\widetilde {M}$ with
respect to the metric $\langle\;,\;\rangle_{\omega_{\phi}}$ is $\phi^{\prime}$,
which converges to zero as $a \rightarrow \pm \infty $. Hence we can view
$\overline{M^{\pm}}$ as the completions of $M^{\pm}$. We will maily discuss
$M^+$, the discussion for $M^-$ is identical. Since $\overline{M}^{+}$
is a closed symplectic manifold, and we fix the values of $\phi$ on $[0,1]$,
the following lemma is well-known.
\vskip 0.1in
\noindent
{\bf Lemma 3.8} {\it There exist constants $\varepsilon_0>0$ and $C>0$
independent of $\phi$ such that the following holds, for each
$\varepsilon<\varepsilon_0$ and each metric ball $D_p(\varepsilon)$ centered
at $p \in M^{+}_0$ and of radius $\varepsilon$. Let $u^{\prime}:\Sigma^{\prime}
\rightarrow M^+$ be a
J-holomorphic map. Suppose $u^{\prime}(\Sigma^{\prime}) \subseteq
D_p(\varepsilon), u^{\prime}(\partial\Sigma^{\prime}) \subset
\partial D_p(\varepsilon),\partial \Sigma^{\prime} \neq \emptyset$ and
$p \in u^{\prime}(\Sigma^{\prime})$.Then
$$\int\limits_{\Sigma^{\prime}}u^{{\prime}*}\omega_{\phi}>C\varepsilon^2 .
\leqno(3.42)$$ }
\noindent
{\bf Lemma 3.9 } {\it There is a constant $\hbar_0>0$ such that
for every finite energy $J$-holomorphic map $u=(a,\widetilde{u}):
 \stackrel{\circ}{\Sigma}\rightarrow {\bf R}\times \widetilde{M} $ with
$\widetilde{E}(u)\neq 0 $ we have $\widetilde{E}(u)\geq \hbar_0$.}
\vskip 0.1in
\noindent
{\bf Proof:} Consider the $\widetilde{J}$-holomorphic map
$$\widetilde{v}=\pi\circ \widetilde{u}:\stackrel{\circ}{\Sigma} \rightarrow Z .$$
$\widetilde{v}$ extends to a $\widetilde{J}$-holomorphic curve
from $\Sigma $ to $Z$. Then the assertion follows from a standard result for
compact symplectic manifolds. $\;\;\;\Box$
\vskip 0.1in
\noindent
Let $(\Sigma; y_1,...,y_l, p_1,...,p_{\nu}) \in {\cal M}_{g,l+\nu}$,
and $u:\stackrel{\circ}{\Sigma} \rightarrow M^{+}$ be a
$J$-holomorphic map. Suppose that $u(z)$ converges to a $k_i$-periodic
orbit $x_{k_i}$ as $z$ tends to $p_i$.
By using the removable singularities theorem we get a $J$-holomorphic map
$\bar{u}$ from $\Sigma $ into $\overline{M}^{+}$. Let $A=[\bar{u}(\Sigma)]$.
It is obvious that
$$E_{\phi}(u)=\omega_{\phi}(A) \leqno(3.43)$$
which is independent of $\phi $.
We fix a homology class $A\in{H_{2}(\overline{M}^{+},\Z)}$ and a fixed set
$\{k_1,...,k_{\nu}\}$.
We denote by ${\cal{M}}_{A}(M^{+},g,l,{\bf k})$ the moduli space
of all $J$-holomorphic curves representing the homology class $A$ and converging
to ${\bf k}$-periodic orbits. There is a natural map
$$P: {\cal{M}}_{A}(M^{\pm},g,l,{\bf k})
\rightarrow \prod_{i=1}^{\nu} S_{k_i}$$
given by
$$P(u,\Sigma;{\bf y},\{p_1,...,p_{\nu}\})=
\left \{x_{k_1},...,x_{k_{\nu}}\right \}.\leqno(3.44)$$
We will identify $S_k$ with $Z$.
For $J$-holomorphic maps in
${\cal{M}}_{A}(M^{+},g,l,{\bf k})$, there is a uniform bound independent of
$\phi $ on the energy $E_{\phi}(u)$.
Similarly, for a map $u$ from $\Sigma$ into
${\bf R}\times \widetilde{M}$, we define
${\cal{M}}_{A}({\bf R}\times \widetilde{M},g,l,{\bf k}^+, {\bf k}^-)$ and
$P^{\pm}$ in the obvious way, where $\pm$ denote {\em positive} or {\em negative} ends.
For any $J$-holomorphic curves in
${\cal{M}}_{A}({\bf R}\times \widetilde{M},g,l,{\bf k}^+, {\bf k}^-)$, there
is a uniform bound on $\widetilde{E}$.

Following McDuff and Salamon [MS] we introduce the notion of
singular points for a sequence $u_{i}$ and the notion of mass of singular
points. Suppose that $({\Sigma}_{i};{\bf y}_i, {\bf p}_i)$ are stable curves
and converge to $(\Sigma;{\bf y}, {\bf p})$ in
${\overline{\cal{M}}}_{g,l+\nu}$. A point
$q\in{\Sigma-\{double \;points\}}$ is called regular for $u_{i}$ if there
exist $q_{i}\in{\Sigma}_{i}$, $q_{i}\rightarrow q$,
and $\epsilon>0$ such that the sequence $|du_{i}|_{h_{i}}$ is uniformly
bounded on $D_{q_{i}}(\epsilon,h_{i})$, where $|du_{i}|_{h_{i}}$ denotes the
norm with respect to the metric $\langle,\rangle $ on $N$ and the metric $h_{i}$ on
$\Sigma_{i}$. A point $q\in \Sigma -\{double \;points\}$ is called singular
for $u_{i}$ if it is not regular. A singular point $q$ for $u_{i}$ is called
rigid if it is singular for every subsequence of $u_{i}$. It is called tame
if it is isolated and the limit
$$m_{\epsilon}(q)=\lim_{i \rightarrow \infty}E_{\phi}(u_{i};D_{q_{i}}
(\epsilon,h_{i}))\leqno(3.45)$$
exists for every sufficiently small $\epsilon>0$.
The mass of the singular point q is defined to be
$$m(q)=\lim_{\epsilon\rightarrow 0}{m_{\epsilon}(q)}.\leqno(3.46)$$
\vskip 0.1in
\noindent
{\bf Lemma 3.10} {\it There is a constant $\hbar >0$ such that every rigid
singular point $y$ for $u_i$ has mass
$$m(y)\geq \hbar.\leqno(3.47)$$ }
\noindent
{\bf Proof: } By using the standard rescaling argument (see [MS]) we may
construct for every rigid singular point a nonconstant $J$-holomorphic map
$v:{\bf C}\rightarrow M^{+}(or \;\;{\bf R}\times \widetilde{M})$. In case
$ u({\bf C})\bigcap M^{\pm}_0\neq \emptyset $ we use Lemma 3.9. Now we assume
that $u({\bf C})$ lies in the cylindrical part. By Lemma 3.5
$\widetilde{E}(v)\neq 0$. We have $\widetilde{E}(v)\geq \hbar_0$
by lemma 3.9. Then Lemma 3.10 follows by (3.6) and (3.13). $\;\;\;\Box$
\vskip 0.1in
\noindent
Denote by $P\subset \Sigma $ the set of  singular points for $u_i$,
the double points and the puncture points. By Lemma 3.10 and (3.43), $P$ is
a finite set. By definition, $|du_i|_{h_i}$ is uniformly bounded on every
compact subset of $\Sigma - P$. We assume that
$$u_i(\Sigma_i)\cap M^+_l \neq \emptyset \leqno(3.48)$$
for some positive number $d$ independent of $i$, where
$M^{+}_d =  M^{+}_{0}\bigcup\{[0, d]\times \widetilde{M}\}.$
If there is no $d$ satisfying (3.48), we will make a translation
along ${\bf R}$; then the discussion is similar.
It follows from (3.48) and the uniformly boundness of $|du_i|_{h_i}$ on compact
sets that $u_i$ maps every compact subset of $\Sigma - P$
uniformly into a bounded subset of $M^+$ (recall the definition (3.8) of the metric
$\langle, \rangle$). By passing to a further subsequence we may
assume that $u_i$ converges uniformly  with all derivatives on every compact
subset of $\Sigma - P $ to a $J$-holomorphic map $u:\Sigma - P
\rightarrow M^+.$ Obviously, $u$ is a finite energy $J$-holomorphic map.

We need to study the behaviour of the sequence $u_{i}$ near each
singular point for $u_{i}$ and near each double point. We consider
singular points for $u_i$. The same method with only minor changes applies to
double points. Let $q \in \Sigma $ be a rigid singular point for
$u_{i}$. Then $q$ is a puncture point for $u$. If $q$ is a removable singularity
of $u$, i.e., there are $\epsilon > 0, d>0 $ such that $u(D_q(\epsilon))\subset
M^+_d$, we consider the compact manifold $\overline{M}^{+}$, and
construct bubbles as usual for a compact symplectic manifold (see [RT1], [PW],
[MS]). In the following, we assume that $q$ is a nonremovable singularity.
We may identify each neighborhood of $q_{i}$ in
$\Sigma_{i}$ with a neighborhood of 0 in {\bf C}, say $D(1)$. The sequence
$u_{i}$ is considered to be a sequence of $J$-holomorphic maps from $D(1)$
into $M^+$, and $q_{i}\in D(1)$, $q_{i}\rightarrow 0$. Without loss of generality,
we may assume that 0 is the unique singular point and is tame (see [MS]).
For $|z|$ sufficiently small, $u(z)$ lies in the cylindrical part.
It is convenient to use cylindrical coordinates $z=e^{s+2 \pi it}.$ We write
$$u_i(s,t)=(a_i(s,t), \widetilde{u}_i(s,t))=(a_i(s,t), \theta_i(s,t), {\bf x}_i(s,t))$$
$$u(s,t)=(a(s,t),\widetilde{u}(s,t))=(a(s,t), \theta(s,t), {\bf x}(s,t)).$$
Let $S^{2}$ denote the standard sphere  in ${\bf R^3}$ with two
distinguished points $0=(0,0,0)$ and its antipodal north pole $(0,0,1)$. For any $\epsilon >0$,
through a conformal transformation
of $S^2$ we may assume that $u_i$ are defined in the
disk $D(\epsilon)=\{z | |z| < \epsilon \}$ such that the center of mass
of the measure $|du_i|^2 $ is on the $z$-axis.
By Theorem 3.7, we have
$$\lim_{s \rightarrow -\infty} \widetilde{u}(s,t)=x(kt)$$
in $C^{\infty}(S^1)$, where $x(\;,\; )$ is a $k$-periodic orbit on $\widetilde{M}$.
Choosing $\epsilon $ small enough we have $|\frac{\partial \theta}{\partial t} - k |< 1/10,\;\;
|\frac{\partial a}{\partial s} - k |< 1/10$ for all $(s,t)\in (-\infty, \log\epsilon)
\times S^{1}$. By choosing $i$ large enough we may assume that
$|\frac{\partial \theta_i}{\partial t}(\log\epsilon,t) - k |< 1/5, \;\;
|\frac{\partial a_i}{\partial s}(\log\epsilon,t) - k |< 1/5$ for all $i$ and all $t \in S^1 $.
For every $i$ there exists $\delta^{\prime}_i>0$ such that
$$E_{\phi}(u_i; D(\delta^{\prime}_i))=m_0 - \frac{1}{2}\hbar_0.$$
By definition of the mass $m_0$, the sequence $\delta^{\prime}_i$ converges to $0$.
We discuss two cases: \\
{\bf Case 1}: There is an $L>0$ such that for all $i>L$
$$|\frac{\partial \theta_i}{\partial t}(s,t) - k |\leq 1/2, \;\;
|\frac{\partial a_i}{\partial s}(s,t) - k |\leq 1/2 \leqno(3.49)$$
hold for all $\log\delta^{\prime}_i\leq s \leq\log\epsilon$, $t \in S^1 $; we set
$\delta_i=\delta^{\prime}_i$. \\
{\bf Case 2}: Otherwise, we could find a subsequence, still denoted by $u_i$,
and $\delta_i > \delta^{\prime}_i$ such that for all $i$ (3.49) holds for all
$\log\delta_i\leq s \leq\log\epsilon , \;\;t \in S^1 $,
and there is $t_0\in S^1$ such that
$$|\frac{\partial \theta_i}{\partial t}(\delta_i,t_0) - k |= 1/2, \;\; or \;\;
|\frac{\partial a_i}{\partial s}(\delta_i,t_0) - k |= 1/2. \leqno(3.50)$$
Since $u_i$ converges uniformly with all derivatives to $u$ on any compact set
of $D(\epsilon)-\{0\}$, $\delta_i $ must converge to 0.
Put $z=\delta_i w = \delta_i e^{r+ 2\pi it}$, and
define the $J$-holomorphic curve $v_i(r,t)$ by
$$v_i(r,t)=(b_i(r,t),\widetilde{v}_i(r,t))=
\left(a_i(\log\delta_{i}+ r,t) -a_i(\log\delta_i,t_0) ,
\widetilde{u}_i(\log\delta_{i}+ r,t)\right).\leqno(3.51)$$
where $t_0\in S^1$ is a fixed point.
\vskip 0.1in
\noindent
{\bf Lemma 3.11 } {\it Suppose that 0 is a nonremovable singular point of $u$.
Define the $J$-holomorphic map $v_i$ as above. Then there exists a
subsequence (still denoted by $v_i$) such that
\begin{description}
\item[(1)] The set of singular points $\{w_1,\cdot \cdot \cdot,w_N\}$ for
$v_i$ is
finite and tame, and is contained in the disc $D(1)=\{w \mid \mid w
\mid \leq 1\};$
\item[(2)] The subsequence $v_i$ converges with all derivatives uniformly on
every compact subset of \linebreak ${\bf C}\backslash\{w_1,\cdot \cdot \cdot,w_N\}$
to a nonconstant
(J,i)-holomorphic map $v:{\bf C}\backslash\{w_1,\cdot \cdot \cdot,w_N\}
\rightarrow {\bf R}\times \widetilde{M};$
\item[(3)]$\lim_{s \rightarrow -\infty} \widetilde{u}(s,t)=
\lim_{r \rightarrow \infty} \widetilde{v}(r,t);$
\item[(4)] $\widetilde{E}(v)+\sum\limits_{1}^{N}m(w_i)=m_0. $
\end{description}}
\vskip 0.1in
\noindent
{\bf Proof: }
{\bf (1) } By definition of $m_0$ and $\delta_i$ we have
$$\lim_{i \rightarrow \infty} \sup \widetilde{E}(v_{i}; D(R)-D(r)) \leq
\frac{\hbar_0}{2}$$
for any $R>r>1$. We claim that $|dv_i|$ is uniformly bounded on ${\bf C}-D(r)$.
Otherwise, by a standard rescaling argument (see [MS]) we could
construct a nonconstant $J$-holomorphic curve
$v:{\bf C}\rightarrow {\bf R}\times \widetilde{M}$ with
$\widetilde{E}(v; {\bf C}) \leq\frac{\hbar_0}{2}$. It follows from Lemma 3.9 that
$\widetilde{E}(v; {\bf C})=0$. Then $v$ must be constant by Lemma (3.5), a
contradiction.
\vskip 0.1in
\noindent
{\bf (2) } The proof is standard (see [MS]).
\vskip 0.1in
\noindent
{\bf (3) } Consider the $\widetilde{J}$-holomorphic map
$$f_i=\pi \widetilde{u_i}:\Sigma \rightarrow Z .$$
Write $A(r,R)=D(R)-D(r)$. Since $E(f_i, A(R\delta_i, \epsilon )) \leq
\frac{1}{2}\hbar_0$, by a standard result about the energy of a pseudo-holomorphic curve on an
arbitrarily long cylinder (see [MS]), there exists a $T_0>0$ such that for
$T>T_0$
$$E(f_i; A(R\delta_i e^T, \epsilon e^{- T})) \leq
\frac{C}{T}E(f_i; A(R\delta_i, \epsilon))\leqno(3.52)$$
and
$$\int\limits_{S^1}d(f_i(\epsilon e^{- T+it}),
f_i(R\delta_ie^{T+it}))dt\leq\frac{C}{\sqrt{T}}\sqrt{E(f_i;A(R\delta_i,\epsilon))}.\leqno(3.53)$$
We choose $T > C$. It follows from (3.52) that
$$\widetilde{E}(u_i; A(R\delta_i ,\epsilon )) \leq
\frac {1}{1 - \frac {C}{T}}\left(\widetilde{E}(u_i; A(\epsilon e^{- T}, \epsilon))
+ \widetilde{E}(v_i;A(R , Re^T))\right).\leqno(3.54)$$
Since $u_i \rightarrow u $ and $v_i \rightarrow v $ uniformly on compact sets,
the above inequality implies that
$$\lim_{\epsilon \rightarrow 0, R\rightarrow \infty}lim_{i\rightarrow
\infty}\widetilde{E}(u_i; A(\epsilon, R\delta_i))=0.\leqno(3.55)$$
It is easy to see that $v$ is a finite energy $J$-holomorphic curve.
Suppose that $v$ converges to a $k^{\prime}$-periodic orbit $x^{\prime}$. From (3.53) and (3.55), we obtain
$$\int\limits_{S^1}\widetilde{d}(x(kt),x^{\prime}(kt))dt=0.$$
 From (3.49), it follows that $|k-k^{\prime}|\leq 1/2$, and therefore $k=k^{\prime}$.
\vskip 0.1in \noindent {\bf (4) } For $\epsilon >0$ sufficiently
small and $R$ large enough, $$E_{\phi}(u_i, D(\epsilon))=
E_{\phi}(u_i, D(\epsilon)- D(R\delta_i)) + E_{\phi}(v_i, D(R)).$$
Obviously, $$\lim_{R\rightarrow
\infty}\lim_{i\rightarrow\infty}E_{\phi}(v_i, D(R))=
\widetilde{E}(v)+\sum\limits_{1}^{N}m(w_i).$$ By (3.49),
$$\lim_{\epsilon \rightarrow 0}\lim_{i\rightarrow\infty}\min
\{a_i(s,t)|_{A(\delta_i,\epsilon)}\}\rightarrow \infty.$$ Together
with (3.55), we have $$\lim_{\varepsilon \rightarrow 0,
R\rightarrow \infty}\lim_{i\rightarrow \infty}E_{\phi}(u_i,
D(\epsilon)- D(R\delta_i))=0.$$ Then (4) follows. $\Box$ \vskip
0.1in \noindent It is possible for $v$ to have
$\widetilde{E}(v)=0$. As in the case of compact manifolds (see
[PW]) we call $v$ a ghost bubble. \vskip 0.1in \noindent {\bf
Lemma 3.12} {\it If $\widetilde{E}(v)=0$ then $N \geq 2.$ } \vskip
0.1in \noindent {\bf Proof: } For Case 1, we use the same argument
as in [PW]. We now prove our lemma for Case 2. If $0$ is the only
singularity, then $v$ is a $J$-holomorphic map from ${\bf R}\times
S^1$ to ${\bf R}\times\widetilde{M}$ with $\widetilde{E}(v)=0$. It
must be $(kr + d, kt + \theta_0)$ by Lemma 3.5. But from (3.50),
it follows that $$|\frac{\partial \theta}{\partial t}(0,t_0) - k
|= 1/2, \;\; or \;\; |\frac{\partial b}{\partial r}(0,t_0) - k |=
1/2.$$ We get a contradiction. Hence there is a singular point
$z\neq 0$. The lemma follows from  symmetry.$\Box$. \vskip 0.1in
\noindent {\bf Remark 3.13 }{\it We need also study sequences of
$J$-holomorphic curves in ${\bf R}\times\widetilde{M}$. The method
is the same as we did for $M^+$. In particular, if $q$ is a
removable singular point of $u$, we collapse the $S^1$ action on
$\widetilde{M}$ at $\infty$ and $-\infty$ to obtain a compact
manifold and construct a bubble as usual procedure for compact
manifold (see [RT1],[MS],[PW]). If $q$ is a nonremovable singular
point of $u$, we construct bubble as above.}
 \vskip 0.1in
\noindent For any $r>0$, we construct $M_r$ as follows: The parts
$|a^{\pm}| \geq 3r $ are cut out from $M^{\pm}$, and the
remainders are glued along a collar of length 2r by the formula
$$a^{-}=a^{+}-4r \leqno(3.56)$$ Let $r\rightarrow \infty$. We may
consider the limit $M_{\infty}$ of $M_r$ to be a compact manifold
obtained by gluing $M^{+}$ and $M^{-}$ along $\widetilde{M}$ at
$+\infty$ and $-\infty$. Let $u_r: \Sigma_r\rightarrow M_r$ be a
$J$-holomorphic curve. Suppose that $\Sigma_r \rightarrow \Sigma =
\Sigma^{+}\cup \Sigma^{-}$. Both $\Sigma^{+}$ and $\Sigma^{-}$ may
contain several connected nodel Riemann surfaces, and $\Sigma^{+}$
and $\Sigma^{-}$ may have several intersect points. For
simplicity, we assume that both $\Sigma^{+}$ and $\Sigma^{-}$ have
one component and $\Sigma^{+}$ and $\Sigma^{-}$ intersect at one
point $p$. For the general case the treatment is identical. We may
consider $\Sigma_r$ to be obtained from $\Sigma$ by resolving the
singularity using the parameter $\delta_r$, $\delta_r \rightarrow
0,$ i.e., $\Sigma_r = \Sigma_1 \#_{\sqrt{\delta_r}}\Sigma_2$ with
gluing formulas $t=t^{\prime}$, $s^{\prime} + \log\delta_r = s.$
The same method as we used before, with only minor changes ,
applies.

\subsection{\bf Compactness theorems }

We introduce a terminology. Let $\Sigma_1$ and $\Sigma_2$
join at $p$, and $(u_1,u_2): \Sigma_1\wedge\Sigma_2 \rightarrow
{\bf R}\times \widetilde{M}$ a map. Choose holomorphic cylindrical
coordinates $z_1=(s_1,t_1)$ on $\Sigma_1$ and $z_2=(s_2,t_2)$
on $\Sigma_2$ near $p$ respectively. Suppose that
$$\lim_{s_1 \rightarrow -\infty} \widetilde{u_1}(s_1,t_1)= x_1(k_1t_1)$$
$$\lim_{s_2 \rightarrow +\infty} \widetilde{u_2}(s_2,t_2)= x_2(k_2t_2).$$
We say $u_1$ and $u_2$ converge to a same periodic orbit, if $k_1=k_2$,
and $P(x_1)=P(x_2)$, where $P$ denotes the projection to $Z$.
 \vskip 0.1in
 \noindent
{\bf Definition 3.14 } {\it Let $(\stackrel{\circ}{\Sigma};{\bf
y},{\bf p})$ be a Riemann surface of genus g with $l$ marked
points ${\bf y}$ and $\nu$ ends ${\bf p}.$ A log stable
holomorphic map with $\{k_1,...,k_{\nu}\}$-ends from
$(\stackrel{\circ}{\Sigma};{\bf y},{\bf p})$ into $M^{+}$ is an
equivalence class of continuous maps u from
$\stackrel{\circ}{\Sigma}^{\prime}$ into $(M^{+})^{\prime}$,
modulo the automorphism group $stab_u$, the translations, and
$S^1$-action on ${\bf R}\times \widetilde{M}$, where
$\stackrel{\circ}{\Sigma}^{\prime}$ is obtained by joining chains
of ${\bf P^1}s$ at some double points of $\Sigma$ to separate the
two components, and then attaching some trees of ${\bf P^1}$'s;
$(M^{+})^{\prime}$ is obtained by attaching some ${\bf R}\times
\widetilde{M}$ to $M^{+}$.  We call components of
$\stackrel{\circ}{\Sigma}$ principal components and others bubble
components. Furthermore,
\begin{description}
\item[(1)] If we attach a tree of ${\bf P^1}$ at a marked point $y_i$ or a
puncture point $p_i$, then $y_i$ or $p_i$ will be replaced by a point
different from the intersection points on a component of the tree. Otherwise,
the marked points or puncture points do not change;
\item[(2)]  $\stackrel{\circ}{\Sigma'} $ is a connected curve with
normal crossings ;
\item[(3)] Let $m_{j}$ be the number of points on ${\Sigma}_{j}$ which are
nodal points or  marked points or puncture points. If
$u|_{{\Sigma}_{j}}$ is constant or $\pi u|_{\Sigma_j}$ is constant
(for maps into ${\bf R}\times \widetilde{M}$), then $m_{j}+2g_{j}\geq3$;
\item[(4)] The restriction of u to each component is J-holomorphic.
\item[(5)]  $u$ converges exponentially to some periodic orbits
$(x_{k_1},...,x_{k_{\nu}})$ as the variable tends to the puncture
$(p_1,...,p_{\nu})$; more precisely, $u$ satisfies (3.39)-(3.41);
\item[(6)] Let $q$ be a nodal point of $\Sigma^{\prime}$ . Suppose $q$ is
the intersection point of ${\Sigma}_{i}$ and ${\Sigma}_{j}$. If $q$ is a
removable singular point of $u$, then u is continuous at $q$;
If $q$ is a nonremovable singular point of $u$, then
$u|_{\Sigma_{i}}$  and $u|_{\Sigma_{j}}$ converge exponentially to the same
periodic orbit on $\widetilde{M}$ as the variables tend to the nodal point $q$.
\end{description}}
\vskip 0.1in If we drop the condition (4), we simply call $u$ a
log stable map. Put $m=l+\nu$ and let $T_m$ be as in introduction.
Let ${\overline{\cal{M}}}_{A}(M^{+},g,T_m)$ be
the space of equivalence classes of log stable holomorphic maps
with ends, and ${\overline{\cal{B}}}_{A}(M^{+},g,T_m)$ be the
space of stable maps with ends.
$\overline{\cal{M}}_{A}(M^{+},g,T_m)$ has an obvious
stratification indexed by the combinatorial type of the domain
with the following data:
\begin{description}
\item[(1)] the topological type of the domain as an abstract Riemann surface
with marked points and puncture points;
\item[(2)] sets of periodic orbits corresponding to each puncture point;
\item[(3)] a decomposition of the integral 2-dimensional class $A = \sum A_i$.
\end{description}
Suppose that ${\cal D}^{J,A}_{g,T_m}$ is the set of indices.
\vskip 0.1in \noindent {\bf Lemma 3.15 }{\it ${\cal
D}^{J,A}_{g,T_m}$ is a finite set.} \vskip 0.1in \noindent {\bf
Proof: } By Lemma 3.10, 3.12 and a standard argument, one can show
that there are finitely many combinatorial types of the domain as
an abstract nodal surface and integral 2-dimensional classes
$A_i$. The class $A$ determines the bound on $T_m$.
Since for a $J$-holomorphic curve in
$${\cal{M}}_{A_i}({\bf R}\times\widetilde{M},g_i,m_i,{\bf k^+}, {\bf k^-})$$
the difference
$$\sum_j k^{+}_j - \sum_j k^{-}_j $$ is determined by the class $A_i$, we
conclude that the set of periods is finite. $\Box$

By using the above lemmas one immediately obtains \vskip 0.1in
\noindent {\bf Theorem 3.16 } {\it Let $\Gamma_i=(u_i,\Sigma_i;
{\bf y}_i,{\bf p}_i) \in {\cal M}_{A}(M^{+},g,T_m)$ be a
sequence. Then there is a subsequence which  {\em weakly
converges} to a log stable $J$-holomorphic curve in
${\overline{\cal{M}}}_{A}(M^{+},g,T_m)$. Here, by  weak
convergence, we mean the Gromov-Uhlenbeck convergence with
possible translation and $S^1$-action on ${\bf R}\times
\widetilde{M}$.}

 \vskip 0.1in
 \noindent
 {\bf Corollary 3.17 } {\it
${\overline{\cal{M}}}_{A}(M^{+},g,T_m)$ is compact. }
\vskip
0.1in
\noindent
 To define stable $J$-holomorphic maps in
$M_{\infty}$ we need to extend log stable maps into $M^{\pm}$ to
include nonconnected components.
Suppose that $\Sigma^{\pm}$ has $l^{\pm}$ connected components
$\Sigma^{\pm}_i,i=1,...,l^{\pm}$ of genus $g^{\pm}_i$ with
$m^{\pm}_i$ marked points and ${\nu}^{\pm}_i$ ends. We have
$$\sum{\nu}^{+}_i= \sum{\nu}^{-}_i=\nu,\;\;\;\;\;\sum{m}^{+}_i = m^+,
\;\;\;\;\sum{m}^{-}_i=m^-,\;\;\;\;\; m^+ + m^- = m.$$
Put
$$\overline{\cal{M}}_{{\bf
A}^{\pm}}(M^{\pm},{\bf g}^{\pm},T_{{\bf m}^{\pm}})=
\bigoplus_{i=1}^{l^{\pm}}{\overline{\cal{M}}}_{A^{\pm}_i}(M^{\pm},g^{\pm}_i,
T_{m^{\pm}_i}),$$ where ${\bf
A}^{\pm}=\left\{A^{\pm}_1,...,A^{\pm}_{l^{\pm}}\right\}, {\bf
g}^{\pm}=\left\{g^{\pm}_1,...,g^{\pm}_{l^{\pm}}\right\}, {\bf
m}^{\pm}=\left\{m^{\pm}_1,...,m^{\pm}_{l^{\pm}}\right\}, {\bf
k}=\left\{{\bf k}_1,...,{\bf k}_{l^{\pm}}\right\}.$ The maps $e_j^{\pm}$
are extended to ${\overline{\cal{M}}}_{{\bf A}^{\pm}}(M^{\pm},{\bf
g}^{\pm},T_{m^{\pm}})$ in a natural way.
 \vskip 0.1in
\noindent
{\bf Definition 3.18 } {\it A stable $J$-holomorphic map
of genus $g$ and class $A$ into $M_{\infty}$ is a triple
$(\Gamma^{-}, \Gamma^{+}, \rho)$, where
$\Gamma^{\pm}\in \overline{\cal{M}}_{{\bf A}^{\pm}}(M^{\pm},{\bf g}^{\pm},T_{m^+})$
and $\rho:\{p^+_1,...,p^+_{\nu}\}\rightarrow
\{p^-_1,...,p^-_{\nu}\} $ is a one-to-one map satisfying
\begin{description}
\item[(1)] If we identify $p^{+}_i$ and $\rho(p^{+}_i)$ then $\Sigma^{+}\bigcup
\Sigma^{-}$ forms a connected closed Riemann surface of genus $g$;
\item[(2)] $\widetilde{u}^+(z)$ and $\widetilde{u}^-(w)$ converge to the same
periodic orbit when $z\rightarrow p^+_i$ and $w \rightarrow \rho(p^{+}_i)$
respectively;
\item[(3)] $(\Gamma^{-}, \Gamma^{+}, \rho)$ represents the homology class $A$.
\end{description}}
\vskip 0.1in \noindent {\bf Remark 3.19 } {\it If
$\{p^{\pm}_1,...,p^{\pm}_{\nu}\}$ in the above definition is the
empty set then $(\Sigma,u)$ is a log stable $J$-holomorphic map
into $M^+$ or $M^-$. } \vskip 0.1in Denote by
${\overline{\cal{M}}}_{A}(M_{\infty},g,m)$ the moduli space of
stable $J$-holomorphic maps in $M_{\infty}$.
Using Lemmas in subsection 3.2 we immediately obtain the following convergence
theorem:
\vskip 0.1in \noindent
{\bf Theorem 3.20 } {\it
${\overline{\cal{M}}}_{A}(M_{\infty},g,m)$ is compact. } \vskip
0.1in \noindent {\bf Theorem 3.21 } {\it Let $\Gamma_{r} \in
{\overline{\cal{M}}}_{A} (M_{r},g,m)$ be a sequence. Then there is
a subsequence which {\em weakly converges} to a stable
$J$-holomorphic map in
${\overline{\cal{M}}}_{A}(M_{\infty},g,m)$.} \vskip 0.1in
\noindent

\section{\bf  Log invariants}

In this section we define the log invariants for a pair $(V,B)$, where $V$
is a compact symplectic manifold, and $B$ a codimension two symplectic
submanifold of $V$. Since the symplectic structure of a tubular neighborhood
of $B$ is modeled on a neighborhood of $Z$ in $\overline{M}^+$. So we
consider $(\overline{M}^+, Z)$.
To define the log GW-invariants we use the virtual neighborhood
technique developed in \cite{R5}. We remark that other
constructions \cite{FO}, \cite{LT3}, \cite{S} can be applied to
our case as well.

\subsection{\bf  Stabilization equations }

We will describe $\overline{\cal B}_{A}(M^+,g,T_m),$
$\overline{\cal F}_{A}(M^{+},g,T_m)$, then we construct a finite
dimensional $V$-bundle ${\cal E}$ over $\overline{\cal
B}_{A}(M^+,g,T_m)$, following Siebert's construction, and define a
stabilization equation ${\cal S}_e=0$. Note that we need only to
consider a neighborhood ${\cal U}$ of $\overline{\cal{M}}_{A}
(M^{+},g,T_m)$ in the configuration space $\overline{\cal
B}_A(M^{+},g,T_m)$ of $C^{\infty}$-log stable (holomorphic or not)
maps. We first consider the case that $\Sigma$ is a smooth
component. Denote by $N$ the one of $M^+$ and ${\bf R}\times
\widetilde{M}$. We consider ${\cal B}_{A_i}(N,g_i,T_{m_i})$. Let
$b = (u,(\Sigma,j);{\bf p}) \in {\cal B}_{A_i}(N,g_i,T_{m_i})$.
Here $j$ is a complex structure ( including marked points ), which
is standard near each puncture point. We introduce the holomorphic
cylindrical coordinates on $\Sigma$ near each puncture points
$p_i$. Then we may consider $(\Sigma,j)$ as a Riemann surface with
some cylindrical ends. Note that we will consider a nodel point of
two Riemann surfaces as a puncture point of both Riemann surface
and use the cylindrical model. So there are two types of puncture
points: orbit type puncture points and non-orbit puncture points,
which must be nodel points. By using the removable singularities
theorem we may consider $u$ as a $J$-holomorphic map from
$(\Sigma,j)$ into $\overline{M}^+$ or $\Re$, where we denote by
$\Re$ the space obtained from ${\bf R}\times \widetilde{M}$ by
collapsing the $S^1$-action on the $\pm \infty$ ends. Then for any
non-orbit puncture point $p$, $u(p)$ lies in  $\overline{M}^+ - Z$
( or $\Re-Z$), while for any orbit type puncture point $p$, $u(p)$
lies in $Z$. For each puncture point $p_i$, we choose the
holomorphic cylindrical coordinates $(s,t)$ on $\Sigma$ near each
$p_i$. Over the tube the linearized operator
$$D_{u}=D\bar{\partial}_{J}(u):C^{\infty}
(\Sigma;u^{\ast}TN)\rightarrow \Omega^{0,1}(u^{\ast}TN)$$ takes
the following form $$D_{u}=\frac{\partial}{\partial
s}+J_0\frac{\partial} {\partial t}+S = \bar{\partial}_{J} +
S.\leqno(4.1)$$ For each orbit type puncture point we choose a
local frame field $e_1=\frac{\partial}{\partial a}, e_2=X_H$ and
$e_3,...,e_{2n+2}\in \xi$ near the periodic orbit $x_i(k_it)$.
Denote by $\bigtriangledown$ the Levi-Civita connection with
respect to the metric $\langle \;,\;\rangle$ . Then
$$S_{k\ell}=\left<e_k,\bigtriangledown_se_{\ell} +
J\bigtriangledown_te_{\ell}+(\bigtriangledown_{e_{\ell}}J)
\frac{\partial u}{\partial t}\right>.\leqno(4.2)$$ Recall that for
a $J$-holomorphic map $u$, $D_u$ is independent of the choice of
connection. Note that $S_{1\ell}=S_{\ell1}=S_{2\ell}=S_{\ell2}=0$
. Since $\frac{\partial \widetilde{u}}{\partial s}\rightarrow 0$,
$\frac{\partial u}{\partial t}\rightarrow X_H$ exponentially and
uniformly in $t$ as $s\rightarrow \pm \infty $, we have $$|S|\leq
Ce^{-\delta s}\leqno(4.3)$$ for some constant $C>0$ for $s$ big
enough. Therefore, the operator $H_s=J_0\frac{d}{dt}+S$ converges
to $H_{\infty}=J_0\frac{d}{dt}$. Obviously, the operator $D_u$ is
not a Fredholm operator because over each orbit puncture end the
operator $H_{\infty}=J\frac{D}{dt}$ has zero eigenvalue. The $\ker
H_{\infty}$ consists of constant vectors. This is also true for
each non-orbit puncture end. To recover a Fredholm theory we use
weighted function spaces. We choose a weight $\alpha$ for each
end. Fix a positive function $W$ on $\Sigma$ which has order equal
to $e^{\alpha |s|} $ on each end, where $\alpha$ is a small
constant such that $0<\alpha<\delta $ and over each end
$H_{\infty}- \alpha = J_0\frac{d}{dt}- \alpha $ is invertible. We
will write the weight function simply as $e^{\alpha |s|}.$ For any
section $h\in C^{\infty}(\Sigma;u^{\ast}TN)$ and section $\eta \in
\Omega^{0,1}(u^{\ast}TN)$ we define the norms
$$\|h\|_{1,p,\alpha}=\left(\int_{\Sigma}(|h|^p+ |\nabla
h|^p)d\mu\right)^{1/p}
+\left(\int_{\Sigma}e^{2\alpha|s|}(|h|^2+|\nabla h|^2)
d\mu\right)^{1/2} \leqno(4.4)$$
$$\|\eta\|_{p,\alpha}=\left(\int_{\Sigma}|\eta|^p d\mu
\right)^{1/p}+ \left(\int_{\Sigma}e^{2\alpha|s|}|\eta|^2
d\mu\right)^{1/2} \leqno(4.5)$$ for $p\geq 2$, where all norms and
covariant derivatives are taken with respect to the  metric
$\langle\;\;\rangle$ on $u^{\ast}TN$ defined in (3.8), (3.9), and
the metric on $\Sigma$. Denote $${\cal C}(\Sigma;u^{\ast}TN)=\{h
\in C^{\infty}(\Sigma;u^{\ast}TN); \|h\|_{1,p,\alpha}< \infty
\},\leqno(4.6)$$ $${\cal C}(u^{\ast}TN\otimes \wedge^{0,1})
=\{\eta\in \Omega^{0,1}(u^{\ast}TN); \|\eta\|_{p,\alpha}< \infty
\}.\leqno(4.7)$$ Denote by $W^{1,p,\alpha}(\Sigma;u^{\ast}TN)$ and
$L^{p,\alpha}(u^{\ast}TN\otimes \wedge^{0,1})$ the completions of
${\cal C}(\Sigma;u^{\ast}TN)$ and ${\cal C}(u^{\ast}TN\otimes
\wedge^{0,1}) $ with respect to the norms (4.4) and (4.5)
respectively.

For each puncture point $p_i,i=1,...,\nu $,let $h_{i0}\in \ker H_{i\infty}$. Put
$H_{\infty}=(H_{1 \infty},...,H_{\nu \infty})$, $h_{0}=(h_{1 0},...,h_{\nu 0}).$
We choose a normal coordinate system near each non-orbit type puncture, and
choose a Darboux coordinate system near each orbit type puncture.
$h_0$ may be considered as a vector field in the coordinate neighborhood.
We fix a cutoff function $\varrho$:
\[
\varrho(s)=\left\{
\begin{array}{ll}
1, & if\ |s|\geq d, \\
0, & if\ |s|\leq \frac{d}{2}
\end{array}
\right.
\]
where $d$ is a large positive number. Put
$$\hat{h}_0=\varrho h_0.$$
Then for $d$ big enough $\hat{h}_0$ is a section in $C^{\infty}(\Sigma; u^{\ast}TN)$
supported in the tube $\{(s,t)||s|\geq \frac{d}{2}, t \in S^1 \}$.
Denote
$${\cal W}^{1,p,\alpha}=\{h+\hat{h}_0 | h \in
W^{1,p,\alpha},h_0 \in kerL_{\infty}\}.$$
The operator $D_u:{\cal W}^{1,p,\alpha}\rightarrow L^{p,\alpha}$
is a Fredholm operator so long as $\alpha$ does not lie in the spectrum
of the operator $H_{i \infty}$ for all $i=1,\cdots,\nu$.
\vskip 0.1in
\noindent
{\bf Remark 4.1 }{\it The index $ind (D_u,\alpha)$ does not change if
$\alpha$ is varied in such a way that $\alpha$ avoids the spectrum of
$L_{\infty i}.$ Conversely, the index will change if $\alpha$ is moved
across an eigenvalue. We will choose $\alpha$ slightly larger than zero such
that at each end it does not across the first positive eigenvalue.}
\vskip 0.1in
\noindent
The proof of the following lemma is almost the same as in \cite{R5}, we omit it.
\vskip 0.1in
\noindent
{\bf Lemma 4.2 } ${\cal B}_{A_i}(N,g_i,T_{m_i})$ is a Hausdorff Frechet V-manifold
for any $2g_i+m_i\geq 3$ or $g_i=0,m_i<3,A_i\neq 0.$
\vskip 0.1in
\noindent
{\bf Remark 4.3} Let $b=(u,\Sigma,j;{\bf p})
\in {\cal{M}}_{A_i}(N,g_i,T_{m_i})$. We describe the neighborhoods of $b$.
\vskip 0.1in
\noindent
{\bf 1: $N=M^+$.} When $2g_i+m_i\geq 3$, $\Sigma$ is stable
and ${\cal M}_{g_i,m_i}$ is a
V-manifold. Hence, the automorphism group $Aut_{\Sigma}$ of $\Sigma$ is finite.
Denote by $O_{j}$ a neighborhood of complex structures on $(\Sigma,j)$.
Note that we change the complex structure in a compact set
$K_{\pm,deform}$ of $\Sigma$ away from the puncture points.
A neighborhood ${\cal U}_b$ of $b$ can be described as
$$O_{j}\times \{\exp_u(h + \hat{h}_0); h\in {\cal C}(\Sigma;u^{\ast}TM^{\pm}),
h_0\in \ker H_{\infty}, \|h\|_{1,p,\alpha} + |h_0| <\epsilon \}/stab_u.$$
For the case $g_i=0, m_i\leq 2, A_i\neq 0$, $\Sigma$ is no longer stable and
the automorphism group $Aut_{\Sigma}$ is infinite.
One must construct a slice $W_u$ of the action $Aut_{\Sigma}$ such that
$W_u/stab_u $ is a neighborhood of $(\Sigma,u)$. We can write again
${\cal U}_b=O_{j}\times W_u/stab_u $ with $O_{j}=point$
( see \cite{R5}).
\vskip 0.1in
\noindent
{\bf 2: $N={\bf R}\times \widetilde{M}$.} We must mod the group
$\C^*$ generated by the $S^1$-action and the translation along ${\bf R}$.
We fix a point $y_0\in \Sigma $ different from the marked
points and puncture points. Fix a local coordinate system $a,\theta,w$ on
${\bf R}\times\widetilde{M}$ such that $u(y_0)=(0,0,0)$. We use
${\cal C}(\Sigma;u^{\ast}TN)=\{h \in {\cal C}(\Sigma;u^{\ast}TN)
|h(y_0)=(0,0,*)\}$ instead of ${\cal C}(\Sigma;u^{\ast}TN)$, then the
construction of ${\cal U}_b$ is the same as for $M^+$.
\vskip 0.1in
\noindent
We can define a bundle
$${\cal F}_{A_i}(N,g_i,T_{m_i})\rightarrow {\cal B}_{A_i}(N,g_i,T_{m_i})$$
whose fiber at $b=(u,\Sigma,j;{\bf p})$ is an infinite dimensional
vector space ${\cal C}(u^{\ast}TN\otimes \wedge^{0,1}).$
\vskip 0.1in
\noindent
For any $D\in {\cal D}^{J,A}_{g,T_m}$,let ${\cal B}_D(M^+,g,T_m)
\subset \overline{\cal B}_A(M^+,g,T_m)$ be the set of $C^{\infty}$ stable maps
whose domain and the corresponding fundamental class of each component and
the sets of periodic orbits corresponding to each puncture point have type
$D$. Then, ${\cal B}_D(M^+,g,T_m)$ is a strata of $\overline{\cal B}_A(M^+,g,T_m)$.
To define the bundle
${\cal F}_{D}(M^+,g,T_{m})\rightarrow {\cal B}_{D}(M^+,g,T_{m})$ and
describe a neighborhood ${\cal U}_b$ for any $b\in {\cal B}_D(M^+,g,T_m)$,
we consider two  simple cases, the discussion for general case is the same.
\vskip 0.1in
\noindent
{\bf Case 1}. Assume that $D$ has two components $(\Sigma_1,j_1)$ and $(\Sigma_2,j_2)$
joining at $p$, and any $b\in {\cal B}_D(M^+,g,T_m)$ has the form:
$b=(u_1,u_2; \Sigma_1 \wedge \Sigma_2,j_1,j_2)$,
where $u_i: \Sigma_i \rightarrow M^+$ are $C^{\infty}$ stable map with
$u_1(p)=u_2(p).$ A neighborhood ${\cal U}_b$ of $b$ can be described as
$$O_{j_1}\times O_{j_2}\times
\left\{\left(\exp_{u_1}(h_1 + \hat{h}_{10}),\exp_{u_2}(h_2 + \hat{h}_{20})\right)|
h_i\in {\cal C}(\Sigma;{u_i}^{\ast}TM^{+})\right.,$$
$$\left.h_{10}=h_{20}\in T_{u_1(p)}M^+, \|h_i\|_{1,p,\alpha} + |h_0| <\epsilon \right\}/stab_u.$$
The fiber of ${\cal F}_{D}(M^+,g,T_{m})\rightarrow {\cal B}_{D}(M^+,g,T_{m})$
at $b$ is ${\cal C}(u_1^{\ast}TM^+\otimes \wedge^{0,1})\times
{\cal C}(u_2^{\ast}TM^+\otimes \wedge^{0,1}).$
\vskip 0.1in
\noindent
{\bf Case 2}. Assume that $D$ has two components $(\Sigma_1,j_1)$ and $(\Sigma_2,j_2)$
 joining at $p$, and any $b\in {\cal B}_D(M^+,g,T_m)$ has the form:
$b=(u_1,u_2; \Sigma_1\wedge\Sigma_2,j_1,j_2)$,
where $u_1: \Sigma_1 \rightarrow M^+$ and
$u_2 :\Sigma_2\rightarrow {\bf R}\times \widetilde{M}$
are $C^{\infty}$ stable maps such that
$u_1$ and $u_2$ converge to a same $k$-periodic orbit when $z_1$ and $z_2$
converge to $p$. A neighborhood ${\cal U}_b$ of $b$ can be described as
$$O_{j_1}\times O_{j_2}\times
\left\{(\exp_{u_1}(h_1 + \hat{h}_{10}),\exp_{u_2}(h_2 + \hat{h}_{20}))|
h_1\in {\cal C}(\Sigma;{u_1}^{\ast}TM^{+}),
h_2\in {\cal C}(\Sigma;{u_2}^{\ast}({\bf R}\times \widetilde{M})),\right.$$
$$\left. h_2(y_0)=(0,0,*),\;\;\;\;\;
\tilde{h}_{10}=\tilde{h}_{20}\in TZ,\;\;\;\; \|h_i\|_{1,p,\alpha} + |h_{i0}|
<\epsilon \right\}/stab_u,$$
where we write $h_0 = (*,*,\tilde{h}_0)$ with $\tilde{h}_{0}\in TZ$.
The fiber of ${\cal F}_{D}(M^+,g,T_{m})\rightarrow {\cal B}_{D}(M^+,g,T_{m})$
at $b$ is ${\cal C}(u_1^{\ast}TM^+\otimes \wedge^{0,1})\times
{\cal C}(u_2^{\ast}({\bf R}\times \widetilde{M})\otimes \wedge^{0,1}).$
\vskip 0.1in
\noindent
The following lemma is obviously holds.
\vskip 0.1in
\noindent
{\bf Lemma 4.4: }{\it ${\cal B}_D(M^+,g,T_m)$ is a  Hausdorff Frechet
V-manifold. The bundle
${\cal F}_{D}(M^+,g,T_{m})\rightarrow {\cal B}_{D}(M^+,g,T_{m})$ is smooth.}
\vskip 0.1in
    \noindent
    Now we construct a finite dimensional $V$-bundle ${\cal E}_D$ over
${\cal B}_{D}(M^+,g,T_{m})$. The construction imitates Siebert's
construction. First of all, we can slightly deform $\omega$ such
that $[\omega]$ is a rational class. By taking multiple, we can
assume $[\omega]$ is an integral class on $\overline{M}^+$.
Therefore, it is the Chern class of a complex line bundle $L$ over
$\overline{M}^+$. We choose a unitary connection
$\bigtriangledown$ on $L$. There is a line bundle associated with
the domain of stable maps called dualized tangent sheaf $\lambda$.
The restriction of $\lambda$ on $(\Sigma,j)$ is
$\lambda_{(\sigma,j)}$, the sheaf of meromorphic 1-form with at
worst simple pole at the special points ( marked points, puncture
points and intersection points ) and for each intersect point $p$,
say $\Sigma_1 $ and $\Sigma_2$ intersects at $p$, $$ Res_p(
\lambda_{(\Sigma_1,j_1)})= Res_p( \lambda_{(\Sigma_2,j_2)}).$$ For
simplicity, we consider two types of stratas as in {\bf Case 1}
and {\bf Case 2}. The general case is identical. For the case {\bf
1} we consider the line bundle $L$, while for the case {\bf 2} we
use $L \cup L'$ instead of $L$, where $L'= p^*(L|_Z)$, and $p; \Re
\rightarrow Z$ is the projection. $L \cup L'$ is a line bundle
over $\overline{M}^+\cup_Z \Re$. To simplify notations we simply
write $L \cup L'$ as $L$. Then $u^*L$ is a line bundle with a
unitary connection over $\Sigma$. It is well-known in differential
geometry that $u^*L$ has a holomophic structure compatible with
the unitary connection. Therefore, $u^*L\otimes
\lambda_{(\Sigma,j)}$ is a holomorphic line bundle. Moreover, if
$D^{\prime}$ is not a ghost component, $\omega(D^{\prime})>0$
since it is represented by a $J$-holomorphic map. Therefore,
$C_1(u^*L)(D^{\prime})>0$. For ghost component $D^{\prime}$,
$\lambda_{D^{\prime}}$ is positive. By taking the higher power of
$u^*L\otimes\lambda_{(\Sigma,j)}$ if necessary , we can assume
that $u^*L\otimes\lambda_{(\Sigma,j)}$ is very ample. Hence,
$H^1(\Sigma, u^*L\otimes \lambda_{(\Sigma,j)})=0$. Therefore,
${\cal E}_b=H^0(\Sigma, u^*L\otimes \lambda_{(\Sigma,j)})$ is of
constant rank (independent of $b$). For sections of $u^*L\otimes
\lambda_{(\Sigma,j)}$ and $u^*L\otimes \lambda_{(\Sigma,j)}\otimes
\wedge^{0,1}$ we may define norms $\|\;\|_{1,p,\alpha}$ and
$\|\;\|_{p,\alpha}$, and define $W^{1,p,\alpha}(\Sigma;u^*L\otimes
\lambda_{(\Sigma,j)})$ and $L^{p,\alpha}(u^*L\otimes
\lambda_{(\Sigma,j)}\otimes \wedge^{0,1})$ in a similar way as
before. Then $H^0(\Sigma, u^*L\otimes \lambda_{(\Sigma,j)})$ is
the $ker \bar{\partial}$, where $$\bar{\partial}:
W^{1,p,\alpha}(\Sigma;u^*L\otimes \lambda_{(\Sigma,j)})
\rightarrow L^{p,\alpha}(u^*L\otimes \lambda_{(\Sigma,j)}\otimes
\wedge^{0,1})$$ is the $\bar {\partial}$-operator. For any $D\in
{\cal D}^{J,A}_{g,T_m}$ we put ${\cal E}_D= \cup_{b\in {\cal
B}_{D}(M^+,g,T_{m})}{\cal E}_b$. It is easy to prove that
    \vskip
0.1in \noindent {\bf Lemma 4.5: }{\it ${\cal E}_D \rightarrow
{\cal B}_{D}(M^+,g,T_{m})$ is a smooth Frechet V-bundle.} \vskip
0.1in \noindent Finally, we define $$\overline{\cal
F}_{A}(M^{+},g,T_m)|_{{\cal B}_{D}(M^+,g,T_{m})} = {\cal
F}_{D}(M^+,g,T_{m})\leqno(4.8)$$ $$\overline{\cal
E}_{A}(M^{+},g,T_m)|_{{\cal B}_{D}(M^+,g,T_{m})} = {\cal
E}_{D}(M^+,g,T_{m}).\leqno(4.9)$$ We use gluing arguments to patch
different strata together ( see below for the details). Note that
$\overline{\cal F}$ (resp. $\overline{\cal E}$ ) is not a bundle
in general because there is no local trivialization. However, for
the gluing argument, it behaves as an infinite dimensional ( resp.
finite dimensional) vector bundle. By abusing the notation, we
will call $\overline{\cal F}$ and $\overline{\cal E}$ bundles.
\vskip 0.1in \noindent Next, we define a map $\eta: \overline{\cal
E}\rightarrow \overline{\cal F}$, and a stabilization ${\cal S}_e$
of $\bar{\partial}$ over a neighborhood ${\cal U}$ of
$\overline{\cal{M}}_{A} (M^{+},g,T_m)$: $${\cal S}_e:
\overline{\cal E}\rightarrow \overline{\cal F}$$ $${\cal S}_e =
\bar{\partial}_J + \eta $$ such that
\begin{description}
\item[(1)] The restriction of $\eta$ to every strata is a smooth bundle map;
\item[(2)] For any $b=(u,(\Sigma,j);{\bf p}) \in {\cal U}$,
$D{\cal S}_e = D\bar{\partial}_J + \eta $
is surjective.
\end{description}
Here we consider $\bar{\partial}_J(u)$ as a map from ${\cal E}_b$ into
${\cal F}_b$ in a natural way.
For any $b=(u,\Sigma,j;{\bf p}) \in \overline{\cal{M}}_{A} (M^{+},g,T_m)$,
$coker D_u$ is a representation space of $stab_u$. Hence it can be decomposed
as sum of irreducible representations. There is a result in algebra saying
that the irreducible factors of group ring contain all the irreducible
representations of finite group . Hence, it is enough to find a copy of group
ring in $E_b$. This is done by algebraic geometry. Recall that
${\cal E}_b=H^0(\Sigma, u^*L\otimes \lambda_{(\Sigma,j)})$.
By taking power of $u^*L\otimes \lambda_{(\Sigma,j)}$ if necessary, we can
assume that $u^*L\otimes \lambda_{(\Sigma,j)}$ induces an embedding
of $\Sigma$ into $CP^N$ for some $N$. Furthermore, since $u^*L\otimes
\lambda_{(\Sigma,j)}$ is invariant under
$stab_u$, $stab_u$ also acts effectively naturally on $CP^N$. Pick any point
$x_0\in im(\Sigma)\subset CP^N$ such that $\sigma_i(x_0)$ are mutually different
for any $\sigma_i\in stab_u$. Then, we can find a homogeneous polynomial $F$ of
some degree, say $\ell$ such that $F(x_0)\neq 0, F(\sigma_i(x_0))=0$ for
$\sigma_i\neq 1$.
Notes that $\sigma_i^*F$ generates a group ring.
Next, we note that $F\in H^0(O(\ell))$. By pull back over $\Sigma$, $F$
induces a section
of $H^0(\Sigma,(u^*L\otimes \lambda_{(\Sigma,j)})^{\ell})$. Therefore, if we replace
$u^*L\otimes \lambda_{(\Sigma,j)}$ by $(u^*L\otimes \lambda_{(\Sigma,j)})^{\ell}$
and redefine
${\cal E}_b=H^0(\Sigma, (u^*L\otimes \lambda_{(\Sigma,j)})^{\ell})$,
${\cal E}_b$ contains a copy of group ring.
    Now, it is clear how to construct a map $\eta$. For
$b\in \overline{\cal{M}}_{A} (M^{+},g,T_m)$, we decompose $Coker D_u$ as
irreducible representations. Pick a factor and project group ring to this irreducible
factor. If there is more
than one factors, we take direct sum of ${\cal E}$ to get a surjective map $\eta_b$ to
$Coker D_u$.
\vskip 0.1in
\noindent
Now we extend the map $\eta_b$ to a neighborhood of $b$. If $\Sigma$ is smooth,
we may choose a neighborhood $O(b)$ of $b$
and take local travilizations for both ${\cal E}$ and
${\cal F}$ over $O(b)$, and extend the map $\eta_b$ to $O(b)$ in a natural way.
By choosing $O(b)$ so small we
may assume that $\eta$ is surjective on $O(b)$.
Then we use cut-off function to extend $\eta$ over whole ${\cal U}$.
In the following we assume that
$\Sigma$ is a nodel surface. We will consider two types of stratas {\bf Case 1}
and {\bf Case 2} as above. We use the holomorphic cylindrical coordinates $(s_i,t_i)$
on $\Sigma_i$ near $p$, and write
$$\Sigma_1-\{p\}=\Sigma_{10}\bigcup\{[0,\infty)\times S^1\},$$
$$\Sigma_2-\{p\}=\Sigma_{20}\bigcup\{(-\infty,0]\times S^1\}.$$
We discuss the two cases separately.
\vskip 0.1in
\noindent
{\bf Case 1.} For any $(r,\tau)$ we construct a surface
$\Sigma_{(r)} =\Sigma_1 \#_{(r)} \Sigma_2 $ as follows, where and later we
use $(r)$ to denote gluing parameters. We cut off the
part of $\Sigma_i$ with cylindrical
coordinate $|s_i|>3r$ and glue the remainders along
the collars of length $2r$ of the cylinders with the gluing formulas:
$$s_1=s_2 + 4r\leqno(4.10)$$
$$t_1=t_2 + \tau.\leqno(4.11)$$
We glue the map $(u^{+},u^{-})$ to get a map $u_{(r)}$
from $\Sigma_{(r)}$ to $M^+$ as follows. Set\\
\[
u_{(r)}=\left\{
\begin{array}{ll}
u_1 \;\;\;\;\; on \;\;\Sigma_{10}\bigcup\{(s_1,t_1)|0\leq s_1 \leq
r, t_1 \in S^1 \}    \\  \\
u_1(p)=u_2(p) \;\;\;\;\;
on \;\;\{(s_1,t_1)| \frac{3r}{2}\leq s_1 \leq
\frac{5r}{2}, t_1 \in S^1 \}  \\   \\
u_2 \;\;\;\;\; on \;\;\Sigma_{20}\bigcup\{(s_2,t_2)|0\geq s_2
\geq - r, t_2 \in S^1 \}    \\    \\
\end{array}
\right.
\]
To define the map $u_{(r)}$ in the remaining part we fix a smooth cutoff
function $\beta : {\bf R}\rightarrow [0,1]$ such that
\[
\beta (s)=\left\{
\begin{array}{ll}
1 & if\;\; s \geq 1 \\
0 & if\;\; s \leq 0
\end{array}
\right.
\]
and $|\beta^{\prime}(s)|\leq 2.$
We assume that $r$ is large enough such that $u_i$ maps
the tube $\{(s_i,t_i)||s_i|\geq r,t_i\in S^1 \}$ into a normal coordinate
domain of $u_i(p)$. We define\\
$$u_{(r)}= u_1(p)+ \left(\beta(3-\frac{2s_1}{r})(u_1(s_1,t_1)-u_1(p)) +
\beta(\frac{2s_1}{r}-5)(u_2(s_1,t_1)- u_1(p))\right).$$
\vskip 0.1in
\noindent
{\bf Case 2.} We choose a Darboux coordinate system $a_i,\theta_i,w_i$ on the cylinder near the
periodic orbit. Suppose that
$$a_i(s_i,t_i)-ks_i-\ell_i\rightarrow 0\;\;\;\;
\theta(s_i,t_i)-kt_i-\theta_{i0}\rightarrow 0$$
In view of the condition $u_2(y_0)=(0,0,*)$ for some point $y_0$, we may
assume that $\ell_2=0,\;\; \theta_{20}=0$.
For any $(r,\theta_0)$ we glue $M^+$ and ${\bf R}\times \widetilde{M}$ with
parameter $(r,\theta_0)$ to get again $M^+$ with gluing formula:
$$a_1=a_2 + 4kr + \ell\leqno(4.12)$$
$$\theta_1 = \theta_2 + \theta_0 \;\;mod \;1.\leqno( 4.13)$$
Now we construct a surface
$\Sigma_{(r)} =\Sigma_1\#_{(r)} \Sigma_2 $ with gluing formulas:
$$s_1=s_2 + 4r\leqno(4.14)$$
$$t_1=t_2 + \frac{\theta_{0}-\theta_{10} + n}{k}\leqno(4.15)$$
for some $n \in Z_k$.
To get a pregluing map $u_{(r)}$ from $\Sigma_{(r)}$ we set\\
\[
u_{(r)}=\left\{
\begin{array}{ll}
u_1 \;\;\;\;\; on \;\;\Sigma_{10}\bigcup\{(s_1,t_1)|0\leq s_1 \leq
r, t_1 \in S^1 \}    \\  \\
\left(ks_1, x(kt_1)\right) \;\;\;\;\;
on \;\;\{(s_1,t_1)| \frac{3r}{2}\leq s_1 \leq
\frac{5r}{2}, t_1 \in S^1 \}  \\   \\
u_2 \;\;\;\;\; on \;\;\Sigma_{20}\bigcup\{(s_2,t_2)|0\geq s_2
\geq - r, t_2 \in S^1 \}.    \\    \\
\end{array}
\right.
\]
We assume that $r$ is large enough such that $u_i$ maps the tube
$\{(s_i,t_i)||s_i|\geq \frac{r}{k},t_i \in S^1 \}$ into a domain
with Darboux coordinates $(a_i,\theta_i, w_i)$. We write in in
terms the coordinate system $(a_1,\theta_1, w_1)$
$$u_{(r)}=\left(a_{(r)}, \widetilde{u}_{(r)}\right)=
\left(a_{(r)}, \theta_{(r)}, w_{(r)} \right)$$ and define\\
$$a_{(r)}= ks_1  + \left(\beta(3-\frac{2s_1}{r}) (a_1(s_1,t_1)-
ks_1-\ell) + \beta(\frac{2s_1}{r}-5) (a_2(s_1,t_1)- ks_1 -
\ell)\right)$$ $$\widetilde{u}_{(r)}=x(kt_1)+
\left(\beta(3-\frac{2s_1}{r}) (\widetilde{u}_1(s_1,t_1)-x(kt_1))+
\beta(\frac{2s_1}{r}-5)
(\widetilde{u}_2(s_1,t_1)-x(kt_1))\right).$$ \vskip 0.1in
\noindent It is easy to check that in both cases $u_{(r)}$ is a
smooth function. By using the exponential decay of $u_i$ one can
easily prove that $u_{(r)}$ are a family of approximate
J-holomorphic map, precisely the following lemma holds. \vskip
0.1in \noindent {\bf Lemma 4.6 }
$$\|\bar{\partial}_{\Sigma_{\xi,r}}(u_{r})\|_{p,\alpha,r}\leq
Ce^{-(\delta-\alpha)r} .\leqno(4.16)$$ {\it The constants C in the
above estimates are independent of $r$.} \vskip 0.1in \noindent
Let $b_{(r)}= (u_{(r)},\Sigma_{(r)},j_1,j_2)$. For any section
$\xi\in C^{\infty}(\Sigma_{(r)};
u_{(r)}^*L\otimes\lambda_{(\Sigma_{(r)},j_1,j_2)})$, we denote
$$\beta\xi = (\xi_1, \xi_2) = \left( \xi\beta(2\alpha r + 1
-\alpha s_1), \xi\beta(2\alpha r + 1 + \alpha
s_2)\right).\leqno(4.17)$$ We may consider $\beta\xi$ as a section
of the bundle $u^*L\otimes\lambda_{(\Sigma,j)}$, where
$u=(u_1,u_2)$, $\Sigma=\Sigma_1\wedge\Sigma_2$, $j=(j_1,j_2)$.
Denote by $\pi$ the projection from $W^{1,p,\alpha}(\Sigma,
u^*L\otimes \lambda_{(\Sigma,j)})$ into $H^0(\Sigma,
u^*L\otimes\lambda_{(\Sigma,j)})$. We have a map $$i_{r}: {\cal
E}_{b_{(r)}}\rightarrow {\cal E}_b$$ $$i_{r}(\xi)= \pi
(\beta(\xi)).\leqno(4.18)$$ The proof of the following lemma is
similar the proof of Lemma 4.9, we omit it. \vskip 0.1in \noindent
{\bf Lemma 4.7:} {\it $i_{r}$ is a isomorphism for $r$ big enough.
} \vskip 0.1in \noindent By multiplying a cut-off function we may
assume that $$\eta_b = 0 \;\;\; for |s_i|\geq R \leqno(4.19 )$$
and $D_u + \eta_b$ remains surjective for $R$ big enough. Consider
a neighborhood $U_{\epsilon}(b_{(r)})=O_{j_1,\epsilon}\times
O_{j_2,\epsilon}\times\{\exp_{u_{(r)}}h | \|h\|_{1,p}<\epsilon
\},$ where $O_{j_i,\epsilon}$ is a $\epsilon$- neighborhood of
$j_i$ with respect to a metric in the space of complex structures.
Let $b'= (u'_{(r)},\Sigma_{(r)},j'_1,j'_2) \in
U_{\epsilon}(b_{(r)})$. To simplify notations, without loss of
generality, we assume that $j'_i = j_i$. The paralell transport to
$u_{(r)}$ along geodesics with respect to the unitary connection
on $L$ induces a map $T$ from $u_{(r)}^{\prime *}L\otimes
\lambda_{(\Sigma_{(r)},j_1,j_2)}$ into $u_{(r)}^*L\otimes
\lambda_{(\Sigma_{(r)},j_1,j_2)}.$ For any $\xi \in {\cal E}_{b'}$
we define $$\eta_{b'}(\xi)= \eta_b(I_{r}(T\xi)).$$ We define in
the neighborhood $U_{\epsilon}(b_{(r)})$ $${\cal S}_e =
\bar{\partial}_J + \eta. $$ Then $$D{\cal S}_e = D\bar{\partial}_J
+ \eta.$$ For any $\eta \in
C^{\infty}(\Sigma_{(r)};u_{(r)}^{\ast}TM^+\otimes \wedge^{0,1}))$,
let $\eta_{i}$ be its restriction to the part
$\Sigma_{i0}\bigcup\{(s_i,t_i)|\;|s_i| < 3r\}$, extended by zero
to yield a section over $\Sigma_i$. Define
$$\|\eta\|_{p,\alpha,r}=\|\eta_{1} \|_{\Sigma_1,p,\alpha} +
\|\eta_{2} \|_{\Sigma_2,p,\alpha}.\leqno(4.20)$$ Denote the
resulting completed spaces by $L_{r}^{p,\alpha}$. For any section
$(\xi, h)\in {\cal E}_{b_{(r)}}\times
C^{\infty}(\Sigma_{(r)};u_{(r)}^{\ast}TM^+)$, denote $$h_0
=\int_{S^1}h(2r,t)dt \leqno(4.21)$$ $$h_1 = h\beta(2\alpha r + 1
-\alpha s_1)\leqno(4.22)$$ $$h_2 = h\beta(2\alpha r + 1 + \alpha
s_2).\leqno(4.23)$$ We define $$\|(\xi, h)\|_{1,p,\alpha,r}=
|\xi|_{1,p,\alpha} + \|h_1\|_{\Sigma_1,1,p,\alpha} +
\|h_2\|_{\Sigma_2,1,p,\alpha} + |h_0|\leqno(4.24)$$ where
$|\;\;|_{1,p,\alpha}$ is the norm on ${\cal E}_{b_{(r)}}$. Denote
the resulting completed spaces by ${\cal W}_{r}^{1,p,\alpha}$.
\vskip 0.1in \noindent {\bf Lemma 4.8 } {\it $D{\cal
S}_{e,b_{(r)}}$ is surjective for $r$ large enough and $\epsilon$
small enough . Moreover, there are a right inverses $Q_{b_{(r)}}$
such that $$D{\cal S}_{e,b_{(r)}}Q_{b_{(r)}}=Id \leqno(4.25)$$
$$\|Q_{b_{(r)}}- Q_b\|\leq \frac{C}{r} \leqno(4.26)$$ $$ \left
\|\left(\frac{\partial}{\partial r}Q_{b_{(r)}}\right)
\left|_{|s_i|\leq r}\right\|\leq \frac{C}{r^2}\right.\leqno(4.27)
$$ for some constant $C>0$ independent of $(r)$ and $\epsilon$.
Here and later the inequality (4.27) means for any $\eta$ , $$
\left\|\left(\frac{\partial}{\partial r}Q_{b_{(r)}}\right)\eta
\left|_{|s_i|\leq r}\right\|_{1,p,\alpha,r} \leq
\frac{C}{r^2}\|\eta\|_{p,\alpha,r}\right. \leqno(4.28) $$ } \vskip
0.1in \noindent {\bf Proof:} The proof is similar to the proof in
\cite{MS}. We first construct an approximate right inverse
$Q_{b_{(r)}}^{\prime}$ such that the following estimates hold:
$$\|Q_{b_{(r)}}^{\prime}\|\leq C \leqno(4.29)$$ $$\|D{\cal
S}_{e,b_{(r)}}Q_{b_{(r)}}^{\prime} - Id \|\leq \frac{1}{2}.
\leqno(4.30)$$ Then the operator $D{\cal
S}_{e,b_{(r)}}Q_{b_{(r)}}^{\prime}$ is invertible and a right
inverse $Q_{b_{(r)}}$ of $D{\cal S}_{e,b_{(r)}}$ is given by
$$Q_{b_{(r)}}=Q_{b_{(r)}}^{\prime}(D{\cal
S}_{e,b_{(r)}}Q_{b_{(r)}}^ {\prime})^{-1}.\leqno(4.31)$$ Given
$\eta \in L_{r}^{p,\alpha}$, we have a pair $(\eta_1,\eta_2)$. Let
$ (\xi, h)= Q_b(\eta_1,\eta_2).$ We may write $h$ as
$(h_1+\hat{h_0},h_2+\hat{h_0})$, and define $$h_{(r)}= \hat{h_0} +
h_1\beta(3/2 - \frac{s_1}{2r}) + h_2\left(1 - \beta(3/2 -
\frac{s_1}{2r})\right).\leqno(4.32)$$ Then we define
$$Q_{b_{(r)}}^{\prime}\eta = (i_r^{-1}(\xi),
h_{(r)}).\leqno(4.33)$$ We must prove that
$$\|\eta_{b_{(r)}}(i_r^{-1}(\xi)) + D_{u_{(r)}}h_{(r)} - \eta
\|_{p,\alpha,r} \leq
\frac{1}{2}\|\eta\|_{p,\alpha,r}.\leqno(4.34)$$ By definition
$\eta_{b_{(r)}}(i_r^{-1}(\xi))=\eta_{b}(\xi)$. Since
$\eta_{b}(\xi) + D_u h = \eta$ the term on the left hand side
vanishes for $|s_i|\leq r$. It suffices to estimate the left hand
side in the annulus $r\leq |s_i| \leq 3r$. Denote
$\beta_1=\beta(3/2 - \frac{s_1}{2r})$, $\beta_2 = \left(1 -
\beta(3/2 - \frac{s_1}{2r})\right)$. Note that in this annulus
$$\beta_1 + \beta_2 =
1,\;\;\;\eta_{b}(\xi)=0,\;\;\;D_{u_i}(h_i+\hat{h_0}) = \eta_i.$$
$$\beta_1D_{u_1}(h_1+\hat{h_0}) +
\beta_2D_{u_2}(h_2+\hat{h_0})=\eta.$$ Since near $u_1(p)=u_2(p)$ (
or near the periodic orbit $x(kt)$), $D = \bar{\partial}_J + S$,
we have $$ D{\cal S}_{e,b_{(r)}}Q_{b_{(r)}}^{\prime}\eta \;-\;
\eta = \eta_{b_{(r)}}(i_r^{-1}(\xi)) + D_{u_{(r)}}h_{(r)} \;-
\;\eta \leqno(4.35)$$ $$=(\bar{\partial}\beta_1)h_1 +
\beta_1(S_{u(r)}-S_{u_1})h_1 + (\bar{\partial}\beta_2)h_2 +
\beta_2(S_{u(r)}-S_{u_2})h_2 - \beta_1S_{u_1}\hat{h_0} -
\beta_2S_{u_2}\hat{h_0} + S_{u(r)}\hat{h_0}.$$ By the exponential
decay of $S$ we get $$\left\|\eta_{b_{(r)}}(i_r^{-1}(\xi)) +
D_{u_{(r)}}h_{(r)} - \eta \right\|_{p,\alpha,r} \leq
\frac{C_1}{r}(\|h_1\|_{p,\alpha} + \|h_2\|_{p,\alpha} + |h_0|)\leq
\frac{C_2}{r}\|\eta \|_{p,\alpha,r}.\leqno(4.36)$$ In the last
inequality we used that $\|Q_b\|\leq C$ and (4.20). Then (4.30)
follows by choosing $r$ big enough. The proof of (4.27) is easy,
we omit it. Here we prove (4.28). From (4.35) and the exponential
decay of $S$ we get $$\left\|\frac{\partial}{\partial
r}\left(D{\cal S}_{e,b_{(r)}}Q_{b_{(r)}}^{\prime}\right) \right\|
\leq \frac{C}{r^2}.\leqno(4.47)$$ From $$ D{\cal
S}_{e,b_{(r)}}Q_{b_{(r)}}^{\prime}\left(D{\cal S}_{e,b_{(r)}}
Q_{b_{(r)}}^{\prime}\right)^{-1}=Id$$ and (4.37) we conclude that
$$\left\|\frac{\partial}{\partial r}\left((D{\cal
S}_{e,b_{(r)}}Q_{b_{(r)}}^{\prime} )^{-1}\right)\right\|\leq
\frac{C}{r^2}.\leqno(4.38)$$ Since for any $\eta \in
L_{r}^{p,\alpha}$ the restriction of $Q_{b_{(r)}}^{\prime}$ to the
parts $|s_i|\leq r$ is $\eta$, from (4.31 ) and (4.38 ) we get $$
\left\|\frac{\partial}{\partial r}Q_{b_{(r)}}\left|_{|s_i|\leq
r}\right.\right\| \leq
\frac{C}{r^2}.\;\;\;\;\;\;\;\;\;\;\;\;\Box$$ \vskip 0.1in
\noindent Thus for every point $b\in \overline{\cal
M}_A(M^+,g,T_m)$ we may choose a neighborhood $O(b)$ such that
$D{\bf S}_e$ is surjective on $O(b)$. Then we use cut-off function
to extend $\eta$ over whole ${\cal U}$. By the compactness of
$\overline{\cal M}_A(M^+,g,T_m)$, we may choose finitely many such
neighborhoods $O(b_j), j=1,2,...,N$ such that $\overline{\cal
M}_A(M^+,g,T_m) \subset \bigcup_j^{N} O(b_j)$. Then we take ${\cal
E}^{N}$ instead of ${\cal E}$ and put $$\eta = \sum_j
\eta_{b_j}.$$ Our stabilization operator is the bundle map $${\cal
S}_e=\bar{\partial}_J + \eta : {\cal E}\rightarrow {\cal
F}.\leqno(4.39)$$ It is easy to see that the restriction of $\eta$
to every strata is a smooth bundle map, and for any
$b=(u,\Sigma,j;{\bf p}) \in {\cal U}$ $D{\cal S}_e =
D\bar{\partial}_J + \eta $ is surjective. \vskip 0.1in \noindent
We continue to discuss the gluing for the ${\bf Case 1}$ and the
${\bf Case 2}$. We define a map $I_{r}: ker D{\cal
S}_{e,b_{(r)}}\rightarrow ker D{\cal S}_{e,b}$ as in (4.17). For
any $h \in ker D{\cal S}_{e,b_{(r)}}$ we denote by $h_i$ the
restriction of $h$ to the part $|s_i|\leq 2r + \frac{1}{\alpha}$.
We get a pair $(h_1,h_2)$. We put $$\beta h = \left(
h_1\beta(2\alpha r + 1 -\alpha s_1), h_2\beta(2\alpha r + 1 +
\alpha s_2)\right)$$ and define $$I_{r}(\xi,h)=(i_r(\xi), \beta
h)- Q_bD{\cal S}_e (i_r(\xi), \beta h)\leqno(4.40)$$ where $Q_b$
is a right inverse of $D{\cal S}_{e,b}$. We also construct a map
$I'_r: ker D{\cal S}_{e,b} \rightarrow ker D{\cal S}_{e,b_{(r)}}$.
Let $(\xi,h)\in ker D{\cal S}_{e,b}$. We write
$h=(h_1+\hat{h_0},h_2+\hat{h_0})$, and define $$h_{(r)}= \hat{h_0}
+ h_1\beta(3/2 - \frac{s_1}{2r}) + h_2\left(1 - \beta(3/2 -
\frac{s_2 + 4r}{2r})\right).\leqno(4.41)$$ We define
$$I'_r(\xi,h)= (i_r^{-1}(\xi), h_{(r)}) - Q_{b_{(r)}}D{\cal
S}_{e,b_{(r)}} (i_r^{-1}(\xi), h_{(r)}).\leqno(4.42)$$
    \vskip 0.1in \noindent {\bf Lemma 4.9:} Both $I_{r}$ and $I'_r$ are
isomorphisms for $r$ big enough. Moreover
$$\left\|\frac{\partial}{\partial r}I'_{r}\left|_{|s_i|\leq
r}\right.\right\| \leq \frac{C}{r^2}.\leqno(4.43)$$ \vskip 0.1in
\noindent {\bf Proof:} The proof is basically a similar gluing
argument as in an unpublished book of Donaldson. The proof is
devides into 2 steps. \vskip 0.1in \noindent {\bf Step 1}. We show
that $I_{r}$ is injective for $r$ big enough. We have from (4.40)
$$\|I_{r}(\xi,h) - (i_{r}(\xi), \beta h)\|_{1,p,\alpha}
 \leq C_1 \| \eta_b(i_{r}(\xi)) + D_u(\beta h)\|_{p,\alpha}$$
$$= C_1 \| \eta_b(i_{r}(\xi)) + (\bar{\partial}\beta)h + \beta \left(D_u h
+ D_{u_{(r)}}h + \eta_{b_{(r)}} (\xi) - D_{u_{(r)}}h - \eta_{b_{(r)}} (\xi)\right) \|_{p,\alpha}$$
$$= C_1 \|(\bar{\partial}\beta)h +
\beta \left( (S_u - S_{u(r)})h\right)\|_{p,\alpha},$$
where we used the fact that $(\xi,h)\in ker D{\bf S}_{e,b_{(r)}}$,
$\eta_b(i_{r}(\xi)) = \beta\eta_{b_{(r)}} (\xi)$ for $r>R$.
Note that
$$S_u = S_{u(r)} \;\;\;\;if \; \; s_1\leq \; r, \; or \;\; s_2 \geq -r $$
$$\beta(2\alpha r + 1 -\alpha s_1)= 1 \;\;\;if \;\;s_1\leq 2r$$
$$\beta(2\alpha r + 1 + \alpha s_2)= 1 \;\;\;if \;\;s_2\geq -2r.$$
By exponential decay of $S$ we have
$$\|(S_u - S_{u(r)})\beta h\|_{p,\alpha} \leq Ce^{-\delta r}\|\beta h\|_{1,p,\alpha}$$
for some constant $C>0$.
Since $(\bar{\partial}\beta(2\alpha r + 1 -\alpha s_1))h$ supports in
$2r \leq s_1 \leq 2r + \frac{1}{\alpha}$, and over this part
$$|\bar{\partial}\beta(2\alpha r + 1 -\alpha s_1)|\leq |\alpha|$$
$$\beta(2\alpha r + 1 + \alpha s_2)= 1,\;\;\;e^{2\alpha|s_1|}\leq e^4e^{2\alpha|s_2|},$$
we obtain
$$\|(\bar{\partial}\beta(2\alpha r + 1 -\alpha s_1))h\|_{p,\alpha}\leq
|\alpha|e^4 \|h_2\|_{p,\alpha} \leq |\alpha|e^4 \|\beta h\|_{p,\alpha}.$$
Similar inequality for $(\bar{\partial}\beta(2\alpha r + 1 + \alpha s_2))h$
also holds. So we have
$$\|(\bar{\partial}\beta)h\|_{p,\alpha} \leq 2|\alpha|e^4 \|\beta h\|_{1,p,\alpha}.$$
Hence
$$\|I_{r}(\xi,h) - (i_{r}(\xi), \beta h)\|_{1,p,\alpha}\leq C_3(|\alpha| +
e^{-\delta r})\|\beta h\|_{1,p,\alpha} \leq 1/2 \|\beta h\|_{1,p,\alpha}\leqno(4.44)$$
for some constant $C_3>0$, here we choosed $r$ big enough and $|\alpha|$ so small
that $C_3(|\alpha| + e^{-c_2r}) < 1/2$. Now suppose that $I_{r}(\xi,h)=0$, then
(4.44) gives us $|i_{r}(\xi)|_{1,p,\alpha}=0$, $\|\beta h\|_{1,p,\alpha}=0$.
It follows that $\xi = 0, \;\; h=0$.
\vskip 0.1in
\noindent
{\bf Step 2}. By a similar culculation as the proof of Lemma 4.8 we obtain
$$\|I'_r(\xi,h) - (i_r^{-1}(\xi), h_{(r)})\|_{p,\alpha,r}\leq \frac{C}{r}
(\|h\|_{p,\alpha} + |h_0|).\leqno(4.45)$$
In particular, it holds for $p=2$. It remains to show that
$\|h_{(r)})\|_{2,\alpha,r}$ is close to $\|h\|_{2,\alpha}$. Denote $\pi$ the
projection into the second component, that is, $\pi(\xi,h)=h$. Then
$\pi(ker D{\cal S}_{e,b})$ is a finite dimentional space. Let $f_i,\;i=1,..,d$
be an orthonormal basis. Then $F=\sum f_i^2e^{2\alpha|s|}$ is an integrable
function on $\Sigma$. For any $\epsilon >0$, we may choose $R_0$ so that
$$\int_{|s_i|\geq R_0}F \leq \epsilon.$$
Then the restriction of $h$ to $|s_i|\geq R_0$ satisfies
$$\|h|_{|s_i|\geq R_0}\|_{2,\alpha}\leq \epsilon\|h\|_{2,\alpha},$$
therefore
$$\|h_{(r)}\|_{2,\alpha,r}\geq \|h|_{|s_i|\leq R_0}\|_{2,\alpha} + |h_0|
\geq (1 - \epsilon)\|h\|_{2,\alpha} + |h_0|,\leqno(4.46)$$
for $r>R_0$. Suppose that $I'_r(\xi,h)=0$. Then (4.45) and (4.46) give us $h=0$,
and so $\xi=0$.
\vskip 0.1in
\noindent
The {\bf step 1} and {\bf step 2} together show that both $I_{r}$ and $I'_r$
are isomorphisms for $r$ big enough. Now we prove (4.43).
By definition $\eta_{(r)}( i_r^{-1}(\xi))$ is independent of $r$, by a direct calculation
we get
$$\left\|\frac{\partial}{\partial r}
\left(D{\cal S}_{e,b_{(r)}}(i_r^{-1}(\xi), h_{(r)})\right)\right\|
\leq \frac{C}{r^2}.\leqno(4.47)$$
Since $h_{(r)}|_{|s_i|\leq r} = h$, using (4.47) we get (4.43). $\Box$

\subsection{\bf  Virtual neighborhood and log invariants }

Let $U_{{\cal S}_e} = ({\cal S}_e)^{-1}(0)$.
Following \cite{R5}, we can assume that
$U_{\S_e}$ is compatible with the stratification of ${\cal E}$ by
defining $\eta_i$ inductively over the strata of
$\overline{\B}_A$. Namely, we define $\eta_i$ on lower strata
first. Then, we extend it over ${\cal E}$ such that $\eta_i$  is
supported in a neighborhood of lower strata. Then, we define
$\eta_{i+1}$ supported away from lower strata. Once $\eta_i$ is
defined in such a fashion.
\vskip 0.1in
\noindent
{\bf Proposition 4.10 }{\it $U_{\S_e}$ has the property\\
{\bf 1.} Each strata of $U_{\S_e}$ is a smooth V-manifold.\\
{\bf 2.} If $\B_{D'}\subset \overline{\B}_D$ is a lower stratum,
$$U_{\S_e}\cap \E|_{\B_{D'}}\subset U_{\S_e}\cap
\E|_{\overline{\B}_D}$$ is a submanifold of codimension
at least 2.}
\vskip 0.1in
\noindent
{\bf Proof:}\\
{\bf 1.} We prove only for the top strata, for lower strata the proof is
the same. Let $b= (u,(\Sigma,j))\in {\cal{M}}_{A} (M^{+},g,T_m)$.
Consider a map
$$ F: {\cal E}_b \times W^{1,p,\alpha}(\Sigma;u^{\ast}TM^+)
\rightarrow L^{p,\alpha}(u^{\ast}TM^+\otimes \wedge^{0,1})$$
$$F(\xi,h)=P_{exp_uh}\left(\bar{\partial}_{J}exp_uh + \eta(T\xi)
\right)$$
where $P_{exp_uh}$ denotes the parallel transport from $exp_uh$ to $u$ along
the geodesics $t \rightarrow exp_u(th)$ and $T$ is a travilization of
$\cal E$ near $b$. Denote by $Q_b$ a right inverse of $D{\cal S}_{e,b}$.
By the Implicity Function Theorem there is smooth map
$$f: kerD{\cal S}_{e,b}\rightarrow L^{p,\alpha}(u^{\ast}TM^+\otimes \wedge^{0,1})$$
such that the zero set of $F$ is locally the form $(\zeta, Q_bf(\zeta))$,
where $\zeta \in kerD{\cal S}_{e,b}$. So $\zeta$ gives a locall coordinate
system of the top strata of $U_{\S_e}$. If we choose another $b'$ we may obtain
another locall coordinate system. We may show the coordinate transformation
is smooth.\\
{\bf 2.} Because of the existance of 2-dimensional gluing paprameters for
both {\bf Case 1} and {\bf Case 2} in above. For general case the situation
is the same.\\
\vskip 0.1in
\noindent
There are two maps, the inclusion map $$I:U_{{\cal
S}_e}\rightarrow {\cal E}\leqno(4.48)$$ and the projection $$\pi:
U_{{\cal S}_e}\rightarrow {\cal U}.\leqno(4.49)$$ $I$ can be viewed
as a section of the bundle $E=\pi^{\ast}{\cal E}$, and $I^{-1}(0)=
\overline{\cal{M}}_{A}(M^{+},g,T_m)$.

\vskip 0.1in \noindent {\bf Lemma 4.11 }{\it I is a proper map.}
\vskip 0.1in The proof is the same as the proof of the compactness
theorem for $\overline{\cal{M}}_{A}(M^{+},g,T_m)$.
 \vskip 0.1in
Using the virtual neighborhood we can define
the log GW-invariants. Recall that we have two natural maps
$$e_i: \overline{\cal B}_A(M^{+},g,T_m) \rightarrow
M^{+} $$ for $i\leq l$ defined by evaluating at marked points
and $$e_j:\overline{\cal B}_A(M^{+},g,T_m)\rightarrow Z$$
for $j>l$ defined by projecting to its
periodic orbit. To define the log GW-invariants, choose  an
$r$-form $\Theta$ on ${\cal E}$ supported in a neighborhood of the zero
section, where $r$ is the dimension of the fiber, such that $$
\int_{E_x}i^{\ast}\Theta=1$$ for any $x\in U_{{\cal
S}_e}$, where $i$ is the inclusion map $E_x \rightarrow E$. We
call $\Theta$ a Thom form. The log GW-invariant can be
defined as
$$\Psi^{(V, B)}_{(A,g,T_m}(\alpha_1,...,
\alpha_{l} ; \beta_{\ell + 1},..., \beta_m)=\int_{U_{{\cal
S}_e}}\prod_i e^*_i\alpha_i\wedge
\prod_j e^*_j\beta_j\wedge I^*\Theta. \leqno(4.50)$$ for
$\alpha_i\in H^*(M^{+}, {\bf R})$ and $\beta_j\in H^*(Z,
{\bf R})$ represented by differential form. Clearly, $\Psi=0$ if
$\sum \deg(\alpha_i)+ \sum \deg (\beta_i)\neq ind$.
We must prove the convergence of the integral (4.50) near each lower strata.
We prove this for {\bf Case 2}, for other case the proof is the same.
We have proved in Lemma (4.8) that there is a uniform right inverse $Q_{(r)}$.
Consider a map
$$ F_{(r)}: {\cal E}_{b_{(r)}} \times W^{1,p,\alpha,r}(\Sigma_{(r)};u_{(r)}
^{\ast}TM^+)
\rightarrow L^{p,\alpha,r}(u_{(r)}^{\ast}TM^+\otimes \wedge^{0,1})$$
$$F_{(r)}(\xi,h)=P_{exp_{u_{(r)}}h}\left(\bar{\partial}_{J}exp_{u_{(r)}}h
+ \eta(T\xi)\right)$$
By the Implicity Function Theorem there is smooth map
$$f_{(r)}: kerD{\cal S}_{e,b_{(r)}}\rightarrow
L^{p,\alpha,r}(u_{(r)}^{\ast}TM^+\otimes \wedge^{0,1})$$
such that the zero set of $F_{(r)}$ is locally the form
$(\zeta, Q_{b_{(r)}}f_{(r)}(\zeta))$,i.e
$$F_{(r)}(\zeta + Q_{b_{(r)}}f_{(r)}(\zeta))=0\leqno(4.51)$$
where $\zeta \in kerD{\cal S}_{e,b_{(r)}}$.
Since there is a isomorphism $I'_r: ker D{\cal S}_{e,b}\rightarrow
ker D{\cal S}_{e,b_{(r)}}$, we have
\vskip 0.1in
\noindent
{\bf Lemma 4.12 } {\it There is a
neighborhood $O$ of $(0,0)$ in $ker D{\cal S}_{e,b}$
and $R_0 > 0$ such that
$$f_{(r.\theta)}\circ I'_{(r,\theta)}: O\times Z_k \rightarrow
U_{{\cal S}_{e,b_{(r,\theta)}}}$$
for $r>R_0$ is a family of orientation preserving local diffeomorphisms.
Moreover,
$(R_0,\infty]\times Z_k \times \left(O/(stab_{u_1}\times
stab_{u_2})\right)/\cong$ is a local  chart around $b$,
where "$\cong$" is an equivalence relation
at $\infty$: $(a_1, b_1, n_1)\cong (a_2, b_2, n_2)$ iff $(a_1,
b_1)=(a_2, b_2)$. }
\vskip 0.1in
\noindent
{\bf Lemma 4.13} { Restricting to the part $|s_i|\leq r$ we have
$$\left\|\frac{\partial}{\partial r}\left(f_{(r)}\circ I'_r(\xi)\right)\right\|
\leq \frac{C}{r^2}.\leqno(4.52)$$}
\vskip 0.1in
    \noindent
    {\bf Proof :} Restricting to the part $|s_i|\leq r$ , $F_{(r)}$ is independent of $r$.
From (4.51) we get $$DF_{(r)}\left(\frac{\partial}{\partial
r}I'_r(\xi) + \frac{\partial}{\partial r}
(Q_{b_{(r)}})f_{(r)}(I'_r(\xi)) +
Q_{b_{(r)}}\frac{\partial}{\partial r}(f_{(r)}
(I'_r(\xi))\right)=0.$$ Using (4.43), (4.28) we get (4.52 ).
$\Box$ \vskip 0.1in \noindent We prove the convergence of the
integral (4.50). Since the marked points $\{y_i\}$ and the end
points ( puncture points ) $\{p_j\}$ vary in a compact set $K$ out
of the gluing part, the inequality (4.52) holds in $K$. By the
standard elliptic estimate we get $$\left|\frac{\partial}{\partial
r}\left(f_{(r)}\circ I'_r(\xi)\right)\right| \leq
\frac{C}{r^2}.\leqno(4.53)$$ By using a similar argument it is
easy to get $$\left|\frac{\partial}{\partial v}\left(f_{(r)}\circ
I'_r(\xi)\right) \right|\leq C\;\;\;\;\;\;\;
\left|\frac{\partial}{\partial v}I'_r(\xi)\right|\leq C \;\;\;
for\;\;other \;\; parameter \;\;v \leqno(4.54)$$ Let $\alpha\in
H^*(M^{+}, {\bf R})$ and $\beta\in H^*(Z,{\bf R})$ represented by
differential form. We may write $$\prod_i
e^*_i\alpha_i\wedge\prod_j e^*_j\beta_j\wedge I^*\Theta =
ydr\wedge dt \wedge d\xi\wedge dj,$$ where $d\xi$ and $dj$ denote
the volume forms of $ker D{\cal S}_{e,b}$ and the space of complex
structures respectively and $y$ is a function. Then (4.53) implies
that $|y|\leq \frac{C}{r^2}$. Then the convergence of the integral
(4.50) follows. \vskip 0.1in \noindent By the same argument as in
\cite{R5} one can easily show that \vskip 0.1in \noindent {\bf
Theorem 4.14 } \vskip 0.1in \noindent {\it (i).
$\Psi^{(M,Z)}_{(A,g,T_m)}(\alpha_1,..., \alpha_l; \beta_{l+1},...,
\beta_m)$ is well-defined, multi-linear and skew symmetric. \vskip
0.1in \noindent (ii).
$\Psi^{(M,Z)}_{(A,g,T_m)}(\alpha_1,...,\alpha_l; \beta_{l+1},...,
\beta_m)$ is independent of the choice of forms $\alpha_i,
\beta_j$ representing the cohomology classes $[\beta_j],
[\alpha_i]$,  and the choice of virtual neighborhoods. \vskip
0.1in \noindent (iii). $\Psi^{(M,Z)}_{(A,g,T_m)}(\alpha_1,...,
\alpha_l; \beta_{l+1},..., \beta_{m})$ is independent of the
choice of $\widetilde{J}$ and $J_{\nu}$.}

\section{\bf  A gluing formula}

We first prove an addition formula for operator $D_u$.
Let $u=(u^+,u^-):(\Sigma^+, \Sigma^-)\rightarrow (M^+,M^-)$ be
$J$-holomorphic curves such that $u^+$ and $u^-$ have $\nu$ ends and they
converge to the same periodic orbits at each end. Note that according
to our convention $\Sigma^{\pm}$ may not be connected (see the end of
Section 4). In this case $Ind(D_{u^{\pm}},\alpha)$ denotes the sum of
the indices of its components. Suppose that
$\Sigma =\Sigma^{+}\bigvee \Sigma^{-}$ has genus $g$ and $[u(\Sigma)] =A$.
Denote by $j$ ( resp. $j^\pm$ ) the complex structure on $\Sigma$ (resp.
$\Sigma^{\pm}$ ), and by $Ind(D_u,j)$ (resp. $Ind(D_{u^{\pm},j^{\pm}},
\alpha)$ ) the corresponding index. From the proof of Lemma 4.9
We get the following index addition formula of Bott-type:
$$Ind(D_{u^+,j^{+}},\alpha) + Ind(D_{u^-,j^{-}},\alpha)
-2(n+1)\nu =2C_1(A) + (n+1)(2-2g),\leqno(5.1)$$
where we used the fact that at every end $\dim \ker L_{\infty}=2(n+1).$
Considering the variation of the complex structures on the
Riemann surfaces we have from (3.63) that
\vskip 0.1in
\noindent
{\bf Theorem 5.1}
$$Ind(D_{u^+},\alpha)+Ind(D_{u^-},\alpha) =2n\nu + 2C_1(A) + (n+1)(2-2g)
+ 6g - 6.\leqno(5.2)$$
\vskip 0.1in
\noindent
{\bf Prof:}   Let us consider a simple case, the general case is identical.
Suppose that $\Sigma^+$ is a connected smooth Riemann surface of genus $g^+$, and
$\Sigma^-$ consists of two components:
$\Sigma^{-}_1$ of genus $g^{-}_1$ and $\Sigma^{-}_2$ of genus $g^{-}_2$. Suppose
that $\Sigma^+$ intersects  $\Sigma^{-}_1$ at $\nu_1$ points, and intersects
$\Sigma^{-}_2$ at $\nu_2$ points. $\Sigma^+$ and $\Sigma^-$ form a Riemann surface
$\Sigma$ of genus
$$ g= g^{+} + g^{-}_1 + g^{-}_2 + \nu_1 + \nu_2 - 2.$$
We add $6g-6$ to both sides of (3.63). Note that
$$Ind(D_{u^+},\alpha) = Ind(D_{u^{+}, j^{+}},\alpha) + 6g^{+} - 6 + 2(\nu_1
+ \nu_2),$$
$$Ind(D_{u^-},\alpha) = Ind(D_{u^{-}_1, j^{-}_1},\alpha) + 6g^{-}_1 - 6 + 2\nu_1 +
Ind(D_{u^{-}_2, j^{-}_2},\alpha) + 6g^{-}_2 - 6 + 2\nu_2.$$
Then (5.2) follows from the above equalities.$\Box $
\vskip 0.1in
\noindent
{\bf Remark 5.2 }{\it Let $u$ be a $J$-holomorphic map from
$(\stackrel{\circ}{\Sigma};y_{1},...,y_{m},
p_1,...,p_{\nu})$ into $M^{\pm}$ such that each end converges to a
periodic orbit. By using the removable singularities theorem we get a
$J$-holomorphic map $\bar{u}$ from $\Sigma $ into $\overline{M}^{\pm}$.
Therefore, we have a natural identification of finite energy pseudo-holomorphic maps
into $M^{\pm}$ and closed pseudo-holomorphic maps into the closed symplectic manifolds
$\overline{M}^{\pm}$. Moreover, the operator $D_u$ is identified with the operator
$D_{\bar{u}}$ in a natural way. Under this identification, the condition
that $u$ converges
to a $k$-multiple periodic orbit at a marked point $p$ is naturally interpreted as
$\bar{u}$ being tangent to $Z$ at $p$ with order $k$. Since  $\ker {L_{\infty}}$ consists of constant vectors,
we can identify the vector fields in ${\cal W}_{\pm}^{1,p,\alpha}$
along $u$ with the  vector fields in ${\cal W}_{\pm}^{1,p,\alpha}$ along
$\bar{u}$, the space $L_{\pm}^{p,\alpha}$ along $u$ is also identified with
the space $L_{\pm}^{p,\alpha}$ along $\bar{u}$.
In the case of closed manifolds the definitions of $L_{\pm}^{p,\alpha}$ and
${\cal W}_{\pm}^{1,p,\alpha}$ are the same as that of [Liu]).}

Thus we have
\vskip 0.1in
\noindent
{\bf Proposition 5.3 }
$$Ind(D_u,\alpha)=Ind D_{\bar{u}}.\leqno(5.3)$$
\vskip 0.1in
\noindent
We next prove a general gluing formula relating
GW-invariants of a closed symplectic manifold in terms of log
GW-invariants of its symplectic cut.

First, we make a remark about an error
in the draft concerning the homology class of the  general gluing formula. The gluing theorem in
the previous section shows that one can glue two pseudo-holomorphic curves
$(f_+, f_-)$ in $M^+, M^-$ with the same end point to a pseudo-holomorphic curve $f$ in
$M$. Suppose that the homology classes of $f_+, f_-, f$ are $A^+, A^-, A$. Then, we
carelessly wrote $A=A^+ + A^-$. R. Fintushel and E. Ionel pointed to us that in general
the homology class of $f$ depends on the pseudo-holomorphic curve representatives
$f_+, f_-$ instead of the homology classes $A^+, A^-$. One can also understand it as follows. Recall that there is a map
$$\pi: M\rightarrow \overline{M}^+\cup_Z \overline{M}^-.$$
$\pi$ induces a homomorphism
$$\pi_*:  H_2(M, \Z)\rightarrow H_2(\overline{M}^+\cup_Z \overline{M}^-, \Z).$$
Using the Mayer-Vietoris sequence for $(\overline{M}^+, \overline{M}^-,\overline{M}^+\cup_Z \overline{M}^-)$, $(f_+,f_-)$ defines a homology class
$[f^+ +f^-]\in H_2(\overline{M}^+\cup_Z \overline{M}^-, \Z)$. The existence of the glued map $f$ implies $[f^+ +f^-]=\pi_*([f])$.
If $(f'_+, f'_-)$ is another representative and glued to $f'$,
$$\pi_*([f'])=[f'_+ +f'_-]=[f_+ +f_-]=\pi_*([f]).$$
When $\ker \pi_*\neq 0$, $[f], [f']$ could be different from a
vanishing 2-cycle  in $\ker \pi_*$.  For the application of our gluing formula to the main theorems,
there is no vanishing cycle (Lemmas 2.11, 2.13, section 2). Hence,
this problem does not arise.

However, the original statement of our general gluing formula is incorrect.
Instead, our argument yields the following modified statement.
Let $[A]=A+\ker \pi_*$.
Then, we define
$$ \Psi_{([A],\cdots)}=\sum_{B\in [A]}\Psi_{(B, \cdots)}.$$
Note that for $B,B'\in [A]$, $\omega(B)=\omega(B')$. By the
  compactness theorem, there are only finitely many such $B$ to be
  represented by stable $J$-holomorphic maps. Hence the summation
  on the
  right hand side is finite.
  By abuse of the notation, we use $[A]=A^+ + A^-$ to represent the
  set of homology classes of glued maps.
  Then we  replace $\Psi_{(A, \cdots)}$ by $\Psi_{([A], \cdots)}$
  in all the statements in this section and the original proof is
  still valid.

The proof of our gluing formula is similar to the proof of the
composition law of GW-invariants and has two steps. The first step
is to define an invariant for $M_{\infty}$  and prove that it is
the same as the invariant of $M_r$. Then, we write the invariant
of $M_{\infty}$ in terms of log invariants of $M^{\pm}$.

  We first construct a virtual
neighborhood for $M_{\infty}$.
Suppose that $\Sigma^{\pm}$ has $l^{\pm}$ connected components
$\Sigma^{\pm}_i,i=1,...,l^{\pm}$ of genus $g^{\pm}_i$ with
$m^{\pm}_i$ marked points and ${\nu}^{\pm}_i$ ends. The moduli space
${\overline{\cal{M}}}_{[A]}(M_{\infty},g,m)$ consists of the
components indexed by the following data:
\begin{description}
\item[(1)] The combinatorial type of $(\Sigma^{\pm},u^{\pm})$:
$\{A^{\pm}_i,g^{\pm}_i,
m^{\pm}_i, (k^{\pm}_1,...,k^{\pm}_{l^{\pm}})\},i=1,...,l^{\pm};$
\item[(2)] A map $\rho:\{p^+_1,...,p^+_{\nu}\}\rightarrow
\{p^-_1,...,p^{-}_{\nu}\}$,
where $(p^{\pm}_1,...,p^{\pm}_{\nu})$ denote the puncture points of
$\Sigma^{\pm}$.
\end{description}
\vskip 0.1in
Suppose that ${\cal C}^{J,A}_{g,m}$ is the set of indices. Let $C \in
{\cal C}^{J,A}_{g,m}$. Denote by ${\cal M}_C$ the set of stable maps
corresponding to $C$. The following lemma is obvious.
\vskip 0.1in
\noindent
{\bf Lemma 5.4 }{\it ${\cal C}^{J,A}_{g,m}$ is a finite
set.} \vskip 0.1in
We use the same method as in the above
subsection to construct a virtual neighborhood $U_C$ for each component $C$
and get a virtual neighborhood $(U,E,I)$ for
${\overline{\cal{M}}}_{[A]}(M_{\infty},g,m)$  starting inductively
from the lowest stratum.
 But $U$ is usually not a
smooth manifold. To see this, we observe that the configuration
space $\overline{\B}_{[A]}(M_{\infty},g,m)$ can be identified as
$$\cup_{C} \overline{\B}_{C}(M_{\infty},g,m).$$
$\overline{\B}_{C}(M_{\infty},g,m)$,$\overline{\B}_{C'}(M_{\infty},g,m)$
may intersect each other at lower strata. Hence, $U=\cup_C U_C$,
where $U_C$ is a virtual neighborhood of $\M_C$. Note that
$\cal{M}_C$ and $\M_{C'}$ may intersect each other, where the
intersection corresponds to a stable map with some component in
$\R\times \widetilde{M}$. Then,  $U_C$ may intersect each other.
By our construction, $U$ has the same stratification structure as
that of $\overline{\cal{M}}_{[A]}(M_{\infty}, g, m)$. Hence,  for
$C\neq C^{\prime}$, $U_C\bigcap U_{C^{\prime}}$ is a stratum of
$U_C, U_{C'}$ codimension at least 2. \vskip 0.1in
The integration theory can be
extended to such a space in an obvious fashion. Namely, we  take
the sum of integrals over each $U_C$. Choose a Thom-form $\Theta$
of $E$.Then,  we can define GW-invariants
$\Psi_{(M_{\infty},[A],g,m)}$ using the same integral formula. We
can also define GW-invariants $\Psi_{C}$ for each component $C$.
By the same argument as in Section 4 we can show the convergence
of the integrals.
It is easy to see that $$\Psi_{(M_{\infty},[A],g,m)}=\sum_{{\cal
C}^{J,[A]}_{g,T_m}}\Psi_{C}.\leqno(5.4)$$ \vskip 0.1in \noindent
{\bf Remark 5.5 }{\it It is easy to see that \vskip 0.1in
\noindent {\bf (i)} For $C= \{A^+, g^+, m^+\}$, we have
$$\Psi_{C}(\alpha^+)=
\Psi^{(\overline{M}^+,Z)}_{(A^+,g^+,m^+)}(\alpha^+);
\leqno(5.5)$$ {\bf (ii)} For $C= \{A^-, g^-, m^-\}$, we have
$$\Psi_{C}(\alpha^-)=
\Psi^{(\overline{M}^-,Z)}_{(A^-,g^-,m^-)}(\alpha^-).
\leqno(5.6)$$} {\bf Theorem 5.6 }{\it For any r, $0<r<\infty$, we
have $$\Psi_{(M_{\infty},[A],g,m)} =
\Psi_{(M_{r},[A],g,m)}.\leqno(5.7)$$ }
{\bf Proof: }  For each $r$ we have $\overline{\cal{M}}_{A}(M_r,g,T_m)$,
denoted by $\overline{\cal{M}}_r$. We may construct a virtual neighborhood
$(U_r, E_r, I_r)$. Put
$\overline{\cal{M}}_{\{r\}}=\cup_r \overline{\cal{M}}_r\times \{r\}$,
Repeating the previous argument, we may construct a virtual neighborhood
$(U_{\{r\}}, E_{\{r\}}, S_{\{r\}})$ such that
$(U_1, E_1, I_1)$ is a virtual neighborhood at $r=1$
and $(U_{\infty},  E_{\infty},I_{\infty})$ is a virtual neighborhood
at $r=\infty$.In general
$U_{\infty}$ is not a smooth boundary of $U_{\{r\}}$. In order to
apply Stokes' Theorem, we need a clear description of
neighborhoods of each lower stratum of $U_{\infty}$ in $U_{(t)}$.
This is basically a gluing argument. We consider a simple case;
the argument for the general case is similar. Suppose that
$U_C\bigcap U_{C^{\prime}} = {\cal M}^{(+,0,-)}$ , where $${\cal
M}^{(+,0,-)}=\{ (a,b,d)\in {\cal M}^{+}_{{\cal S}_e}
(A^{+},1)\times{\cal M}^{c}_{{\cal S}_e}(A^{c},1)\times {\cal
M}^{-}_{{\cal S}_e}(A^{-},1)$$ $$| P^+(a)=P^-(b),
P^+(b)=P^-(d)\}.$$ Suppose that $(\Gamma^+,\Gamma_0,\Gamma^-)\in
{\cal M}^{(+,0,-)}$, where
$$\Gamma^{\pm}=\left((u^{\pm},\Sigma^{\pm},p^{\pm}),\eta^{\pm}\right),$$
$$\Gamma_0= \left((u_0,\Sigma_0,q^+,q^-), \eta_0 \right).$$ Let
$O_{(+,0.-)}$ be a neighborhood of $(\Gamma^+,\Gamma_0,\Gamma^-)$
in ${\cal M}^{(+,0,-)}$. For any $(\theta, r_1, r_2)$ we glue
$M^+, {\bf R}\times \widetilde{M}$ and $M^-$ to get
$M_{(\theta,r_1,r_2)}$ with the following gluing formulas
$$a^{+}=a_0 + 4r_1,\;\;\; a_0= a^- + 4r_2$$ $$\theta^{+}=\theta^-
= \theta_0 + \theta \;\;mod \;1.$$ Given a complex structure
$\xi=(\xi^+,\xi_0,\xi^-)$, we construct $\Sigma_{(\xi,r_1,r_2)}=
\Sigma^{+} \#_{r_1}\Sigma_0 \#_{r_2}\Sigma^{-} $ and
 a pre-gluing map $u_{(\theta,r_1,r_2)}:\Sigma_{(\xi,r_1,r_2)}
\rightarrow M_{(\theta,r_1,r_2)}$ in a similar way as in Section 4
(Recall that we perturb $u^{\pm}$ and $u_0$ only in the gluing
domain). Then we use a similar method as in Section 4 to prove
that a neighborhood of $(\Gamma^+,\Gamma_0,\Gamma^-)$ in $U_{(t)}$
can be described as $O_{(+,0.-)}\times N$, where $N:S^1\times
(R_1,\infty]\times (R_2,\infty]\rightarrow {\bf R}^4$ is a local
hypersurface given by $$y =
(e^{-r_1}\cos\theta,e^{-r_2}\sin\theta, e^{-r_2}\cos\theta,
e^{-r_2}\sin\theta ).$$ Both $$O_{(+,0.-)}\times N(S^1\times
\{\infty \}\times (R_2,\infty])$$ and $$O_{(+,0.-)}\times
N(S^1\times (R_1,\infty]\times \{\infty \})$$ are its boundary. We
show that the Stokes theorem still applies to such a space.
Without loss of generality we consider a domain $D$, $(0,0,0,0)\in
D\subset N$. The boundary $\partial D$ consists of three parts: \\
1. $\gamma_1$: the boundary of $D$ in the hypersurface $N$; 1.
$\gamma_1$: the boundary of $D$ in the hypersurface $N$; we assume
it to be smooth; \\ 2. $\gamma_2 = D \cap N(S^1\times \{\infty
\}\times (R_2,\infty])$;\\ 3. $\gamma_3 = D \cap N(S^1\times
(R_1,\infty]\times \{\infty \})$.\\ We draw a small ellipsoid
$B(\epsilon_1,\epsilon_2): \{\frac{x_1^2 + x_2^2}{\epsilon_1} +
\frac{x_3^2 + x_4^2}{\epsilon_2}=1\} $ around $(0,0,0,0)$ in ${\bf
R}^4$, then $D - B(\epsilon_1,\epsilon_2)$ is a polyhedron (union
of simplices) in Euclidean space, for which the Stokes Theorem
applies. We use Stokes' Theorem and then let $\epsilon_i
\rightarrow 0$. To simplify notations we don,t consider
$O_{(+,0.-)}$. For any differential form $\alpha$ we have
$$\int_{D}d\alpha= \int_{\partial D}\alpha -
\int_{N_{\epsilon}}\alpha $$ where  $N_{\epsilon}$ denotes $N\cap
\partial B(\epsilon_1,\epsilon_2)$. If we can show that the
integral is convergent as $\epsilon_i \rightarrow 0$, we can
obtain $$\int_{D}d\alpha= \int_{\partial D}\alpha.$$ We first let
$\epsilon_1 \rightarrow 0$, then $\epsilon_2 \rightarrow 0$, By
the estimate (4.53 ), (4.54 ) and the same argument as in the
subsection 4.2 the convergence follows. This argument obviously
works for $O_{(+,0.-)}\times N$. Then we can use Stokes' Theorem.
Theorem 5.6 is proved. $\Box$ \vskip 0.1in \noindent We derive a
gluing formula for the component $C=\{A^+, g^+, m^+, k; A^-, g^-,
m^-, k\}$. For any component $C$ we can use repeatedly this
formula. Choose a homology basis $\{\beta_b\}$ of $H^{\ast}(Z,
{\bf R})$. Let $(\delta_{ab} )$ be its intersection matrix. \vskip
0.1in \noindent {\bf Theorem 5.7 }{\it Let $\alpha^{\pm}_i$ be
differential forms with $deg \alpha^{+}_i=deg \alpha^{-}_i$ even.
Suppose that $\alpha^+_i|_Z=\alpha^-_i|_Z$ and hence
$\alpha^+_i\cup_Z \alpha^-_i\in H^*(\overline{M}^+ \cup_Z
\overline{M}^-, \R)$. Let $\alpha_i=\pi^*(\alpha^+_i\cup_Z
\alpha^-_i)$ (1.11). For $C=\{A^+, g^+, m^+, k; A^-, g^-, m^-,
k\}$, we have the gluing formula
$$\Psi_{C}(\alpha_1,...,\alpha_{m^++m^-}) =k\sum
\delta^{ab}\Psi^{(\overline{M}^+,Z)}_{(A^+,g^+,m^+,k)}
(\alpha^{+}_1,...,\alpha^{+}_{m^+}, \beta_a)
\Psi^{(\overline{M}^-,Z)}_{(A^-,g^-,m^-,k)}(\alpha^{-}_{m^++1},...,
\alpha^{-}_{m^++m^-},\beta_b).\leqno(5.8)$$ where we use
$\Psi^{(\overline{M}^{\pm},Z)}_{(A^{\pm},g^{\pm},m^{\pm},k)}$ to
denote
$\Psi^{(\overline{M}^{\pm},Z)}_{(A^{\pm},g^{\pm},T_{m^{\pm}})}$.}

\vskip 0.1in \noindent {\bf Proof: }  We denote $U_{C}=U|_{C}$,
$U_{(A^{\pm},k)}= (S^{\pm}_e)^{-1}(0)/stb_{u^{\pm}}$, and
$$U_{(A^+,A^-)}=\left\{(a,b)\in U_{(A^+,k)}\times U_{(A^-,k)} |
P^{+}(a)=P^{-}(b)\right\}.$$ There is a natural map of degree k
$$Q:U_{C}\rightarrow U_{(A^+,A^-)}$$ defined by
$$Q(a,b,n)=(a,b),$$ and a map $$P:U_{(A^+,k)}\times U_{(A^-,k)}
\rightarrow Z\times Z$$ defined by $$P(a,b)=(P^+(a),P^-(b)).$$
Note that $${\cal M}_{S_e}(A^{+},A^{-})=P^{-1}(\Delta),$$ where
$\Delta \subset Z\times Z $ is the diagonal. The Poincar\'e dual
$\Delta^{\ast}$ of $\Delta $ is $$\Delta^{\ast}=\Sigma
\delta^{ab}\beta_a \wedge \beta_b .$$ Choose a Thom form $\Theta=
\Theta^{+}\wedge \Theta^{-}$, where $\Theta^{\pm}$ are Thom forms
in $F^{\pm}$ supported in a neighborhood of the zero section . By
perturbing $\alpha^{\pm}_i$ we may assume that $\alpha_i$ are
smooth forms. Then
$$\Psi_{C}(\alpha_1,...,\alpha_{m^++m^-})=\int_{U_{C}}
\prod_{1}^{m^+}\alpha_i\wedge \prod_{m^+}^{m^++m^-}\alpha_j\wedge
I^{\ast}\Theta $$ $$= k
\int_{U_{(A^+,A^-)}}\prod_{1}^{m^+}\alpha_i\wedge
\prod_{m^+}^{m^++m^-}\alpha_j\wedge I^{\ast}\Theta $$
$$=k\int_{U_{(A^{+},k)}\times U_{(A^{-},k)}}
\sum\delta^{ab}\prod_{1}^{m^+}\alpha^{+}_i\wedge
I^{\ast}\Theta^{+}\wedge\beta_a \wedge
\prod_{m^+}^{m^++m^-}\alpha^{-}_j \wedge I^{\ast}\Theta^{-}
\wedge\beta_b $$ $$=k\sum
\delta^{ab}\Psi^{(\overline{M}^+,Z)}_{(A^+,g^+,m^+,k)}(\alpha^{+}_1,
...,\alpha^{+}_{m^+},\beta_a)
\Psi^{(\overline{M}^-,Z)}_{(A^-,g^-,m^-,k)}(\alpha^{-}_{m^++1},...,
\alpha^{-}_{m^++m^-},\beta_b). \;\;\;\;\Box $$ \vskip 0.1in
\noindent For general $C$ with $(k_1,...,k_{\nu})$-periodic orbits
we may easily obtain \vskip 0.1in \noindent {\bf Theorem 5.8 }{\it
$$\Psi_C(\alpha)=|{\bf k}|\sum_{I,J}\Psi^{(\overline{M}^+,Z)}_
{(A^+,g^+,m^+,{\bf k})}(\alpha^{+}, \beta_I)\delta^{I,J}
\Psi^{(\overline{M}^-,Z)}_{(A^-,g^-,m^-,{\bf
k})}(\alpha^{-},\beta_J),\leqno(5.9)$$ where we associate
$\beta_i\delta^{i,j}\beta_j$ to every periodic orbit as in Theorem
5.7, and put $|{\bf k}|=k_1...k_{\nu}$,
$\delta^{I,J}=\delta^{i_1,j_1}...\delta^{i_{\nu},j_{\nu}}$, and
denote by
$\Psi^{(\overline{M}^{\pm},Z)}_{(A^{\pm},g^{\pm},m^{\pm},{\bf
k})}(\alpha^{\pm},\beta_J)$ the product of log invariants
cooresponding to each component.} \vskip 0.1in \noindent For
example, for $C=\{A^+, g^+, m^+, k_1,k_2; ,A_1^-, g_1^-, m_1^-,
k_1, A_2^-, g_2^-, m_2^-, k_2\}$, our formula (5.9) reads:
$$\Psi_C(\alpha)=k_1k_2\sum_{i_1,i_2,j_1,j_2}\Psi^{(\overline{M}^+,Z)}_
{(A^+,g^+,m^+,k_1,k_2)}(\alpha^{+},
\beta_{i_1},\beta_{i_2})\delta^{i_1,j_1} \delta^{i_2,j_2}$$
$$\Psi^{(\overline{M}^-,Z)}_{(A_1^-, g_1^-, m_1^-,
k_1)}(\alpha_1^{-},\beta_{j_1}) \Psi^{(\overline{M}^-,Z)}_{(A_2^-,
g_2^-, m_2^-, k_2)}(\alpha_2^{-},\beta_{j_2}).$$ \vskip 0.1in
\noindent {\bf Remark 5.9: }{\it The definition of log invariants
and gluing formula extend to the case that $Z$ is a disjoint union
of smooth codimension two symplectic submanifolds in an obvious
fashion. Furthermore, if symplectic cutting only obtains one
symplectic manifolds, gluing formula also extends to this case in
an obvious fashion. } \vskip 0.1in

\section{\bf Proofs of the Main Theorems }

\noindent {\bf Proof of Theorem A} \vskip 0.1in Let $M$ be a
3-fold and $M_f$ be obtained by a flop. By Proposition 2.9, after
a perturbation of complex structures to almost complex structures,
$M, M_f$ have the same blow-up $M_b$. Without the loss of
generality, we can assume that $M_b$ is the blow up along a
$\O(-1)+\O(-1)$ rational curve $\Gamma \subset M$ or
$\Gamma_f\subset M_f$.    By proposition 2.10, we can relate $M,
M_b$ by a symplectic cutting such that $\overline{M}^{-}=M_b$, and
$$\overline{M}^{+} = P({\O}(-1) + {\O}(-1) + {\O}).$$ The same
procedure applies to $M_f$. We have
$$\overline{M}^{+}=\overline{M}^{+}_f= P({\O}(-1) + {\O}(-1) +
{\O}),\leqno(6.1)$$ $$\overline{M}^{-}=\overline{M}^{-}_f=
M_b.\leqno(6.2)$$ \vskip 0.1in Now, we compare GW-invariants of
$M, M_f$ with log GW-invariants of $M_b$. By Lemma 2.11, there is
no vanishing 2-cycle in our cases. Now we use the gluing formula
(5.8) to prove our assertion. First, we assume that
$deg(\alpha_i)\geq 4$ and is Poincar\'e dual to a point or a
2-dimensional homology class $\Sigma$. Then, we can choose the
pseudo-submanifold representative of $\Sigma$ such that $\Sigma$
is in $M^-$. Hence, $\alpha_i$ is supported in $M^-$. For
simplicity we consider a special component: $$C=\{A^+,g^+,m^+, k;
A^{-},g^{-},m^{-},k \}.$$ The general case can be treated in the
same way.

The following
argument depends only on an index calculation. The index has an additive
property. If a  relative stable map has more than
one component, we can always construct a pre-gluing map and use
the index of the pre-gluing map. Hence, we can assume that the stable
map under consideration has only two components $(u^+, u^-)$,
where
$u^{\pm}$ is a $J$-map in $M^{\pm}$.

 Since $M$ is a 3-fold, using the addition formula (5.2) for the index we have
$$Ind(D_{u^{-}},\alpha) + Ind(D_{u^{+}},\alpha)-4=Ind(D_u).\leqno(6.3)$$
By (5.3),
$$Ind(D_{u^+},\alpha)=Ind(D_{\bar{u}^+}).\leqno(6.4)$$
Note that $\overline{M}^{+} = P({\O}(-1) + {\O}(-1) + {\O})$.
Next we claim that
$$Ind(D_{\bar{u}^+})\geq 6.\leqno(6.5)$$
Note that $\bar{u}^+$ can be identified as a holomorphic curve $h$ in
$P({\O}(-1) + {\O}(-1) + {\O})$ over $\P^1$ by remark 3.24.
Then,
$$Ind(D_{h})=2(C_1([h])-k+1),\leqno(6.6)$$
where $C_1(M^+)$ is represented precisely by $3Z_{\infty}$, where $Z_{\infty}$
is the infinity divisor of $P(\O(-1)+\O(-1)+\O)$.
A simple index calculation shows that
$C_1([h])=3k$.  Hence, if $k>0$
$$Ind(D_{\bar{u}^+})\geq 6.$$
 From (6.3)-(6.5) we have
$$Ind(D_{u^{-}},\alpha) \leq -2 + Ind(D_u).\leqno(6.7)$$
Hence $Ind(D_{u^{-}},\alpha)< Ind (D_u)$. However, the representatives of $\alpha_i$ (hence $\phi^*(\alpha_i)$)
is supported in
$M^-$. It is clear that $\alpha^+_i=0$ and the only nontrivial terms in
the gluing formula are
$$k\sum \delta^{ab}\Psi^{(\bar{M}^+, Z)}_{(A^+,g^+,m^+,k)}(\beta_a)\Psi^{(\bar{M}^-, Z)}_{
(A^-,g^-,m^-, k)}(\{\alpha^-_i\}, \beta_b).
$$
However, $\sum_i deg (\alpha_i)=Ind (D_u) >Ind (D_{u^-})$. Then for any $\beta_b $,
$$\Psi^{(\overline{M}^-,Z)}_{(A^{-},g^{-},m^{-},k)}(\{\alpha^-_i\}, \beta_b)=0.$$
If $u^{\pm}$ has more than one end, say $\nu$ ends, $Ind(D_{u^-})$ increases faster than $4\nu$.
It will force $Ind(D_{u^-})$ to become even more negative.
It follows that $\Psi_{C}=0$ except for $C= \{A^+, g, m\}$ or
$C= \{A^-, g, m\}$ i.e., $C$ stays completely on the one side.

Suppose that $A\neq n[\Gamma]$. We claim that $C=\{A^-,g,m\}$. If
$C=\{A^+,g,m\}$, $A^+$ is represented by a stable map whose image
is completely inside the total space of $\O(-1)+\O(-1)$. Then, it must be homologous
to $n[\Gamma]$ for some $n$. Hence,  $\pi_*(A)=\pi_*(n[\Gamma])$, where $\pi_*$ is defined in (2.3).
However, there is no
vanishing 2-cycle. This contradicts  the assumption that $A\neq n
[\Gamma]$. Hence,
$$\Psi^M_{(A,g,m)}(\{\phi^*(\alpha_i)\})=\Psi^{(M_b, Z)}_{(A^-,g,m)}(\{\alpha^-_i\}).\leqno(6.8)$$
It is easy to observe that $A^-=\phi_b(A)$, where $\phi_b: H_2(M, Z)\rightarrow H_2(M_b, \Z)$ is defined
in (2.14).
The same argument shows that
$$\Psi^{M_f}_{(\phi(A),g,m)}(\{\alpha_i\})=\Psi^{(M_b, Z)}_{(\phi_b(A),g,m)}(\{\alpha^-_i\})\leqno(6.9).$$
Hence,
$$\Psi^M_{(A,g,m)}(\{\phi^*(\alpha_i)\})=\Psi^{M_f}_{(\phi(A),g,m)}(\{\alpha_i\})\leqno(6.10).$$

When $deg(\alpha_i)=2$, we can use (1.1) to reduce it to the previous
case.

When $A=n[\Gamma]$, we do not need the gluing formula. By the assumption, $\Gamma$
generates an extremel ray. If $f$ is a stable map representing
$A$, $im(f)=\Gamma$. $\Gamma, \Gamma_f$ have isomorphic
neighborhoods as complex manifolds. It follows from the definition
that
$$\Psi^M_{(n[\Gamma], g)}=\Psi^{M_f}_{(n[\Gamma_f],g)}. \Box $$
\vskip 0.1in
{\bf Remark 8.1: }{\it Although it is not needed in our proof,
a similar gluing argument can show that
$$\Psi^{(M_b, Z)}_{(\phi_b(A),g)}(\{\alpha^-_i\})=\Psi^{M_b}_{(\phi_b(A),g)}(\{\alpha^-_i\})\leqno(6.11).$$
}
\vskip 0.1in
Only Corollary A.2 needs a proof. The others are immediate
consequences of Theorem A and Corollary A.2.
\vskip 0.1in
\noindent
{\bf Proof of Corollary A.2}
\vskip 0.1in
 For simplicity, we assume that all the exceptional curves are in the same
homology class $[\Gamma]$. Let $n_{\Gamma}=n_{\Gamma_f}$ be the number of
such curves, $[\Gamma]$ and $[\Gamma_f]$ be the homology classes of the
exceptional curves respectively in $M$ and $M_f$. Then for any $0<k\in {\bf Z}$,
by the formula for multiple cover maps, $\Phi^{M}_{(k[\Gamma],0)} =
\Phi^{M_f}_{(k[\Gamma_f],0)} = \frac{n_{\Gamma}}{k^3}$ (see \cite{V}).
 Then the total 3-point function can be written in
the form
$$\Psi^{M}(\beta_1,\beta_2,\beta_3)
=\beta_1\wedge \beta_2\wedge\beta_3 + \sum_{ A\neq k[\Gamma]}
\Psi^{M}_{(A,0,3)}(\beta_1, \beta_2, \beta_3)q^{A}$$
$$+\sum_{n[\Gamma], n\neq 0}\Psi^{M}_{(A,0,3)}(\beta_1, \beta_2, \beta_3)q^{A}.\leqno(6.12)  $$
As for the last term, it is zero except that $deg (\beta_i)=2$ for
the dimension reasons. If $deg (\beta_i)=2$,
$$\sum_{n[\Gamma], n\neq 0}\Psi^{M}_{(A,0,3)}(\beta_1, \beta_2, \beta_3)q^{n[\Gamma]}
= \frac{q^{[\Gamma]}}{1- q^{[\Gamma]}}
\beta_1([\Gamma])\beta_2([\Gamma])\beta_3([\Gamma])n_{\Gamma}.\leqno(6.13)$$
As for the first term, it is also zero except that $(deg (\beta_1),
deg (\beta_2), deg (\beta_3))=(2,2,2), ( 4,2,0), (6,0,0).$
There is a similar expression for the 3-point function of $M_f$. We see
that, by Theorem A, only the first and last  terms are different under identification $q^A\rightarrow q^{\phi(A)}$. Suppose that
$\beta_i=\varphi^*\alpha_i$. Let's first consider the case that
$(deg (\beta_1), deg (\beta_2), deg (\beta_3))=( 4,2,0), (6,0,0).$ In these cases, the last terms are zero. We claim that
the first terms are the same as well. This is obvious for the case of $(6,0,0)$. For the case $(4,2,0)$, we
can assume that $\alpha_3=1$.  By  definition,
the map $H^4(M_f, \R)\rightarrow H^4(M, \R)$ is Poincar\'e dual to the inverse of the map
$H_2(M, \Z)\rightarrow H_2(M_f, \Z)$. By our construction, $\alpha_1, \beta_1$ are Poincar\'e
dual to the same 2-manifold $\Sigma$ disjoint from the exceptional locus. Hence,
$$\varphi^*\alpha_2\wedge \varphi^*\alpha_1=\varphi^*\alpha_2(\Sigma)=\alpha_2(\Sigma)=\alpha_2\wedge \alpha_1.$$
For the case $(2,2,2)$, both the first term and last term are non-zero. Let $p: M_b\rightarrow M;
p_f: M_b\rightarrow M_f$ be the projection. We observe that $p^*_f\alpha_i-p^*\beta_i=
-\alpha_i(\Gamma_f)Z$, where $Z$ is the exceptional divisor.
A routine calculation shows that
$$\alpha_1\wedge\alpha_2\wedge\alpha_3-\beta_1\wedge\beta_2\wedge\beta_3=\alpha_1(\Gamma_f)\alpha_2(\Gamma_f)\alpha_3(\Gamma_f).\leqno(6.14)$$
For $q^{[\Gamma_f]}\neq 1$,
$$\frac{q^{[\Gamma_f]}}{1-q^{[\Gamma_f]}}+\frac{q^{-[\Gamma_f]}}{1-q^{-[\Gamma_f]}}=-1.\leqno(6.15)$$
Recall that
$$\beta_1(\Gamma)\beta_2(\Gamma)\beta_3(\Gamma)=-\alpha_1(\Gamma_f)\alpha_2(\Gamma_f)\alpha_3(\Gamma_f)\leqno(6.16)$$
and
$$\phi([\Gamma])=-[\Gamma_f].\leqno(6.17)$$
We have
$$\alpha_1\wedge\alpha_2\wedge\alpha_3 + \frac{q^{[\Gamma_f]}}{1- q^{[\Gamma_f]}}
\alpha_1(\Gamma_f)\alpha_2(\Gamma_f)\alpha_3(\Gamma_f)n_{\Gamma_f}
= \beta_1 \wedge\beta_2 \wedge\beta_3 + \frac{q^{\phi([\Gamma])}}
{1- q^{\phi([\Gamma_f])}}\beta_1(\Gamma)\beta_2 (\Gamma)
\beta_3 (\Gamma)n_{\Gamma}.\leqno(6.18)$$
The assertion follows. $\Box$
\vskip 0.1in
\noindent
{\bf Proof of Theorem B}
\vskip 0.1in
The proof is similar to that of Theorem A. By Proposition 2.12, we can perform  a symplectic cutting
such that one
part $\overline{M_e}^{-}$ is $M_b $. The other part $\overline{M_e}^{+} $ is a collection of quadric 3-folds.
Without the loss of generality, we can assume that $\overline{M_e}^+$ has only one
component.
 Again, we first argue the case that
$deg(\alpha_i)\geq 4$, where we choose the support of $\alpha_i$ inside $M^-$. By a simple index calculation, if $h$ is
a holomorphic curve of $Y$ tangent to the infinity divisor with  order $k$, its
index is $2(3k-k+1)\geq 6$ if $k>0$.
By the same argument as in Theorem A we conclude that we need only consider
those $J$-holomorphic curves which don't go through $\widetilde{M}$.
In this case, there is no holomorphic curves in $\overline{M_e}^+$ disjoint from $Z$.
This shows that for $B\neq 0$
$$\Psi^{M_e}_{B,g,m}(\{\alpha_i\})=\sum_{A^-}\Psi^{(M_b, Z)}_{(A^-,g,m)}(\{\alpha^-_i\}),\leqno(6.19)$$
where summation is taken over $A^-$ homologous to $\pi_*(A)$. A
moment of thought tells us that these are the set of $\phi_b(A)$
such that $\phi_e(A)=B$.
Then, the Theorem B follows from (6.19) and (6.8)
Then, we use formula (1.1) to reduce the case $deg(\alpha_i)=2$ to the previous case.
$\Box$

Only Corollary B.2 needs a proof. The others are immediate
consequences of Theorem B and Corollary B.2.

\vskip 0.1in
\noindent
{\bf Proof of Corollary B.2}
\vskip 0.1in

    The proof is a generalization of  Tian's argument \cite{Ti}. The surjective map
$$\varphi: H_2(M, \R)\rightarrow H_2(M_e, \R) \leqno(6.20)$$
induces an injective map
$$\varphi^*: H^2(M_e, \R)\rightarrow H^2(M, \R).\leqno(6.21)$$
By definition, the map on $H^4$ is Poincar\'e dual to  a right inverse of (6.20).
We claim that the ordinary cup product remains  the same after transition. Assume that $\beta_i=\varphi^*\alpha_i$.
Again, we need to consider the
case that $(deg (\beta_1),
deg (\beta_2), deg (\beta_3))=(2,2,2), ( 4,2,0), (6,0,0).$
The case $(6,0,0)$ is obvious. The proof of the case $(4,2,0)$ is similar to that of Corollary of Theorem A.
$\beta_1$ is Poincar\'e dual to $A_1\in H_2(M, \R)$ such that $\alpha_1$ is Poincar\'e dual to $\varphi_*(A_1)$.
Hence,
$$\varphi^*\alpha_2\wedge \varphi^*\alpha_1=\varphi^*\alpha_2(A_1)=
\alpha_2(\varphi_*(A_1))=\alpha_2\wedge \alpha_1.\leqno(6.22)$$
For the case $(2,2,2)$, clearly $\varphi^*(\beta_i)(\Gamma)=0$. Without
 loss of generality, assume that there is a 4-manifold representing
$\varphi^*(\beta_i)$ which is disjoint from $\Gamma$. Hence, it can be viewed as
a submanifold of $M_e$. Clearly, the same 4-manifold represents $\beta_i$.
Hence,
$$\beta_1\wedge \beta_2 \wedge \beta_2=\alpha_1\wedge\alpha_2\wedge\alpha_3.\leqno(6.23)$$
Therefore,
$$\Psi^{M}_{\varphi^*w}(\varphi^*(\beta_1),\varphi^*(\beta_2),\varphi^*(\beta_3))
=\varphi^*(\beta_1)\wedge \varphi^*(\beta_2)\wedge\varphi^*(\beta_3)\leqno(6.24)$$
$$+ \sum_{ A\neq k[\Gamma]}\sum_{m}\frac{1}{m!}
\Psi^M_{(A,0, m+3)}(\varphi^*(\beta_1),\varphi^*(\beta_2),\varphi^*(\beta_3),
\varphi^*w, \cdots, \varphi^*w)q^A$$
$$+\sum_{ k[\Gamma], k\neq 0}\sum_{m}\frac{1}{m!}
\Psi^M_{(A,0, m+3)}(\varphi^*(\beta_1),\varphi^*(\beta_2),\varphi^*(\beta_3), w,\cdots, w)q^{k[\Gamma]}.$$
By the previous argument, the first term is the same as $\alpha_1\wedge\alpha_2\wedge\alpha_3$.
The last term is always zero.
Now we change the formal variable by $q^A\rightarrow q^{\varphi_*(A)}$ and apply Theorem B.
Then we prove
$$\Psi^{M}_{\varphi^*w}(\varphi^*(\alpha_1),\varphi^*(\alpha_2),\varphi^*(\alpha_3))=
\Psi^{M_e}_w(\alpha_1, \alpha_2, \alpha_3).\leqno(6.25)$$

\end{document}